\documentclass[reqno]{amsart}
\usepackage[utf8]{inputenc}
\usepackage{enumitem}

\usepackage{amsmath,amssymb,amsbsy,amsfonts,amsthm,latexsym,amsopn,amstext,mathtools,leftidx,
            amsxtra,euscript,amscd,verbatim,epsf,colordvi,graphics}
\usepackage{color}
\usepackage{pgf,tikz}
\usetikzlibrary{decorations.pathreplacing}
\usetikzlibrary{arrows}
\usepackage{mathrsfs}
\usepackage{graphicx}

\usepackage{hyperref}
\usepackage{cleveref}
\usepackage[margin=11pt,font=footnotesize,labelfont=bf]{caption}
\newtheoremstyle{tplain}{3pt}{3pt}{\rmfamily}{}{\bfseries}{.}{0.5em}{}
\theoremstyle{tplain}

\usepackage{pgf,tikz}
\usepackage{tcolorbox}
\usepackage[enableskew]{youngtab}
\definecolor{darkgreen}{cmyk}{1,0,1,0}
\newtheorem{thm}{Theorem}
\newtheorem{lem}{Lemma}
\newtheorem{ex}{Example}
\newtheorem{cor}{Corollary}
\newtheorem{prop}{Proposition}
\newtheorem{obs}{Remark}
\newtheorem{defi}{Definition}
\def \cal{\mathcal}
\def \rm {\mathrm}
\def \mbf {\mathbf}
\def \std {\mathrm{std}}
\newcommand{\oo}{\color{blue}}

\makeatletter
\newcommand*\bigcdot{\mathpalette\bigcdot@{.5}}
\newcommand*\bigcdot@[2]{\mathbin{\vcenter{\hbox{\scalebox{#2}{$\m@th#1\bullet$}}}}}
\makeatother

\newcommand{\YT}[3]{
\vcenter{\hbox{
\begin{tikzpicture}[x={(0in,-#1)},y={(#1,0in)}] 
\foreach \rowi [count=\i] in {#3} {
 \foreach \e [count=\j] in \rowi {
  \draw (\i,\j) rectangle +(-1,-1);
  \draw (\i-0.5,\j-0.5) node {$#2\e$};
 }
}
\end{tikzpicture}
}}
}

\title[Sagan-Stanley skew RSK and ballot switching] { Skew RSK  and the switching on ballot tableau pairs}

\author{Olga Azenhas}
\address{ University of Coimbra, CMUC, Department of Mathematics, Portugal}
\email{oazenhas@mat.uc.pt}



\begin{document}

\begin{abstract}
   In  arXiv:1808.06095 we have  introduced
 the  Knuth class of the word recording a sequence of locations for repeated internal insertion operations in the Sagan-Stanley skew RSK correspondence,  with no prescribed external insertion of new cells, to be  a preserver for the $P$-tableau.
  As a consequence  the Benkart-Sottile-Stroomer switching involution on ballot tableau pairs allows a realization as a recursive internal  insertion procedure.
  This amounts to explain the various  presentations of  Littlewood-Richardson (LR) commuters  and their coincidence predicted by Pak and Vallejo  with contributions by Danilov and Koshevoi. In particular,  the aforesaid presentation provides internal insertion as an alternative to Sch\"utzenberger- Lusztig involution (or evacuation) to constructing the Gelfand-Tsetlin pair in the Henriques-Kamnitzer $\mathfrak{gl}_n$-crystal commuter.   In addition,  the  coincidence of LR commuters solves
  the   Lecouvey-Lenart conjecture, recently further developed by Kumar-Torres,
  on  bijections between the Kwon and Sundaram
  branching models.
\end{abstract}
\maketitle

\tableofcontents


\section{Introduction}

We give a recursive presentation  of the Benkart-Sottile-Stroomer (BSS) switching map \cite{bss96} on ballot tableau (also known as Littlewood-Richardson tableau) pairs based on the Sagan-Stanley internal row  insertion procedure  on skew-tableaux \cite{ss} (see also \cite{rssw}, and \cite{imusa} for recent developments). Our key tool is to observe that the Knuth class of an internal  insertion order word of locations in a skew-tableau  preserves the $P$-tableau in the Sagan-Stanley skew Robinson-Schensted-Knuth (RSK) correspondence    with no prescribed external insertion of new cells.
 Our aim is to provide a clarification of, and, thereby, fulfill  the original question raised by Pak and Vallejo in \cite{pakvallejo}, with further contributions by Danilov and Koshevoi \cite{DK05,DK08}, on the coincidence and the involutive nature of all Littlewood-Richardson (LR) commuters.  Danilov and Koshevoi \cite{DK08} have proven the coincidence of the BSS switching commuter, referred to as $\rho_1$ in \cite{pakvallejo}, with the  Henriques-Kamnitzer crystal commuter in type $A$ and hive commuter,  here denoted $Com_{HK}$ respectively  $Com^h_{HK}$, \cite{HK2,HenKam}.

  LR commuters are combinatorial constructions that express in a bijective manner the symmetry of
 the $\mathfrak{gl}_n$-tensor products. It means to specify combinatorial objects exhibiting the symmetry of the tensorands $V_\mu$ and $V_\nu$ in $V_\mu\otimes V_\nu\simeq V_\nu\otimes V_\mu $, by counting  the multiplicity $c^\lambda _{\mu,\nu}=c_{\nu,\mu}^\lambda$ of the irreducible $\mathfrak{gl}_n$-module $V_\lambda$ in the decomposition into irreducibles of the tensor
product $V_\mu \otimes V_\nu$.
An alternative interpretation of these numbers is that they are the structure constants of the cohomology ring of a Grassmannian in the basis
of Schubert classes. More generally, beyond Cartan type $A$, maps exhibiting the isomorphism of the tensor product of $\mathfrak{g}$-crystals $A\otimes B$ and $B\otimes A$ are called \emph{commuters}.

 Here, we focus on the LR commuter or fundamental symmetry map referred to as $\rho_3$ in \cite[Section 6.1]{pakvallejo} which corresponds to the commutativity bijections originally on ballot tableaux \cite{az1,az2}, and later detailed on ballot tableaux and hives in \cite{akt16, tka18},
  and referred as $\rho^{(n)}$ respectively $\sigma^{(n)}$.
  More precisely, our $\rho_3$ commuter, thanks to the the BSS switching $\rho_1$ presentation on ballot tableau pairs as a recursive Sagan-Stanley internal row insertion, coincides with the Henriques-Kamnitzer $\mathfrak{gl}_n$-crystal commuter \cite{HK2,HenKam}  by providing the associated Gelfand-Tsetlin pair.

    Henriques and Kamnitzer \cite{HK2}, following an idea of A. Berenstein, show that the map $ a\otimes b\mapsto\xi(\xi(b)\otimes\xi(a))$, with $\xi$ the Lusztig-Sch\"utzenberger involution,
does give a crystal isomorphism from $B(\mu)\otimes B(\nu)$ to $B(\nu)\otimes B(\mu)$ and thereby a crystal commuter. In other words, denoting by $LR(\lambda/\mu,\nu)$ the set of ballot tableaux of shape $\lambda/\mu$ and content $\mu$, to conclude that a commuter $\rho$ coincides with the  Henriques-Kamnitzer $\mathfrak{gl}_n$-crystal commuter  \cite{HenKam,HK2}, it is enough, from considerations on highest and lowest weights in a crystal, to show that the commuter $\rho:LR(\lambda/\mu,\nu)\longrightarrow LR(\lambda/\nu,\mu)$  is such that for $T\in LR(\lambda/\mu,\nu)$ with left and right Gelfand-Tsetlin (GT) pair $(G_\mu,G_\nu)$, $\rho(T)$ has left and right  Gelfand-Tsetlin (GT) pattern pair given by $(\xi(G_\nu),\xi(G_\mu))$. The commuter $\rho_3$ does so but replaces the Sch\"utzenberger involution with Sagan-Stanley internal row insertion or its reverse internal insertion. 
  At this point it must be noted that any of  those two  Gelfand-Tsetlin patterns together with  a triple of boundary partitions completely specify an LR tableau and the corresponding LR hive \cite{fultonbuch,akt16,krv21}.

 The coincidence of Littlewood-Richardson (LR) commuters is instrumental   on a  Lecouvey-Lenart conjecture \cite{leclen}, recently further developed by Kumar-Torres \cite{sathishtorres,sathishkumarb}, on  bijections between the Kwon \cite{kwon18} and Sundaram  \cite{sundaram} branching models. As mentioned by Kumar-Torres in \cite{sathishtorres}, the only difference
between their bijection and the bijection conjectured by Lecouvey-Lenart is the
Littlewood–Richardson commuter used. While Lecouvey-Lenart use
Henriques–Kamnitzer $\mathfrak{gl}_n$-crystal commuter \cite{HenKam, HK2}, Kummar-Torres use the one by Kushwaha–Raghavan–Viswanath \cite{krv21} on flagged hives studied in \cite{krv21,krsv24}. The Lecouvey-Lenart conjecture is then positively answered thanks to the coincidence of LR commuters. We also note that a major fact in the settling of this conjecture is  that  Kumar-Torres bijection  restricts to  tableaux satisfying the Sundaram condition and those whose evacuation  satisfy the Kwon condition by considering and recognizing that they can be embedded in  the 
Kushwaha–Raghavan–Viswanath \cite{krv21} bijection on flagged hives. In other words, denoting by ${}^-LR^\lambda_{\mu,\nu}$ the set of left companions of $LR(\lambda/\mu,\nu)$, the Kumar-Torres bijection shows that the left  companion of a Sundaram LR tableau is a tableau in ${}^-LR^\lambda_{\mu,\nu}$ satisfying the Kwon condition. 
See Section \ref{conjecture}.

 Other  realisations for the  Benkart-Sottile-Stroomer (BSS) switching  commuter on ballot tableau pairs  (that is, a tableau pair $(Y,T)$, written $Y\cup T$, with $Y$ the Yamanouchi tableau of shape $\mu$ and $T$ a ballot tableau of shape $\lambda/\mu$ with $\mu\subseteq \lambda$), denoted $\rho_1$ in
 \cite{pakvallejo}, are based on   compositions of Sch\"utzenberger involutions \cite[Section 3]{lee}, \cite{pakvallejo}, or on tableau sliding, as in the Thomas-Yong infusion involution \cite{tyong2,tyong1,infusion}. The latter is realised in \cite{tyong2} via Fomin's {\em jeu de taquin} growths \cite[Chapter 7, Appendix 1]{stanley}.
   Beyond type $A$, Lenart \cite{Le}   realises    the Henriques-Kamnitzer  $\mathfrak{g}$-crystal commuter  \cite{HK2,HenKam} via van Leeuwen's {\em jeu taquin} \cite{lee} generalising  the Fomin's growth diagram presentation of {\em jeu de taquin} on Young tableaux { for a broader set of root systems beyond the type A.
    For further details on coincidence of various LR commuters, $\rho_1, \rho_2, \rho'_2=\rho_2^{-1}$
in \cite{pakvallejo},
the Henriques--Kamnitzer $\mathfrak{gl}_n$--crystal commuter  \cite{HK2, HenKam} and hive commuter \cite{HenKam}, as well
 the commutativity bijection of Danilov and Koshevoi \cite{DK05, DK08} on arrays, the Knutson-Tao-Woodward puzzles \cite{knutson} and the mosaic model  \cite{mosaic}, we refer the reader to \cite[Section 12]{akt16}, \cite[Introduction]{tka18} and \cite{azkoma25}.

 \subsection{Sagan-Stanley internal insertion and our results} On skew tableaux there are two types of insertion  \cite{ss}: external and internal  both of which are based on the usual Schensted insertion operation. However the corresponding procedures on a  skew tableau $T$ are different. The former proceeds very similarly to the usual Schensted insertion. The later has two main steps, firstly one chooses an  \emph{inner corner}  of $T$ (see Section \ref{sec:innercorner} for the definition) and  bumps its entry, and, secondly, one inserts the bumped entry externally, in the  row immediately  below, in the usual manner. Eventually the bumping route lands at the end of some row of $T$ where the last bumped entry settles and thus added at the end of  that row of $T$.   The \emph{internal row insertion operation}, denoted $\phi_i$ if the row coordinate of the vacated inner corner is $i$ (see Definition \ref{def:operator}), is two-fold, adds one box, the vacant box, to the inner shape which in turn expands the outer shape in one box, the last bumped entry, but, contrary to the external insertion, without contributing with a new element to the multiset of entries of $T$.

 The internal insertion procedure on a skew-tableau is an iteration of  the internal row insertion operation and thus requires {\em a priori} in each iterative step an inner corner of that skew tableau. Such information is encoded by a second skew-tableau, sharing the inner border with the first,  in the Sagan-Stanley skew-RSK correspondence  \cite{ss}. On its turn the instructions that it provides can be translated into a { \em companion word}, Definition \ref{comp},  listing the row coordinates of the entries, in the standard order, of the second skew-tableau. This word is {\it the internal row insertion order word} of the first tableau. The internal insertion procedure  is not  independent of the order of the chosen inner corners. However the Knuth class of the \emph{ companion word} of the second skew-tableau provides  a set of internal insertion order words preserving  the $P$-tableau in the Sagan-Stanley skew RSK correspondence,  when the matrix  prescribing external insertion of new cells is empty \cite{ss,rssw}, as shown in Theorem \ref{1}.

Knuth relations on the companion word of the second skew tableau is a partial contribution to the question under what conditions is the $P$-tableau preserved in the skew RSK (see question $(3)$ of \cite[ Section 9]{ss}).
 The internal insertion procedure  in general is not independent of a particular sequence of chosen inner corners. Recently Imamura-Mucciconi-Sasamoto \cite{imusa} observed the same property for the invariance of $P$-tableau  under Knuth relations on the companion word of the second skew tableau. Although  Knuth relations do not capture completely the invariance of the $P$-tableau in the Sagan-Stanley skew RSK  correspondence, (see Example \ref{ex:sufficientcond}), they are enough for the purpose of our paper.
 A nice observation \cite{rssw} is that the rectification of  a skew-tableau of inner shape, say $\mu$, can  either be calculated by using {\em jeu de taquin} or the internal insertion procedure by choosing, in the Sagan-Stanley internal skew procedure, an arbitrary  second skew-tableau of inner shape $\mu$, and \emph{outer} shape an appropriate rectangle. That is,  the rectification  does not depend on the order of {\em jeu de taquin moves} nor on the  internal insertion order words provided by the mentioned rectangular skew tableaux.
   In fact, it turns out that  companion words of rectangle tableaux are { anti-Yamanouchi words}, therefore, Knuth equivalent when of the same content.

Two words $\pi$ and $\pi'$ are Knuth equivalent if and only if their
$P$-tableaux under RSK correspondence are equal $P(\pi) = P(\pi')$ \cite{stanley}.    Theorem \ref{1} (Theorem \ref{th:U}) below is a natural generalization of this property for the skew RSK  rephrased in Theorem \ref{th:internal} for the internal insertion location words in the Sagan-Stanley internal insertion correspondence.

The \emph{companion word} $\mathcal{R}(U)$ of a skew tableau $U$ defines an internal row insertion operator $\phi_{\mathcal{R}({U})}$ for any skew tableau $T$ with the same inner shape as $U$. See  Definition \ref{comp} respectively  Definition \ref{def:operator} and its extensions \eqref{def:operatorbar}, \eqref{internalinsertionwordoperator}. If $U$ is a skew tableau, with inner shape $\mu=(\mu_1,\dots,\mu_n)$, on the alphabet $[n]$, the companion word of $U$ factorizes into $n$ maximal row words (possibly with some empty factors) $\cal{R}(U)=\cal{R}_n\cdots \cal{R}_1$ where $\cal R_i$ is the  row word defined by the row coordinates of the $i$-cells in $U$ for $i=1,\dots,n$, and the internal insertion operator $\phi_{\mathcal{R}({U})}$ factorizes accordingly $\phi_{\mathcal{R}({U})}=\phi_{\mathcal{R}_n}\circ\cdots \circ\phi_{\mathcal{R}_1}$. We extend the action of $\phi_{\mathcal{R}({U})}$ on $T$ to the pair $Y\cup T$ with $Y=Y_{\mu}$ the Yamanouchi tableau of shape $\mu$, and denote it $\bar \phi_{\mathcal{R}({U})}$, by  filling the vacated cell of $T$ under the action of $\phi_i$ with $i$:  $\bar\phi_i(Y\cup T)= Y_{(\mu_1,\dots,\mu_i+1,\mu_{i+1},\dots, \mu_n)}\cup \phi_i(T)$.
 We refer the reader to Sections \ref{sec:tabasic} and  \ref{sec:innercorner} for an explanation of
undefined terms.

Let $YT(\lambda/\mu)$ be the set of all semistandard tableaux (SSYT) of shape $\lambda/\mu$.
\begin{thm}\label{1} [Theorem \ref{th:U},  Sagan-Stanley row internal insertion operators and Knuth relations]. Let $T\in YT(\alpha/\mu)$,  $U, U'\in YT(\beta/\mu)$ and $P(T,U)$ respectively $P(T,U')$ the  corresponding $P$-tableaux in the Sagan-Stanley internal insertion correspondence. Then

$(a)$ $U$ and its standardization $\mathrm{std}\, U$ have the same companion word, ${\mathcal{R}({U})}=\mathcal{R}(\mathrm{std}\, U)$.

$(b)$ $P(T,U)=\phi_{\mathcal{R}(U)}(T)$.

$(c)$ $P(T,U)=P(T,\mathrm{std}\,U)$
 and  internal row insertion commutes with standardization  
$\mathrm{std}(P(T,U))=$ $P(\mathrm{std}\, T, U)$ $=P(\mathrm{std}\, T,\mathrm{std}\, U).$

$(d)$ $ P(T,U)  =P(T,U')=\phi_{\mathcal{R}(U)}(T)=\phi_{\mathcal{R}(U')}(T)$ whenever $\mathcal{R}(U)\equiv \mathcal{R}(U')$  are Knuth equivalent.

\end{thm}
Finally, Knuth equivalence of internal row insertion order words in Theorem \ref{1}, $ (d)$, means  (see Proposition \ref{propp:knuth})  that the composition  of internal row insertion operators  satisfy Knuth relations.
The  Knuth relations  satisfied by the Sagan-Stanley internal row insertion operators are the key fact  in Theorem \ref{2} (Theorem \ref{th:main})  to show  that Benkart-Sottile-Stroomer switching map \cite{bss96} on ballot tableau pairs,  denoted $\rho_1$,  can  be rephrased  in  the language of Sagan-Stanley internal row insertion operations or, equivalently, reverse internal insertion operations as next theorem { states}.

Let $n\ge 1$ and as usual put $[n]=\{1,\dots,n\}$. The set  ${\cal{LR}}^{(n)}$ denotes the set of all ballot semistandard tableau pairs $Y_\mu\cup T$,  say  $\mu=(\mu_1,\dots,\mu_n)$, and $T$ a ballot tableau of skew shape $\lambda/\mu$, for some $\mu\subseteq \lambda=(\lambda_1,\dots,\lambda_n)$. The switching map on ${\cal{LR}}^{(n)}$ is denoted by $\rho_1^{(n)}$. For $1\le i\le n$, let $(Y_\mu\cup T)^{[i]}:=Y_{(\mu_1,\dots,\mu_i)}\cup T^{[i]}\in {\cal{LR}}^{(i)}$  be the  restriction of $Y_\mu\cup T$ to the first $i$ rows with $T^{[i]}$ of  shape $(\lambda_1,\dots, \lambda_i)/(\mu_1,\dots,\mu_i)$.
(See Section \ref{sec:companionwords} and Definition \ref{def:restrictiontab} for precise definitions.)

\begin{thm} \label{2} [Main Theorem \ref{th:main} ]   Let $n\ge 1$ and $Y_\mu\cup T\in {\cal{LR}}^{(n)}$ with  $T$ a ballot tableau of shape $\lambda/\mu$ and weight $\nu$.
For $1\le i\le n$, let $(Y_\mu\cup T)^{[i]}\in {\cal{LR}}^{(i)}$ 
with $T^{[i]}$ of weight $\nu^{(i)}$. Consider the $i$th row word of $T^{[i]}$ where $V_i$ is the  row subword  restricted to the  entries in $[i-1]$, and  $\widehat \nu_i=\lambda_i-\mu_i-|V_i|$ is the number of entries equal to $i$. Put  $(Y_\mu\cup T)^{[0]}=Y_{\nu^{0}}\cup H^{(0)}:=\emptyset$, $\nu^{0}:=0$,
$\bar\phi_{\emptyset}=id$  and  $\rho_1^{(0)}(\emptyset):=\emptyset$.
Then, for $i=1,\dots,n$, it holds
\begin{align}
\rho_1^{(i)}[(Y_\mu\cup T)^{[i]}]&=(\bar\chi_i^{\mu_i}\circ\bar\phi_{V_i}\circ\bar\omega_i^{\widehat\nu_i})\circ\rho^{(i-1)}[(Y_\mu\cup T)^{[i-1]}]
\label{introd:mainrecursion}
\\
&=\bar\chi_i^{\mu_i}\circ\bar\omega_i^{\widehat\nu_i}\circ\bar\phi_{V_i}(Y_{\nu^{(i-1)}}\cup H^{(i-1)})
=Y_{\nu^{(i)}}\cup H^{(i)}\in {\cal{LR}}^{(i)},\label{introd:mainrecursionx+}
\end{align}
where $\bar\omega_i^{\widehat\nu_i}$ adds the $i$th row word $i^{\widehat\nu_i}$ to $ Y_{\nu^{(i-1)}}$, $\bar\chi_i^{\mu_i}$ adds the row word $i^{\mu_i}$ at the end  of the $i$th row of $\bar\phi_{V_i}\circ\bar\omega_i^{\widehat\nu_i}(Y_{\nu^{(i-1)}}\cup H^{(i-1)})$  and  $H^{(i)}\equiv Y_{(\mu_1,\dots,\mu_i)}$.
In particular, all   bumping routes of $\bar\phi_{V_i}$ are pairwise disjoint and terminate in the $i$th row.
\end{thm}
This theorem is illustrated in Section \ref{subsec:skew}.
We observe that the recursive internal insertion presentation of switching $\rho_1^{(n)}$, supplied with \emph{add operators} $\bar\omega_i$ and $\bar\chi_i$, on the ballot pair $Y_\mu\cup T$ in \eqref{introd:mainrecursion}, \eqref{introd:mainrecursionx},
\begin{align}\label{barho}
\rho_1^{(n)}(Y_\mu\cup T)&=(\bar\chi_n^{\mu_n}\circ\bar\phi_{V_{n}}\circ\bar\omega_n^{\widehat\nu_n})\circ\cdots \circ (\bar\chi_2^{\mu_2}\circ\bar\phi_{V_{2}}\circ\bar\omega_2^{\widehat\nu_2})
\circ(\bar\chi_1^{\mu_1}\circ\bar\omega_1^{\widehat\nu_1})\,(\emptyset)
\\
&=Y_\nu\cup H,\; H\equiv Y_\mu,\nonumber
\end{align}
 also produces the \emph{companion tableau} $G_\nu(T)$ or the Gelfand-Tsetlin (GT) pattern of type $\nu$ and content $\lambda-\mu$ of $T$, defined by the nested sequence of partitions $\nu^{(1)}\subseteq \nu^{(2)}\subseteq \cdots\subseteq\nu^{(n)}=\nu$ such that $\nu^{(i)}$ is the content of the ballot tableu $T^{[i]}\in {\cal{LR}}^{(i)}$, for $i=1,\dots,n$. Since $\rho_1^{(n)}$ is an involution,
$\rho_1^{(n)}(Y_\nu\cup H)=Y_\mu\cup T$, this allows another presentation of the switching commuter $\rho_1^{(n)}$ via reverse internal row insertion that we call \emph{  deletion operator} $\rho^{(n)}$ in \cite{akt16}. Deletion operator $\rho^{(n)}$ ($\rho_3$ in \cite{pakvallejo}) just reverses the process as in \eqref{introd:mainrecursion}, \eqref{introd:mainrecursionx} and gives
\begin{align}\rho^{(n)}(Y_\nu\cup H)=Y_\mu\cup T=\rho_1^{(n)}(Y_\nu\cup H)\label{barholga}\end{align}
by producing the GT pattern $G_\nu(T)$ of type $\nu$  and content $\lambda-\mu$ of $T$. This reverse process coincides with  the deletion operations  as explained for $\rho^{(n)}$ in  \cite{az1,akt16} and translated for hives as $\sigma^{(n)}$ in \cite{akt16,tka18}.

For $i=1,\dots,n,$, let $\bar\theta_i:=\bar\chi_i^{\mu_i}\bar\phi_{V_{i}}\bar\omega_i^{\widehat\nu_i}$ in \eqref{introd:mainrecursion}, \eqref{introd:mainrecursionx}. The operator defined by
\begin{align}\bar\rho^{(n)}(Y_\mu\cup T):=\bar\theta_n\cdots\bar\theta_2\bar\theta_1\,(\emptyset)=(\bar\chi_n^{\mu_n}\circ\bar\phi_{V_{n}}\circ\bar\omega_n^{\widehat\nu_n})\cdots (\bar\chi_2^{\mu_2}\circ\bar\phi_{V_{2}}\circ\bar\omega_2^{\widehat\nu_2})\circ(\bar\chi_1^{\mu_1}\circ\bar\omega_1^{\widehat\nu_1})\,(\emptyset).
\end{align}
is the translations   of the operator $\bar\sigma^{(n)}=(\sigma^{(n)})^{-1}$ on hives \cite{akt16,tka18} to ballot tableaux. That is, $\bar\rho^{(n)}=(\rho^{(n)})^{-1}$. From Theorem \ref{2} and \eqref{barho}, one has

\begin{align}\rho_1^{(n)}(Y_\mu\cup T)=\bar\rho^{(n)}(Y_\mu\cup T)=Y_\nu\cup H, \mbox{ $H\equiv Y_\mu$ and shape $\lambda/\nu$}.
\end{align}


Since $\rho_1^{(n)}$ is an involution, 
 it follows
\begin{align}\rho_1^{(n)}=\bar\rho^{(n)}=\rho^{(n)}.\end{align}

Now, for $i=1,\dots,n,$ let  $$\delta_i:={(\bar\theta_i)}^{-1}=(\bar\omega_i^{\hat\nu_i})^{-1}(\bar\phi_{V_{i}})^{-1}(\bar\chi_i^{\mu_i})^{-1}$$ where the action of this operator is realised through  deletion operations in the reverse process of \eqref{introd:mainrecursion}. Note, from Theorem \ref{2}, all   bumping routes of $\bar\phi_{V_i}$ are pairwise disjoint and terminate in the $i$th row, hence $(\bar\phi_{V_{i}})^{-1}$ is a reverse internal insertion operation and starts in row $i$. Since, $\bar\rho^{(n)}$ is reversible by reverse row internal insertion, it defines $\rho^{(n)}$ via the GT pattern $G_\nu$ that it produces.


We write $$\rho^{(n)}(Y_\nu\cup H)=Y_\mu\cup T$$  in the sense that  $\rho^{(n)}$ is defined by the production
of  the  GT pattern of type $\nu$ of $T$ given by the sequence of inner shapes in
$$Y_\nu\cup H,\,\delta_n(Y_\nu\cup H),\,\delta_{n-1}\delta_n(Y_\nu\cup H), \dots,\,\delta_2\cdots\delta_{n-1}\delta_n(Y_\nu\cup H),\,\delta_1 \delta_2\cdots\delta_{n-1}\delta_n(Y_\nu\cup H)=\emptyset.$$

 Hence $$\rho^{(n)}(Y_\nu\cup H)=\rho_1^{(n)}(Y_\nu\cup H)$$



\noindent  We also write $$\rho^{(n)}(Y_\mu\cup T)=Y_\nu\cup H$$  in the sense that  $\rho^{(n)}$ is defined by the production
of  the  GT pattern of type $\mu$ of $H$ given by the sequence of inner shapes in
$$Y_\mu\cup T,\,\delta_n(Y_\mu\cup T),\, \dots,\,\delta_2\cdots\delta_{n-1}\delta_n(Y_\mu\cup T),\,\delta_1 \delta_2\cdots\delta_{n-1}\delta_n(Y_\mu\cup T)=\emptyset.$$

 Then
$$\rho_1^{(n)}(Y_\mu\cup T)=\bar\rho^{(n)}(Y_\mu\cup T)=\rho^{(n)}(Y_\mu\cup T)=Y_\nu\cup H.$$

 Thereby $\rho^{(n)}$ and  $\bar\rho^{(n)}$ just provide another method to compute switching  on ballot tableau pairs as well as the GT pattern pair in the Henriques--Kamnitzer crystal $Com_{HK}$ and hive $Com^h_{HK}$ commuters \cite{HK2,HenKam}.
From the bijection, denoted $\varphi$, between hives and ballot tableaux  \cite{fultonbuch,akt16,tka18,krv21,sathishtorres}, one has the following corollary. (We warn the reader that a hive allows several representations, namely,  vertex representation, as in \cite{kt1}, \cite{fultonbuch} (and \cite{krv21} with a flag condition),  edge representation
as introduced by~\cite{KTT1}, and gradient representation as in \cite{akt16}.)
\begin{cor} The following commuters on ballot tableaux or hives coincide and are involutions
\begin{align}\rho_1=\rho_2=\rho'_2=\rho=\bar\rho=\rho_3=Com_{HK}\\
\sigma=\bar\sigma=Com^h_{HK}\\
\varphi\sigma=Com_{HK}.
\end{align}

\end{cor}
\subsection{ Lecouvey-Lenart and Kumar-Torres bijections between Sundaram and Kwon branching models coincide}
{ Let $LR(\lambda/\mu, \nu)$ be te set of ballot (or LR) tableaux of shape $\lambda/\mu$ and content $\nu$. Given $T\in LR(\lambda/\mu, \nu)$, we may associate a pair $(G_\mu(T),G_\nu(T))$ of  semistandard tableaux (or GT), called the left and right companions of $T$, of shape $\mu$ and content the reverse of $\lambda-\nu$, respectively of shape $\nu$ and content $\lambda-\mu$.}
The right companion map $c$ \cite{sathishtorres} induces a bijection  between the set $ LR(\lambda/\mu, \nu)$ of ballot tableaux of skew shape $\lambda/\mu$ and content $\nu$
 and the set $LR^\lambda_{\mu,\nu}$ of $\mu$-dominant semi-standard Young tableaux $G_\nu$ of shape $\nu$ and
content $\lambda-\mu$.  The left companion { map $c^-$ map injects each  $T\in LR(\lambda/\mu, \nu)$   to its left companion tableau $G_\mu(T)$ of shape $\mu$ and content reverse  $\lambda-\nu$. The left companion tableau $G_\mu$ of $ T\in LR(\lambda/\mu, \nu)$ can be characterized by the $\nu$-dominance of its contre-tableau, that is, the contre-tableau of $G_\mu$, with  shape $\mu$ and weight $\lambda-\nu$, is $\nu$-dominant.
The set ${}^-LR^\lambda_{\mu,\nu}$ denotes the set of  semi-standard Young tableaux  of shape $\mu$ and
content the reverse $\lambda-\nu$ whose contre-tableau is $\nu$-dominant.  The rectification of the contre-tableau of $G_\mu$ is exactly
 $\xi (G_\mu)$ that is still $\nu$-dominant because  Knuth equivalence (or rectification) preserves $\nu$-dominance. Therefore $\xi({}^-LR^\lambda_{\mu,\nu})=LR^\lambda_{\nu,\mu}$.


\begin{cor} The Henriques-Kamnitzer  symmetry  $ LR^\lambda_{\mu,\nu}\overset{\sim}\longrightarrow LR^\lambda_{\nu,\mu}$ can be defined by

\begin{align}G_\nu \overset{c^{-1}}\rightarrow T \overset{c^-}\rightarrow G_\mu(T)\rightarrow \xi( G_\mu(T)),\label{hksym}\end{align}    where $\xi$ is the Lusztig-Sch\"utzenberger involution. That is,
 $G$ is the left companion of $T\in  LR(\lambda/\mu, \nu)$ if and only if $\xi(G)$ is the right companion of $\rho_1(T)=\rho^{(n)}(T)$.  Moreover, $\xi(G_\mu)$ can be calculated by the reverse Sagan-Stanley internal insertion: it is the GT pattern of shape $\mu$ and content $\lambda-\nu$ produced  by the sequence of inner shapes in

\begin{align}Y_\mu\cup T,\,\delta_n(Y_\mu\cup T),\, \dots,\,\delta_2\cdots\delta_{n-1}\delta_n(Y_\mu\cup T),\,\delta_1 \delta_2\cdots\delta_{n-1}\delta_n(Y_\mu\cup T)=\emptyset.
\end{align}
\end{cor}


 The commuter \eqref{hksym} concerning ballot tableaux can be translated for hives because $LR^\lambda_{\mu,\nu}$ and $Hive(\mu, \nu,\lambda)$ are in bijection thanks to \cite{fultonbuch}. (We refrain from defining here hives and refer the reader to \cite{fultonbuch,krv21,sathishtorres}.) It follows then that from a hive $h\in Hive(\nu,\mu, \lambda)$ we can injectively obtain simultaneously a
$\nu$-dominant tableau $P_\mu$ of shape $\mu$ and weight $\lambda-\nu$, that is, $P_\mu\in LR^{\lambda}_{\nu,\mu}$
and a
$\mu$-dominant contretableau $C(P_\nu)$ of shape $\nu$ and weight $\lambda - \mu$, that is, $ P_\nu\in {}^-LR^\lambda_{\nu,\mu}$ and $\xi(P_\nu)\in LR^\lambda_{\mu,\nu}$.

Fom  the coincidence of LR commuters and the work of Kumar-Torres \cite{sathishtorres}, \cite{sathishkumarb} on flagged hives  by Kushwaha–Raghavan–Viswanath \cite{krv21,krsv24},  the Lecouvey-Lenart conjecture \cite{leclen} on  bijections between the Kwon \cite{kwon18} and Sundaram  \cite{sundaram} branching models is settled. We refer the reader to Section \ref{conjecture} for notation and relevant definitions.
The Lenart-Lecouvey conjecture says that the  bijection defined by the Hendriques-Kamnitzer LR commuter between $LR^\lambda_{\mu,\nu}$
 and $LR^\lambda_{\nu,\mu}$
restricts
to a bijection between $LRS(\lambda/\mu, \nu)$, the set of LR tableaux in $ LR(\lambda/\mu, \nu)$ satisfying
the Sundaram property, and $LRK^\lambda_{\nu,\mu}$, the set of tableaux in $ LR^\lambda_{\nu, \mu}$ such that their {Sch\"utzenberger evacuation} $evac_{2n}$ (or Lusztig-Sch\"utzenberger involution $\xi$) within the crystal $B(\mu,2n)$,  satisfy the Kwon
property. This amounts to say that the left companions of $LRS(\lambda/\mu, \nu)$ are Kwon tableaux.  Kumar-Torres \cite{sathishtorres,sathishkumarb} show then that the flagged hive by Kushwaha–Raghavan–Viswanath when restricted to a Sundaram LR tableau exhibits  its GT pattern pair.  

Let
$LRS(\lambda,\mu):=\bigcup LRS(\lambda/\mu,\nu),$  and $LRK(\lambda,\mu):=\bigcup LRK^\lambda_{\nu,\mu}$ where in both cases the union is taken over all even partitions $\nu$, that is, $\nu_{2i-1}=\nu_{2i}$, $i\ge 1$.
\begin{thm}\cite{sathishtorres,sathishkumarb} The bijection of Kushwaha–Raghavan–Viswanath \cite{krv21} between $LR^\lambda_{\mu,\nu}$ and
$ LR^\lambda_{\nu,\mu}$
restricts to a bijection between $LRS(\lambda/\mu, \nu)$ and $LRK^\lambda_{\nu,\mu}$.
\end{thm}

\begin{cor}[Corollary \ref{cor:llkt}] The Kumar-Torres bijection
\begin{align*}LR(\lambda/\mu, \nu)\overset{\sim}\longrightarrow  LR_{\mu,\nu}^\lambda
\overset{U}\longrightarrow LR_{\nu,\mu}^\lambda\end{align*}
where $U$ is the Kushwaha–Raghavan–Viswanath symmetry, and the Lecouvey-Lenart bijection
\begin{align*}LR(\lambda/\mu, \nu)\overset{\sim}\longrightarrow  LR_{\mu,\nu}^\lambda
\overset{U'}\longrightarrow LR_{\nu,\mu}^\lambda\end{align*}
where $U'$ is the LR commuter defined by Henriques–Kamnitzer, coincide. Thereby both restrict to a bijection
between $LRS(\lambda/\mu, \nu)$ and $LRK^\lambda_{\nu,\mu}$ and to $$LRS(\lambda, \mu)\overset{\sim}\longrightarrow  LRK(\lambda, \mu).$$
\end{cor}

   The Littlewood-Richardson commuter based on internal (or reversal) row insertion operations was first introduced in \cite{az1} and called $\rho_3$ in \cite{pakvallejo}. The involutive nature of this commuter for tableaux and hives was  completely detailed in \cite{akt16, tka18} without making recourse of the BSS switching involution. Its coincidence with tableau switching was foreseen in \cite{az2} which we fulfill here.

\subsection{Organization of the paper} The paper is organized in seven sections. In Sections \ref{sec:tabasic} and  \ref{sec:innercorner}, we fix the basic notation to work with,  introduce  our main definitions  and  recall the skew RSK internal insertion correspondence.  In particular, in Section \ref{sec:lrcompanion} we recall the companion pair of a ballot tableau and its importance on characterizing LR commuters. In Section \ref{sec:preserver} we provide a preserver for the $P$-tableau.
Theorem \ref{1} is proved in Section \ref{sec: prooftheorem1} as a consequence of several lemmata.
In Section \ref{switchingllkt}  the recursive presentation of tableau switching on ballot tableau pairs in terms of the Sagan-Stanley internal row insertion is worked out;    the Lecouvey-Lenart conjecture is  settled as a consequence of the coincidence of LR commuters and  the major contribution by Kumar-Torres bijection. In Section \ref{recursionswitch} the recursion on ballot tableau pairs is shown in Theorem  \ref{recursion}.
Theorem \ref{2}  (Main Theorem \ref{th:main})  is proved in Section \ref{sec:proof2}.

\medskip
\bigskip

\emph{This paper is an extension  of the arXiv preprint \cite{az18}, also announced in \cite{tka18}, with further results and applications.}

\bigskip
\section*{Acknowledgements}The author thanks the hospitality of the University of Vienna where her sabbatical leaving took place in the academic year 2015-2016. 
She also benefitted of many fruitful discussions with Ronald C. King and Itaru Terada  while enrolled in the previous projects \cite{akt16,tka18}. In particular, she owes  to Itaru Terada to pointing out the relationship between the Sagan-Stanley internal insertion operations in \cite{ss} and those used by the author in \cite{az1,az2} as already mentioned in \cite{akt16,tka18}. She also thanks to Bruce Sagan for informing that no  conditions for the equality of the $P$-tableau in the skew RSK correspondence were known at the time of her  arxiv post \cite{az18};
 and to Sathish Kumar and Jacinta Torres for letting her know their work \cite{sathishtorres, sathishkumarb}  in Sapporo during the FPSAC25 conference}

This work was financially supported by the Fundação para a Ciência e a Tecnologia (Portuguese Foundation for Science and Technology) under the scope of the projects UID/00324/2025 \linebreak (https://doi.org/10.54499/UID/00324/2025) (Centre for Mathematics of the University of Coimbra) and by the FCT sabbatical grant SFRH/BSAB/113584/2015.


\section{Preliminary  definitions}\label{sec:tabasic}
\subsection{Basics on Young tableaux} As usual $[n]$  denotes the set of positive integers $\{1,\dots, n\}$, $n\ge 1$,  and, if $1\le d\le n$, $[d,n]$ denotes the set $\{d,\dots,n\}$.  A partition (or a normal shape) $\mu$ is  a non negative integer vector in weakly decreasing order, $(\mu_1\ge \cdots\ge \mu_n\ge 0)$. It is identified with the Young diagram of shape $\mu$, in the English convention, that is, $n$ left justified rows of boxes with $\mu_i$ boxes in the $i$th row, for each  $i$, numbering the rows and the columns in matrix style. The box or cell of the Young diagram in row $i$ and column $j$ will be denoted $(i,j)$ with $1\le j\le \mu_i$.
Partitions are usually   denoted by lowercase Greek letters. We write $|\mu|:=\mu_1+\cdots+\mu_n$ for the number that $\mu$ partitions, and the number of positive parts in this summand  is the length $\ell(\mu)\le n$ of   $\mu$. The unique partition of length 
$0$ is the null partition $(0)$, identified with  $\emptyset$, the unique empty Young diagram. A \emph{corner} of a Young diagram is a cell such that its removal still leaves a Young diagram.

For Young diagrams $\mu\subseteq \lambda$, the skew partition (or skew shape) $\lambda/\mu$ is the set-theoretic difference $\lambda\setminus \mu$.
A semistandard
Young tableau (SSYT) $T$ of shape $\lambda/\mu$ is a filling of the boxes of $\lambda/\mu$ over a  finite subset $[n]$ of the   positive integers, such that the labels of each row weakly increase from left to right and the labels of each column strictly increase from top to
bottom. The skew tableau $T$ comprises  an \emph{inner border}, defined by the unfilled \emph{inner shape} $\mu$, a filled skew shape $\lambda/\mu$, and an \emph{outer  border}, defined by the \emph{outer shape} $\lambda$. The labels of $T$ are often referred to elements or entries of $T$.  We denote by  $\{T\}$ \emph{the multiset of entries of $T$} counting the number of repetitions of each entry. If $\mu=\emptyset$,  $T$ is of partition (normal) shape.
The unique empty skew tableau $\mu/\mu:=\emptyset_\mu$, is the Young diagram of shape $\mu$.

 The  {\em (row) reading word}  of $T$ is the word  $w(T)$ on the alphabet $[n]$, read left to right across rows of $T$ taken in turns from bottom to top. When needed, we also consider the \emph{Kashiwara reading word} of $T$, also called  north-western column reading  word of $T$, that is,  the word read from right to left across columns from top to bottom.  The \emph{content or weight} of $T$ is the content of its reading word $w(T)$, that is, the vector $\gamma=(\gamma_1,\dots,\gamma_n)$ where $\gamma_i$ records the number of $i$'s in $T$, for all $i$ in the given alphabet. The length $|w(T)|$ of $w(T)$ is the number of letters which appear in $w(T)$. Equivalently, 
 $|\lambda|-|\mu|=|w(T)|=|\gamma|:=\gamma_1+\cdots+\gamma_s+\cdots$. As usual, given the words $u$ and $v$ over an alphabet,  $uv$ denotes their concatenation. A word is said to be a \emph{row word} if its letters weakly increase from right to left.
  The set of all SSYT's of shape $\lambda/\mu$ is denoted by $YT(\lambda/\mu)$. If we want to emphasize that the labels of the entries range on the set $\{1,\dots,n\}$ then we also write $YT(\lambda/\mu, n)$.
 For an illustration see \eqref{ex:readingword}.

  Noting that $T\in YT(\lambda/\mu, n)$ is also realized by a sequence of nested partitions $\mu\subseteq \lambda^{(1)}\subseteq\cdots\subseteq \lambda^{(n)}=\lambda$ where $ \lambda^{(m)}/\mu$ defines  the filling of the boxes of $T$ on the alphabet $[m]$, for  $1\le m\le n$,  the \emph{restriction of $T$  to the alphabet $[m]$},  $T_{|[m]}\in YT(\lambda^{(m)}/\mu, n)$, is the subtableau of $T$ of content $(\gamma_1,\dots,\gamma_m)$, precisely realized by the subsequence $\mu\subseteq \lambda^{(1)}\subseteq\cdots\subseteq \lambda^{(m)}$.

   A tableau in  $YT(\lambda/\mu)$ with $|\gamma|$ boxes is said to be  {\em standard}  if  the entries are the
numbers from $1$ to $|\gamma|$, each occurring once. The {\em standard order} of the boxes on a semistandard Young tableau
 is given by the numerical ordering of the labels with priority, in the case of equality,
given by rule  \emph{ southwest=smaller}, \emph{northeast=larger}.

Given $U\in YT(\lambda/\mu)$ of weight $\gamma$, the \emph{standardisation} of $U$ is the \emph{standard tableau} $\mathrm{std}~ U$ obtained by the standard order of the boxes of $U$. This means to renumber the entries of $U$ in numerical order from $1$ to $|\gamma|$, and, in case of equal entries,  regard those to the left as smaller than those to the right. The tableau $U$ is easily recovered from $\mathrm{std}~ U$ and its weight $\gamma$. For an illustration see Example \ref{ex:companionword}.

\subsection{Basic calculus on Young tableaux}\label{sec:knuth}
Recall the \emph{ Schensted row insertion} (here also called external Schensted row insertion)  takes a SSYT  $T$ of partition shape,
and an element  $m$ in the $T$-alphabet, and constructs a new tableau, denoted $T \bigcdot m$. For a word $w=w_1\cdots w_s$ on the  $T$-alphabet, we  recursively  define the new tableau $T \bigcdot w=(\cdots ((T\bigcdot  w_1)\bigcdot w_2)\bigcdot\cdots)\bigcdot w_s$ \cite{fulton}.
Recall the  elementary \emph{Knuth transformations} on words over a totally ordered alphabet and its compatibility with row Schensted insertion \cite{fulton, stanley}. For letters $x,y,z$ in a  totally ordered alphabet, an elementary Knuth transformation is governed by the rules below. As usual we write $\equiv$ for Knuth equivalence between words,
\begin{align}\label{knuthrow1}yzx \equiv yxz, \mbox{   if  $yz \bigcdot x=\YT{0.16in}{}{
{y,z},
}\bigcdot \YT{0.16in}{}{{x},
}=\YT{0.15in}{}{
{x,z},
{y},
}=\YT{0.16in}{}{{y},
}\bigcdot\YT{0.16in}{}{
{x,z},
}=y\bigcdot xz$ , equivalently, $x < y \le z$},
\end{align} and

\begin{align}\label{knuthrow2}xzy \equiv  zxy \mbox{ if $xz \bigcdot y=\YT{0.16in}{}{
{x,z},
} \bigcdot \YT{0.16in}{}{{y},
}=\YT{0.15in}{}{
{x,y},
{z},
}=\YT{0.16in}{}{{z},
} \bigcdot \YT{0.15in}{}{
{x,y},
}=z \bigcdot xy$, equivalently, $x \le y < z$}.
\end{align}


 Two words $w$ and $w'$ are said to be \emph{Knuth equivalent}, $w\equiv w'$, if they can be transformed into each other by a sequence of elementary Knuth transformations.
 Two skew-tableaux $T$ and $U$  are Knuth equivalent $T\equiv U$ if and only if $w(T)\equiv w(U)$. Equivalently, $T$ and $U$  have the same \emph{rectification}, that is, the insertion tableaux $P(w(T))=P(w(U))$ obtained by row Schensted insertion of the words $w(T)$ respectively  $w(U)$ \cite{fulton, stanley}.

\begin{obs}The row reading and column reading words of $T$ are Knuth equivalent. The Kashiwara reading word  and the reverse row reading word of $T$ are Knuth equivalent. The $P$-tableau of the Schensted row insertion of a  word $w=w_1 w_2\cdots w_s$ equals the $P$-tableau of the \emph{Schensted column insertion} of the reverse word of $w$, $w_s\cdots w_2 w_1$. That is, $P(w)=w_1\bigcdot w_2 \bigcdot\cdots \bigcdot w_s= w_s\leftarrow \cdots \leftarrow w_2\leftarrow w_1$ where $\leftarrow$ means \emph{Schensted column insertion}. For instance, in \eqref{knuthrow1}, $yz \bigcdot x=x\leftarrow z\leftarrow y=y\bigcdot xz=z\leftarrow x\leftarrow y$. (Similarly for \eqref{knuthrow1}.)
Thereby, the $P$-tableau of the reading word of $T$, and the $P$-tableau of the Kashiwara reading word in the {Schensted column insertion} is the same.
\end{obs}
If nothing in contrary is said we always consider Schensted row insertion.

Either for the  Schensted row insertion or  Schensted column insertion version of RSK, we have the following result:
 two words $\pi$ and $\pi'$  are Knuth  equivalent if and only if their
$P$-tableaux under  RSK correspondence  are equal $P(\pi) = P(\pi')$ \cite{fulton, stanley}. In the usual RSK,  Knuth relations  completely characterize the words having the same $P$-tableau.

\subsection{Ballot tableaux and companion pairs}\label{sec:lrcompanion} We follow closely the references \cite{akt16,tka18,azkoma25} and we refer to them for additional details. A SSYT tableau is  said to be  {\em ballot} or \emph{Littlewood-Richardson}(LR) tableau  if the content of each suffix of its   reading word, that is, the content of the subword read backwards from  the end to any letter,   is a partition. Such  a word is also called Yamanouchi, ballot or reverse lattice word. The ballot tableau of shape $\mu$  is also called \emph{Yamanouchi tableau}, $Y_\mu$. In other words, it is the unique tableau of shape and content $\mu$, that is, the entries of row $i$ consist of $i$'s.  Denote by $LR(\lambda/\mu,\nu)$, where $\mu,\nu\subseteq \lambda$ are partitions, the subset of $YT(\lambda/\mu)$ consisting of ballot tableaux  of shape $\lambda/\mu$ and content $\nu$. That is $T\in LR(\lambda/\mu,\nu)$ if and only if $T\equiv Y_\nu$.

\begin{ex}For instance, in \eqref{ex:readingword}, $H$ and $T$ are SSYT's of shape $(4,3,2)/(2,1,0)$, and $Y_{(3,2,1)}$ is the Yamanouchi tableau of shape $(3,2,1)$.  The  reading words of $H, \,T$ and $Y$ are $121312,\, 231211$ and $ 322111$ respectively. The two last words are Yamanouchi but not the first. Therefore $T$ and $Y$ are ballot but  $H$ is not,
\begin{align}\label{ex:readingword}
H=\YT{0.16in}{}{
{,,1,2},
{,1,3},
{1,2},
}\quad
T=\YT{0.16in}{}{
{,,1,1},
{,1,2},
{2,3},
} \quad
Y_{(3,2,1)}=\YT{0.16in}{}{
{1,1,1},
{2,2},
{3},
} \quad T\equiv Y_{(3,2,1)}.
\end{align}
\end{ex}

A \emph{key tableau} of shape $\pi$ is the tableau $Y_\alpha$ of  shape $\pi$ whose content $\alpha$ is a permutation of its shape $\pi$ or the columns are nested as sets. In particular,  when $\alpha$ is the reverse of $\pi$, $rev\,\pi$,  we say that $Y_{rev\,\pi}$ is  the \emph{anti-Yamanouchi} tableau of shape $\pi$. For instance, for $n=4$ and $\mu=(2,1,0,0)$ respectively $\nu=(3,2,1,0)$, one has the anti-Yamanouchi tableaux of shapes $\mu$ respectively $\nu$ where the weight is the reverse of the shape

\begin{align}
Y_{(0,0,1,2)}=\YT{0.16in}{}{
{3,4},
{4},
}=\mathrm{evac}_4Y_{(2,1,0,0)},\quad Y_{(0,1,2,3)}=\YT{0.16in}{}{
{2,3,4},
{3,4},
{{4}},
}=\mathrm{evac}_4Y_{(3,2,1,0)}.
\end{align}


Let $B(\pi, n)$ be the crystal of  tableaux of shape $\pi$ on the alphabet $[n]$, and $\xi=\mathrm{evac}_n$ the Sch\"utzenberger-Lusztig involution (or Sch\"utzenberger evacuation) on that crystal. The highest weight element of $B(\pi,n)$ is $Y_\pi$ and the lowest element is $\xi(Y_\pi)=Y_{rev\,\pi}$.

 A SSYT tableau $T$ is  said to be  {\em anti-ballot} or an \emph{opposite Littlewood-Richardson}(LR) tableau  if its reading  word $w(T)$ is anti-Yamanouchi, which means the content of each prefix, that is, the content of the subword read left to right  from  the beginning to any letter   is a reverse partition.
  A tableau $T\in YT(\lambda/\mu,n)$ of content $\mathrm{rev}\nu$ is said to be anti-ballot if $T\equiv evac_n(Y_\nu)= Y_{\mathrm{}rev\,\nu}$.
 Denote by $opLR(\lambda/\mu, rev\nu)$, where $\mu,\nu\subseteq \lambda$ are partitions, the subset of $YT(\lambda/\mu,n)$ consisting of  tableaux $T$ of shape $\lambda/\mu$ and content $rev\nu$  such that $T\equiv evac_n(Y_\nu)=Y_{\mathrm{rev}\nu}$.
 \begin{ex} Let $n=4$ and $T\in LR(\lambda/\mu,\nu)$ as in the previous example. Recalling that  \emph{reversal} \cite{bss96,azkoma25} is the version of Sch\"utzenberger evacuation for   skew-tableaux, one has that  \emph{reversal} of $T$ is in $opLR(\lambda/\mu, rev\nu)$ with reading word $443324$ equals to
 \begin{align}\YT{0.16in}{}{
{, ,2,4},
{,3,3},
{{4},4},
}\equiv Y_{(0,1,2,3)}.
\end{align}
\end{ex}


Given $n\ge 1$ and a ballot tableau $T\in LR(\lambda/\mu,\nu)$ with $\ell(\lambda)\le n$, one associates a \emph{pair} of \emph{companion tableaux} or \emph{Gelfand-Tsetlin patterns} (GT) $(G_\mu(T),G_\nu(T))$, \emph{left companion} respectively  \emph{right companion} uniquely determined by $T$.  (We identify a Gelfand-Tsetlin pattern of base (or type) $\kappa$ with its natural tableau presentation of shape $\kappa$.)

\begin{obs} In fact $T\in LR(\lambda/\mu,\nu)$ is completely specified either by its left or right companion tableau. In addition, they are linearly bijectively related, see \cite{pakvallejo, azkoma25}.
\end{obs}

The left companion $G_\mu(T)$ of shape $\mu$ and content $rev(\lambda-\nu)$ is  obtained from $T$ by recording the sequence of partitions $\mu^{(n-r+1)}$ giving the
shapes occupied by the entries including the empty entries identified with 0 $< r$ in rows $r, r + 1, \dots , n$ of $T$, for $r = 1, 2, \dots, n$. We get then the nested sequence of partitions $\mu^{(1)}\subseteq \mu^{(2)}\subseteq \cdots \subseteq\mu^{(n)}=\mu$ defining $G_\mu(T)$.

 The right  companion $G_\nu(T)$ of shape $\nu$ and content $\lambda-\mu$ is  obtained from $T$ by recording the sequence of partitions $\nu^{r}$ giving the
shapes occupied by the positive entries $\le r$ in rows $1, 2, \dots , r$ of $T$, for $r = 1, 2, \dots, n$.
We get then the nested sequence of partitions $\nu^{(1)}\subseteq \nu^{(2)}\subseteq \cdots \subseteq\nu^{(n)}=\nu$ defining $G_\nu(T)$. Equivalently, the row $r$ of $G_\nu(T)$ records the row coordinates of the $r$-cells of $T$ for $r=1,\dots,n$ \cite{akt16}.

\begin{ex} Let $n=6$, $\lambda=(6,5,5,4,3,0)$, $\nu=(4,4,3,2,0,0) $, $\mu=(4,3,2,1,0,0)$ and $T\in  LR(\lambda/\mu,\nu)$ as below. We illustrate $T$ with its GT pattern pair $(G_\mu(T),G_\nu(T))\in B(\mu,6)\otimes B(\nu,6)$
\begin{align} &T=\YT{0.15in}{}{
 {{},{},{},{},{1},{1}},
 {{},{},{},{1},{2}},
 {{},{},{1},{2},{3}},
 {{},{2},{3},{4}},
 {{2},{3},{4}},
},\quad
G_\mu(T)=  \YT{0.15in}{}{
 {2,2,{2},{4}},
 {3,{3},{5}},
 {4,{6}},
 {{6}},
}, \quad
\mbox{$\vcenter{\hbox{\begin{tikzpicture}[x={(1cm*0.4,-1.7320508cm*0.4)},y={(1cm*0.4,1.7320508cm*0.4)}]
\path(0,0)--node[pos=0.45]{$4$}(1,1);
\path(1,1)--node[pos=0.45]{$3$}(2,2);
\path(2,2)--node[pos=0.45]{$2$}(3,3);
\path(3,3)--node[pos=0.45]{$1$}(4,4);
\path(4,4)--node[pos=0.45]{$0$}(5,5);
\path(5,5)--node[pos=0.45]{$0$}(6,6);
\path(0,1)--node[pos=0.45]{$4$}(1,2);
\path(1,2)--node[pos=0.45]{$3$}(2,3);
\path(2,3)--node[pos=0.45]{$1$}(3,4);
\path(3,4)--node[pos=0.45]{$0$}(4,5);
\path(4,5)--node[pos=0.45]{$0$}(5,6);
\path(0,2)--node[pos=0.45]{$4$}(1,3);
\path(1,3)--node[pos=0.45]{$2$}(2,4);
\path(2,4)--node[pos=0.45]{$1$}(3,5);
\path(3,5)--node[pos=0.45]{$0$}(4,6);
\path(0,3)--node[pos=0.45]{$3$}(1,4);
\path(1,4)--node[pos=0.45]{$2$}(2,5);
\path(2,5)--node[pos=0.45]{$0$}(3,6);
\path(0,4)--node[pos=0.45]{$3$}(1,5);
\path(1,5)--node[pos=0.45]{$0$}(2,6);
\path(0,5)--node[pos=0.45]{$0$}(1,6);
\end{tikzpicture}}}$
}\nonumber
\\& \mbox{GT pattern of type $\mu$ and weight $rev(\lambda-\nu)=(0,3,2,2,1,2)$ defined by the nested sequence of partitions}\nonumber\\
\nonumber\\
&\mu=\mu^{(6)}=(4,3,2,1,0,0)\supseteq \mu^{(5)}=(4,3,1,0,0)\supseteq \mu^{(4)}=(4,2,1,0)\supseteq\nonumber \\ &\supseteq\mu^{(3)}=(3,2,0)\supseteq \mu^{(2)}=(3,0)\supseteq \mu^{(1)}=(0)
\end{align}
and \begin{align}
 & G_{\nu}(T)=
\YT{0.15in}{}{
 {1,1,{2},{3}},
 {2,{3},{4},{5}},
 {3,{4},{5}},
 {4,5},
}
\end{align}
The Kashiwara reading word (right to left across columns top to bottom) of
\begin{align}\label{lowestweight}
G_\mu(T)\otimes Y_{\rm{rev}\nu}=G_\mu(T)\otimes Y_{(0,0,2,3,4,4)}\equiv Y_{rev\lambda}
\end{align}

\begin{align}
&G_\mu(T)\otimes Y_{(002344)}\equiv 4252362346\otimes 5645634563456\\
&\equiv G_\mu(T)\leftarrow 5\leftarrow 6\leftarrow 4\leftarrow 5\leftarrow 6\leftarrow 3\leftarrow 4\leftarrow 5\leftarrow 6\leftarrow 3\leftarrow 4 \leftarrow 5\leftarrow 6\\
&=\YT{0.16in}{}{
{2, 2,2,3,4,6},
{3,3,3,4,5},
{{4},4,4,5,6},
{5,5,5,6},
{6,6,6},
}= Y_{rev\lambda}\end{align}

\begin{obs}\label{dualword}$G_\mu\otimes Y_{\rm{rev}\nu}\equiv Y_{rev\lambda}$ is equivalent $ Y_{\nu}\otimes C(G_\mu)\equiv Y_{\lambda}$ where $C(G_\mu)$ is the contre-tableau of $G_\mu$. Note the word of the contre-tableau is the dual word of $w(G(\mu))$, that is, if $w=w_1w_2\cdots w_s$ is a word in the alphabet $[m]$ then the dual is $(m-w_s)\cdots (m-w_2)(m-w_1)$. Equivalently, $w(C(G_\mu))$ is $\nu$-dominant, and  $w(Y_\nu)w(C(G_\mu))$ is a Yamanouchi word of weight $\lambda$ (for details see \cite[Section 5.1, Appendix A.1]{fulton} and \cite[Appendix]{fultonbuch}),
\begin{align}\label{contre}Y_\nu\otimes C(G_\mu)\equiv Y_\lambda.\end{align}
Therefore, since $\xi(G_\mu)$  is the rectification of $C(G_\mu)$, $\xi(G_\mu)\equiv C(G_\mu)$ is $\nu$-dominant,
\begin{align*}&Y_\lambda=Y_\nu\leftarrow C(G_\mu)=Y_{(443200)}\leftarrow \YT{0.16in}{}{
{, ,,1},
{,,1,3},
{,2,{4},4},
{3,5,5,5},
}.
\end{align*}
\end{obs}


The Kashiwara reading word (right to left along columns top to bottom) of $Y_\mu\otimes G_\nu(T)$ is the   Yamanouchi word of weight $\lambda$ whose column insertion gives $Y_\lambda$:
\begin{align}&Y_\mu\otimes G_\nu(T)\equiv 1121231234\otimes 3524513451234\equiv \\
&Y_{(4,3,2,1,0,0)}\leftarrow 3\leftarrow5\leftarrow 2\leftarrow 4\leftarrow 5\leftarrow 1\leftarrow 3\leftarrow 4\leftarrow 5\leftarrow 1\leftarrow 2\leftarrow 3\leftarrow 4 \\
&=Y_\lambda
\end{align}
\end{ex}

\begin{obs}In the row reading word option the corresponding versions are $Y_{rev \nu}\bigcdot G_\mu(T)\equiv Y_{rev \lambda}$ respectively $G_\nu(T)\bigcdot Y_\mu\equiv Y_\lambda$ \cite[Section 5.2, Exercise 3]{fulton}.
\end{obs}

We now collect a few facts that characterize the set of GT-patterns of shape $\mu$ and weight $rev(\lambda-\nu)$ respectively the set of GT-patterns of shape $\nu$ and weight $\lambda-\mu$ each of which specifying the set $LR(\lambda/\nu,\mu)$. We need some notation. We denote the set of those GT-patterns by ${}^-{LR}^\lambda_{\mu,\nu}$ respectively $LR^\lambda_{\mu,\nu}$. Recall that
\begin{align}\label{tensordecomp}B(\mu)\otimes B(\nu)\cong \bigoplus B(\lambda)^{c_{\mu,\nu}^\lambda},
\end{align}
where the sum is taken for all partitions $\lambda\supseteq \mu,\nu$, and $c_{\mu,\nu}^\lambda=\#L(\lambda/\mu,\nu)$.

The set $LR^\lambda_{\mu,\nu}$ (${}^-LR^\lambda_{\mu,\nu}$)
 of semistandard tableaux $ G_\nu$ ($G_\mu$) of shape $\nu$ ($\mu$) and content $\lambda-\mu$ ($\rm{rev}(\lambda-\nu)$) are those in the crystal
$ B(\nu)$ ($ B(\mu)$) such that $Y_\mu \otimes G_\nu$ ($G_\mu\otimes Y_{\rm{rev}\nu}$) is the highest (lowest) weight element of weight $\lambda$ ($\rm{rev}\lambda$) in  a connected component of $B(\mu)\otimes B(\nu)$ isomorphic to $B(\lambda)$.

Therefore, $G_\nu\in LR^\lambda_{\mu,\nu}$ if and only if it satisfies  the left $\mu$-dominance (the Kashiwara readi the word of $G_\nu$ concatenated on the left with the canonical Yamanouchi word of weight $\mu$ (the Kashiwara reading word of $Y_{\mu}$) is  $\equiv Y_{\lambda}$.

On the other hand, $G_\mu\in {}^-LR^\lambda_{\mu,\nu}$  if and only if it satisfies  the right $\rm{rev}\nu$-dominance, that is, the word of $G_\mu$ concatenated on the right with the canonical anti-Yamanouchi word of weight $\rm{rev}\nu$ ( the Kashiwara reading word of $Y_{\rm{rev}\nu}$) is  $\equiv Y_{\rm{rev}\lambda}$. Equivalently, from \eqref{contre}, $Y_\nu\otimes C(G_\mu)\equiv Y_\lambda$, it amounts to say that $\xi(G_\mu)\in LR^\lambda_{\nu,\mu}$.

We know  from \cite{HenKam} that $LR^\lambda_{\mu,\nu}$, ${}^-LR^\lambda_{\mu,\nu}$ and $LR(\lambda/\nu,\mu)$ are in bijection.
In fact, from \cite{HenKam}, $T\in LR(\lambda/\nu,\mu) $ has companion pair $(G_\mu,G_\nu)$ if and only if  $Y_\mu \otimes G_\nu\equiv Y_\lambda$ and  $G_\mu\otimes Y_{rev \nu}\equiv Y_{rev\lambda}$ are the highest and lowest weight elements of a same connected component of $B(\mu)\otimes B(\nu)$ isomorphic to $B(\lambda)$.

 We know  from \cite{pakvallejo,HenKam} that $LR^\lambda_{\mu,\nu}$, ${}^-LR^\lambda_{\mu,\nu}$ and $LR(\lambda/\nu,\mu)$ are in bijection.
In fact, from \cite{HenKam}, the pair $(G_\mu,G_\nu)\in LR^\lambda_{\mu,\nu}\times {}^-LR^\lambda_{\mu,\nu}$ is the companion pair of $ T\in LR(\lambda/\mu,\nu)$ if and only if  $Y_\mu \otimes G_\nu\equiv Y_\lambda$ and  $G_\mu\otimes Y_{rev \nu}\equiv Y_{rev\lambda}$ are the highest respectively lowest weight elements of a same connected component of $B(\mu)\otimes B(\nu)$ isomorphic to $B(\lambda)$ where  $T\in LR(\lambda/\mu,\nu)$ is the recording tableau in the column insertion \cite{Nak,kwon18} of $G_\mu\leftarrow Y_{rev \nu}$ and $Y_\mu\leftarrow G_\nu$.  In other words, each copy of $B(\lambda)$ is uniquely parameterized by a $ T\in LR(\lambda/\mu,\nu)$.
From \cite{kwon18} we then have an RSK version of \eqref{tensordecomp}
\begin{align}\label{tensordecomp2}B(\mu,n)\otimes B(\nu,n)\cong \bigoplus_{\begin{smallmatrix}\lambda\\
T\in LR(\lambda/\mu,\nu)
\end{smallmatrix}} B(\lambda,n)\times \{T\},
\end{align}
where $\lambda$ is taken over all partitions of $n$ such that $\mu,\nu\subseteq\lambda$.

Furthermore,  from Henriques-Kamnitzer commuter \cite{HenKam}, $(G_\mu,G_\nu)$ is the companion pair of some
$ T\in LR(\lambda/\mu,\nu)$ if and only if $(\xi(G_\nu),\xi(G_\mu)) $ is the companion pair of some $H\in LR(\lambda/\nu,\mu)$. Note that the Henriques-Kamnitzer commuter subsumes Remark \ref{dualword}.

 We then have the right companion bijection $c$
\begin{align}c:LR(\lambda/\mu,\nu)\longrightarrow LR_{\mu,\nu}^\lambda, \;T\mapsto G_\nu(T)\end{align}
and the left companion bijection  $c^-$
\begin{align} c^-:LR(\lambda/\mu,\nu)\longrightarrow {}^-LR_{\mu,\nu}^\lambda\;
 T\mapsto G_\mu(T).
 \end{align}
The Henriques-Kamnitzer commuter \cite{HenKam} $Com_{HK}$ can be written as
\begin{align*}Com_{HK}=\xi\circ c^-\circ c^{-1}
\end{align*} or in its left version

\begin{align*}Com_{HK}^-=\xi\circ c\circ {c^-}^{-1}
\end{align*}

\begin{align}\label{comhk}&Com_{HK}:LR_{\mu,\nu}^\lambda\overset{c^{-1}}\longrightarrow LR(\lambda/\mu,\nu)\overset{c^-}\longrightarrow {}^-LR_{\mu,\nu}^\lambda \overset{\xi}\longrightarrow LR_{\nu,\mu}^\lambda,\; G_\nu\mapsto T\mapsto G_\mu\mapsto \xi(G_\mu),
\end{align}

or

\begin{align}\label{comhk-}&Com^-_{HK}:{}^-LR_{\mu,\nu}^\lambda\overset{{c^-}^{-1}}\longrightarrow LR(\lambda/\mu,\nu)\overset{c}\longrightarrow LR_{\mu,\nu}^\lambda \overset{\xi}\longrightarrow {}^-LR_{\nu,\mu}^\lambda,\; G_\mu\mapsto T\mapsto G_\nu\mapsto \xi(G_\nu).
\end{align}

\begin{ex}For $n=6$ this is an illustration of the inverse of map $c^-$ applied to $G_\mu(T)$ to get $T$:
\begin{align}
& G_\mu(T)=  \YT{0.15in}{}{
 {2,2,{2},{4}},
 {3,{3},{5}},
 {4,{6}},
 {{6}},
},\quad \YT{0.15in}{}{
 {,,{},{}},
 {,{},{}},
 {,{}},
 {{}},
},\quad
\YT{0.15in}{}{
 {2,2,{2},{4}},
 {3,{3},{5}},
 {4},
}\quad \YT{0.15in}{}{
 {,,{},{}},
 {,{},{},1},
 {,{},1},
 {{}},
},\quad
\quad\YT{0.15in}{}{
 {2,2,{2},{4}},
 {3,{3}},
 {4},
}\quad \YT{0.15in}{}{
 {,,{},{}},
 {,{},{},1},
 {,{},1,2},
 {{},2},
 {2},
},\nonumber\\
&\YT{0.15in}{}{
 {2,2,{2}},
 {3,{3}},
}
\quad
\YT{0.15in}{}{
 {,,{},{}},
 {,{},{},1},
 {,{},1,2},
 {{},2,3},
 {2,3},
},
\quad
\YT{0.15in}{}{
 {2,2,{2}},
}
\quad
\YT{0.15in}{}{
 {,,{},{}},
 {,{},{},1},
 {,{},1,2},
 {{},2,3},
 {2,3,4},
}=\tilde T,\quad
\YT{0.15in}{}{
 {,,{},{},1,1},
 {,{},{},1,2},
 {,{},1,2,3},
 {{},2,3,4},
 {2,3,4},
}=T
\end{align}
where $\widehat\nu=(2,1,1,1,0)=\lambda-shape(\tilde T)=\lambda-(4,4,4,3,3)$, and $T$ is obtained from $\tilde T$ by adding $\widehat \nu_i$, $i$'s, to row $i$ of $\tilde T$, for $i=1,\dots,\ell(\lambda)$.
\end{ex}
\subsection{The companion word of a skew-tableau }\label{sec:companionwords}

 \begin{defi} \label{comp}Let $U\in YT(\lambda/\mu,n)$ of weight $\gamma$. The   \textit{companion word} of $U$ is defined to be the word $\cal R({U})=\cal R(\std ~U):= u_{|\gamma|}\cdots u_{2}u_{1} $ listing the row indices of the entries of $\std ~U$, from    the bigger  to the smaller.
Equivalently, to construct $u_{|\gamma|}\cdots u_{2}u_{1} $, for $p=1,\dots,|\gamma|$,   put $u_p=i$, if the number $p$ is in the $i$th row of $\std ~U$. The companion word of $T$ factorizes into $n$ maximal row words (possibly with some empty factors) $\cal{R}(U)=\cal{R}_n\cdots \cal{R}_1$ where $\cal R_i$ is the  row word defined by the row coordinates of the $i$-cells in $T$ for $i=1,\dots,n$.
\end{defi}
\begin{ex}\label{ex:companionword} Let

\begin{align*}U=&\YT{0.16in}{}{
{,,,1,3},
{,,{2},4},
{1,2,3},
}\in YT(\lambda/(3,2,0),4),\quad
V=\YT{0.16in}{}{
{,,,4,6},
{,,{5},7},
{4,5,6},
}\in YT(\lambda/(3,2,0),6),\\
 \std~ U=\std~ V=&\YT{0.16in}{}{
{,,,2,6},
{,,{4},7},
{1,3,5},
}\in YT(\lambda/(3,2,0),7).
\end{align*}
The corresponding companion words are
$$\cal R(U)=\cal R(V)=\cal R(\std U)=\cal R(\std V)=2\,13\,23\,13$$
$$\cal R(U)=\cal {R}_4\cal {R}_3\cal {R}_2\cal {R}_1=2\,13\,23\,13$$
$$\cal R(V)=\cal {R}_7\cal {R}_6\cal {R}_5\cal {R}_4\cal {R}_3\cal {R}_2\cal {R}_1=2\,13\,23\,13\,\emptyset\,\emptyset\,\emptyset$$
$$\cal R(\std U)=\cal R(\std V)=\cal {R}_7\cal {R}_6\cal {R}_5\cal {R}_4\cal {R}_3\cal {R}_2\cal {R}_1=2\,1\,3\,2\,3\,1\,3$$
\end{ex}

\begin{prop} Let $U\in YT(\beta/\mu)$ and $\cal {R}(U)$ its companion word. If $\mu=\emptyset$ then $\cal {R}(U)$ is a Yamanouchi word $\equiv Y_\beta$. If $\mu\neq\emptyset$ and $\beta$ has rectangle shape then $\cal {R}(U)$ is a anti-Yamanouchi word of content $(\beta_1^{\ell(\beta)})-\mu=(\beta_1-\mu_1, \beta_1-\mu_2,\dots,\beta_1- \mu_{\ell(\beta)})$ and thus $\equiv Y_{(\beta_1^{\ell(\beta)})-\mu}$.
\end{prop}
\begin{proof} If $\mu=\emptyset$ it is a consequence of the definitions of companion word and Yamanouchi word \cite{fulton}. If $\mu\neq\emptyset$ and $\beta$ has rectangle shape, it is enough to observe that the word $\cal {R}(U)$ is a reverse Yamanouchi word of content  $(\beta_1-\mu_1, \beta_1-\mu_2,\dots,\beta_1- \mu_{\ell(\beta)})$, that is,  in each prefix of $\cal {R}(U)$ the number of $(i+1)$'s is at least equal to the number of $i$'s, for all $i$.
\end{proof}

\begin{ex}The companion words of $Y_{(1,3,4)}$ respectively $U$
$$\YT{0.16in}{}{
{1,2,2,3},
{2,3,3},
{{3}},
}=Y_{(1,3,4)},\quad \cal{R}(Y_{134})=1223\,112\,1,\quad U=\YT{0.16in}{}{
{1,1,2,3},
{2,3,3},
{{4}},
}\quad \cal{R}(U)=3\,122\,12\,11 
$$
are Yamanouchi words of weight $(4,3,1)$.
\end{ex}

\begin{obs} If $T$ is a ballot tableau of weight $\nu$, $\cal{R}(T)$ is the word of the  companion tableau or Gelfand-Tsetlin pattern $G_\nu$ of $T$.
Let $U\in YT(4,3,2)/(2,1)$ be a ballot tableau of weight $\nu=(3,2,1)$
$$U=\YT{0.16in}{}{
{,,1,1},
{,1,2},
{2,3},
}
$$
The companion word of $U$ is $\cal{R}(U)=3\,23\,112=\cal{R}_3\cal{R}_2\cal{R}_1\equiv G_\nu$ which is precisely the reading word of the \emph{companion tableau}
$G_\nu=\YT{0.16in}{}{
{1,1,2},
{2,3},
{3},
}$ of the ballot tableau $U$.

Let $T$ be  a ballot tableau of weight $\nu=(4,4,3,2)$
\begin{eqnarray}\label{gtobs} T&=&\YT{0.15in}{}{
 {,,,,{1},{1}},
 {,,,{1},{2}},
 {,,{1},{2},{3}},
 {,{2},{3},{4}},
 {{2},{3},{4}},
}
\end{eqnarray}
$\cal{R}(T)=\emptyset\,45\,345\,2345\,1123=\cal{R}_5\cal{R}_4\cal{R}_3\cal{R}_2\cal{R}_1\equiv G_\nu=\YT{0.16in}{}{
{1,1,2,3},
{2,3,4,5},
{3,4,5},
{4,5},
}$
\end{obs}
\begin{ex}
 Companion words of rectangle shapes are anti-Yamanouchi words
\begin{align}
U=\YT{0.16in}{}{
{,,,1},
{,1,3,4},
{{3},3,4,5},
}\quad
\std U=\YT{0.16in}{}{
{,,,2},
{,1,5,7},
{{3},4,6,8},
}&\nonumber
\end{align}
\begin{align}
 \cal {R}(U)=\cal {R}(\std U)=32323312\equiv
 \YT{0.16in}{}{
{,,,1},
{,2,2,2},
{{3},3,3,3},
}\equiv 3^42^31\equiv
\YT{0.16in}{}{
{1,2,2,3},
{2,3,3},
{{3}},
}=Y_{(4^3-(4,3,1))}\nonumber
 \end{align}

\end{ex}
If $\mu\subseteq \lambda\subseteq \gamma$ are partitions, we say that the shape of $\gamma/\lambda$ extends the shape of $\lambda/\mu$ which in turn extends the shape $\mu$. In general, when we say that the SSYT  $V$ extends the SSYT $U$ it is meant that the shape of $V$ extends the shape of $U$ and we write $U\cup V$ for the object formed by gluing $U$ and $V$ together.

 If $S$ and $T$ are
SSYT's of shapes $\mu$ and $\lambda/\mu$,
 $S\cup T$ is said to be of shape $\lambda$. If $Z$ is another SSYT of shape $\gamma/\lambda$, $T\cup Z$ has shape $\lambda/\mu$ and we write $S\cup T\cup Z:=(S\cup T)\cup Z=S\cup (T\cup Z)$.
The tableaux $S$, $T$ and $Z$ are filled in any finite ordered alphabet consisting of positive integers. Sometimes it is convenient to look at  $S\cup T$ as a SSYT of shape $\lambda$ by adding a constant to all entries of $T$, for instance, the biggest entry of $S$. In view of these definitions if $U$ is a tableau we say the subtableaux $U_1$ and $U_2$ decompose $U$ if $U=U_1\cup U_2$ and $U_2$ extends $U_1$.
When $S$ and $T$ are both ballot tableaux with $S$ of normal shape, we say that we have a \emph{ballot tableau pair} of normal shape. In this case $S$ is the Yamanouchi tableau  $Y_\mu$.
\begin{defi}Let ${\cal{LR}}^{(n)}$ denote the set of all ballot semistandard tableau pairs of partition shape where the length of the shape is $\le n$.
\end{defi}

\section {Sagan-Stanley skew RSK and internal row insertion operators} \label{sec:innercorner}
\subsection{Internal insertion words of a skew tableau and their Knuth classes}
Let $T\in YT(\lambda/\mu) $ with $\ell(\lambda)=n$.
An \emph{inner corner} of $T$ is a cell $(i, j)$ such that when added to the Young diagram of $\mu$ still results in a valid Young diagram.
For instance, $(1,4)$, $(2,3)$,  $(4,2)$  and $(5,1)$ are inner corners of $T$ but $(6,1)$ and $(3,3)$ are not,
$T=\YT{0.16in}{}{
{,,},
{,,{1}},
{{},},
{{}},
{3},}$ with $\ell(\lambda)=5$.
 \begin{defi}\label{def:operator}~\cite{ss} Let $T\in YT(\lambda/\mu) $ with $\ell(\lambda)=n$ and $1\le i\le n+1$.  The Sagan-Stanley {\it  internal (row) insertion operator} $\phi_{i}$  is an operation over $ T$ defined whenever $(i,\mu_i+1)$ is an inner corner of $T$. We have to distinguish two cases.
\begin{itemize}
   \item The cell \emph{$(i,\mu_i+1)$ is an empty cell}  which means that   $(i - 1, \mu_i+1)   \not\in \lambda/\mu$. Then $\phi_i$  just adjoins the blank cell $(i,\mu_i+1)$ to the Young diagram of $\mu$, in row $i$ of $T$.

  \item The cell \emph{ $(i,\mu_i+1)\in \lambda/\mu$}.    Then $\phi_i$ vacates the cell $(i,\mu_i+1)$ of $T$, bumps its entry    and inserts the bumped element,
using the usual Schensted row insertion rules, into row $i + 1$ of $T$.
The insertion
then continues in a normal fashion, ending with an element settling at the
end of some row $\le n+1$.
\end{itemize}
\end{defi}

$$\phi_4\phi_1\YT{0.16in}{}{
{,,},
{,,{1}},
{{},},
{3},}=\phi_4\YT{0.16in}{}{
{,,,},
{,,{1}},
{{},},
{3},}=\YT{0.16in}{}{
{,,,},
{,,{1}},
{{},},
{{}},
{3},}$$
 With an internal row insertion operation $\phi_i$ on $T$, no new entry is added to the tableau $T$. Instead  the skew-shape changes by adding one blank box at the end of row $i$ of the inner shape $\mu$, and, if it is the case,  one filled box is added to the outer shape $\lambda$. The new tableau $\phi_iT$ has shape $(\lambda+e_t)/(\mu+e_i)$ with $i\le t\le n+1$, and is Knuth equivalent to $T$.
In particular, if  $ T$ is a ballot tableau, then $\phi_i T$ is also ballot.
Whenever the internal row insertion operator $\phi_i$ is defined on $T$, it  can be easily extended to  the tableau pair $Y\cup T$   with $Y=Y_\mu$, by putting
\begin{eqnarray}\label{def:operatorbar}\bar\phi_i(Y\cup T):=\begin{cases} Y_{(\mu_1,\dots, \mu_i+1,\dots, \mu_n)}\cup \phi_i(T), \;\text{if $1\le i\le n$},\\
 Y_{(\mu_1,\dots, \mu_n,{1})}\cup \phi_{n+1}(T), \;\;\text{if $i=n+1$ and $\mu_{n}>0$}.
 \end{cases}
 \end{eqnarray}
 $$\bar\phi_4\phi_1\YT{0.16in}{}{
{{\color{red}1},{\color{red}1},{\color{red}1}},
{{\color{red}2},{\color{red}2},{1}},
{{\color{red}3},{\color{red}3}},
{3},
}
=\YT{0.16in}{}{
{{\color{red}1},{\color{red}1},{\color{red}1},{\color{red}1}},
{{\color{red}2},{\color{red}2},{1}},
{{\color{red}3},{\color{red}3}},
{{\color{red}3}},
{3},}
$$
If $Y\cup T\in{\cal{LR}}^{(n)}$,  $\bar\phi_i(Y\cup T)\in{\cal{LR}}^{(n+1)}$.

\begin{defi} Given $T\in YT(\alpha/\mu)$, an  internal insertion order word for $T$ is the companion word $\cal{R}(U)$ of any  skew tableau $U$  sharing the inner border with $T$, that is, $U\in  YT(\beta/\mu)$ with $\mu\subseteq \beta$. We say $U\in  YT(\beta/\mu)$ with $\mu\subseteq \beta$  is an internal insertion order tableau for $T$.
\end{defi}

Given an internal insertion order word 
$u_{|\gamma|}u_{|\gamma|-1}\cdots u_2u_1$ of $T\in YT(\alpha/\mu)$ of content $\gamma$,
one defines the corresponding \emph{internal insertion  operator}
   \begin{align}\phi_{u_{|\gamma|}u_{|\gamma|-1}\cdots u_2u_1}:=
   \phi_{u_{|\gamma|}}\circ\phi_{u_{|\gamma|-1}}\circ\cdots \circ\phi_{u_2}\circ\phi_{u_1}.\label{internalinsertionwordoperator} \end{align}
   This is well defined because $\phi_{u_{|\gamma|}u_{|\gamma|-1}\cdots u_2u_1}T=\phi_{\cal{R}({U})}T$ for some $U\in  YT(\beta/\mu)$. Then $P(T,U):=\phi_{\cal{R}({U})}T\in YT(\lambda/\beta)$.

 The  bijection below is  a special case of the Sagan-Stanley skew-RSK correspondence, Theorem 6.11 in \cite{ss}, when the matrix word $\pi=\emptyset$, and is denoted by $\bf SS$. In this situation the  skew-insertion procedure is reduced to the internal (Schensted) row insertion procedure. This correspondence, calculates a bijection between pairs of tableaux $(T,U)$ sharing the inner border and pairs of tableaux $(P,Q)$ sharing the outer border that  preserves the Knuth equivalence class,  $T\equiv P$ and $U\equiv Q$. Also  the outer border of $T$ equals the inner border of $Q$, and the outer border of of $U$ equals the inner border of $P$. (See  also \cite{rssw} for more details and properties.)
 \begin{thm} \cite[Theorem 6.1, Theorem 6.11]{ss}, \cite{rssw}\label{th:internal} (Sagan-Stanley internal row insertion correspondence.) Fix  partitions $\mu\subseteq \alpha,\beta$. There is a bijection,
\begin{eqnarray}\label{P}YT(\alpha/\mu)\times YT(\beta/\mu)&\longrightarrow &\bigcup_{\begin{smallmatrix}\lambda\nonumber\\
|\lambda|=|\alpha|+|\beta|-|\mu|
\end{smallmatrix}}YT(\lambda/\beta)\times YT(\lambda/\alpha)\\
(T,{U})\;\;\;\;\;&\overset{\mbf{SS}}\longrightarrow& \;\;\;\;\;\;\;\;\;\;\;\;(P,Q),
\end{eqnarray}
where $P\equiv T$ and $Q\equiv U$. The $P$-tableau, $P(T,{U}):=P$ is given by $\phi_{\cal{R}({U})}T$.
{ The $Q$-tableau $Q(T,{U}):=Q$ is the recording tableau.}
\end{thm}
\begin{cor}
\label{cor:internal}  Fix  partitions $\mu\subseteq \alpha,\beta$ and the above setting
\begin{enumerate}
\item  If $\beta$ is a rectangle shape with $\ell(\beta)\ge \ell(\mu)$. Then
 $P=\emptyset_\beta\cup \mathrm{rect}(T)$ where  $\mathrm{rect}(T)$ is the rectification of $T$.

\item If $\alpha=\mu$
 and $U\in YT(\beta/\mu)$,
   $P(\emptyset_{\mu}, U)=\emptyset_\beta$ equals to the  Young diagram of shape $\beta$ and $Q(\emptyset_{\mu}, U)=U$.
   \item If $\mu=\emptyset$, $P(T, U)=T$ and $\cal{R}(U)$ is a Yamanouchi word of weight $\beta$.

  \end{enumerate}
\end{cor}

    Let $U$ in the alphabet $[m]$. Set $T^{(0)} = T$ and $U^{(0)} = \emptyset_\lambda$.  For $j = 1, \dots, m,$ let $r_1^{(j)}\ge\cdots\ge  r^{(j)}_{k_j}$
 be the row coordinates of
all $j$-cells in the standard order of $U$ and define
$ T^{(j)} = \phi_{\cal R_{r_1^{(j)}}}\circ\cdots\circ  \phi_{\cal R_{r^{(j)}_{k_j}}}(T^{(j-1)})$. Then define $U^{(j)}$ adding to $U^{(j-1)}$ $j$-cells so that the external shape of $U^{(j)}$ matches that of
$T^{(j)}$. Finally set $P = T^{(m)}$ and $Q = U^{(m)}$.  See Example \ref{ex:ss} below.

The following exhibits a symmetry of the skew RSK map of tableaux proven in \cite{ss} not immediate from the definition, and the compatibility with  standardization. Schensted insertion commutes with standardisation. In addition the effect of each operator $\phi_{u_i}$ on the shape of $T$ will be the same as the effect on the shape of ${\rm std} T$
\begin{prop}

$(a)$ \cite[Theorem 3.3]{ss} If $\mathbf{SS}(P,Q) = (P',Q')$, then $\mathbf{SS}(Q, P) =
(Q', P')$.
That is,  like $P'$ and $P$, also the recording tableau $Q'$ is obtained from $Q$
following a number of internal insertions.

$(b)$ $\mathrm{std}\circ \mathbf{SS}(P,Q) = \mathbf{SS} \circ\mathrm{std}(P,Q)$.

\end{prop}


Unlike for the classical RSK correspondence, a detailed description of properties of Sagan
and Stanley's algorithm has proven to be more challenging to obtain. A recent attempt was undertaken by  Imammura-Mucciconi-Sasamoto \cite{imusa} by iterating the skew RSK correspondence, that is, the space of ordered pairs (P,Q) of tableaux of skew shape, which are built on repeated internal insertions in P whose locations are determined by Q


Given a tableau $T\in YT(\alpha/\mu)$, we would like to characterize the tableaux ${U}, {U}'\in YT(\beta/\mu)$ such that $P(T,{U})=P(T,{U}')$ in \eqref{P}. The theorem below shows that Knuth equivalence on the companion words $\cal{R}({U})$ and $\cal{R}({U}')$ provides  sufficient  conditions for the equality. Given a skew tableau $T$, similarly to the ordinary RSK, Knuth equivalent internal insertion order words of $T$ means a certain Knuth commutation of the corresponding internal insertion operators which gives rise to the same tableau $P$. Our theorem  partially answers our aim and the second part of  question $(3)$ of Section 9, in \cite{ss}.

\begin{thm} (Theorem  \ref{1} )\label{th:U} (Sagan-Stanley row internal insertion operators and Knuth relations.)
 Let $T\in YT(\alpha/\mu)$,  $U, U'\in YT(\beta/\mu)$ and $P(T,U)$ respectively $P(T,U')$ the  corresponding $P$-tableaux in the Sagan-Stanley internal insertion correspondence. Then

$(a)$ $U$ and its standardization $\mathrm{std}\, U$ have the same companion word, ${\mathcal{R}({U})}=\mathcal{R}(\mathrm{std}\, U)$.

$(b)$ $P(T,U)=\phi_{\mathcal{R}(U)}(T)$.

$(c)$ $P(T,U)=P(T,\mathrm{std}\,U)$
 and  internal row insertion commutes with standardization  
$\mathrm{std}(P(T,U))=$ $P(\mathrm{std}\, T, U)$ $=P(\mathrm{std}\, T,\mathrm{std}\, U).$

$(d)$ $ P(T,U)  =P(T,U')=\phi_{\mathcal{R}(U)}(T)=\phi_{\mathcal{R}(U')}(T)$ whenever $\mathcal{R}(U)\equiv \mathcal{R}(U')$  are Knuth equivalent.

\end{thm}

The proof will be delayed until Section \ref{sec: prooftheorem1} and will follow from   two lemmata, Lemma \ref{lem:knuth-intern2}, Lemma \ref{lem:knuth-intern3} and Proposition \ref{propp:knuth}.

\begin{ex}\label{ex:ss} Let $\mu=(1)$, $\alpha=(3,2,0)$ and $\beta=(4,2,1)$. Below
 $U$, $\std U$,  and $U'$, $\std(U')$ are internal insertion order tableaux  of
 $
T=\YT{0.16in}{}{
{,1,3},
{{2},3},
}$. The corresponding internal insertion order words are Knuth equivalent, $\cal R(U)=312211=\cal R(\std U)\equiv \cal R(U')=\cal R(\std U')=132211$, $$\begin{array}{cccccccccc}
 U=\YT{0.16in}{}{
{,1,1,2},
{{2},2},
{3},}&
\std U=\YT{0.16in}{}{
{,1,2,5},
{{3},4},
{6},}
&U'=\YT{0.16in}{}{
{,1,1,3},
{{2},2},
{3},}&\std U'=\YT{0.16in}{}{
{,1,2,6},
{{3},4},
{5},}
\end{array}.$$
Knuth equivalence of internal insertion order words, $\cal R(U)\equiv \cal R(U')$, implies that the corresponding internal insertion operators $\phi_{\cal R(U)}$ and $ \phi_{\cal R(U')}$  commute according to the Knuth relations, that is, $\phi_{\cal R(U)}=\phi_{\cal R(U')}$.

 \begin{eqnarray}\label{stand1}
P(T,U)=\phi_3\phi_{122}\phi_{11}T=\phi_{13}\phi_{22}\phi_{11}T=P(T,U')
=\YT{0.16in}{}{
{,,,},
{,,3},
{,3},
{1},
{2},
}.
\nonumber\end{eqnarray}

\begin{align}Q(T,U)=\YT{0.16in}{}{
{,,,2},
{,,1},
{1,2},
{2},
{3},
}&
Q(T,U')=\YT{0.16in}{}{
{,,,3},
{,,1},
{1,2},
{2},
{3},
}.
\end{align}

$\cal R(U)\equiv \cal R(U')$ and $\cal R(Q(T,U))=513423\equiv \cal R(Q(T,U'))=153423.$

\medskip
\begin{ex}\label{ex:sspq}
For the $T$ and $U$ above one has $n=3$, $\alpha=(3,2,0)$, $\beta=(4,2,1)$, $\lambda=(4,3,2,1,1)$, $(P,Q)\in YT((4,3,2,1,1)/(4,2,1))\times YT((4,3,2,1,1)/(3,2,0))$
\begin{eqnarray}YT(\alpha/(1),3)\times YT(\beta/(1),3)&\longrightarrow &\bigcup_{\begin{smallmatrix}\lambda\nonumber\\
|\lambda|=11
\end{smallmatrix}}YT(\lambda/\beta,3)\times YT(\lambda/\alpha,3)\\
(T,{U})\;\;\;\;\;&\overset{\mbf{SS}}\longrightarrow& \;\;\;\;\;\;\;\;\;\;\;\;(P,Q).
\end{eqnarray}

\begin{align*}
&T^{(0)}=\YT{0.16in}{}{
{,1,3},
{{2},3},
}
\quad U^{(0)}=\emptyset_\alpha=\YT{0.16in}{}{
{,,},
{{},},
}\\
&\overset{\phi_{\cal{R}_1}=\phi_{11}}\longrightarrow
T^{(1)}=\YT{0.16in}{}{
{,,},
{{1},3,3},
{2}
}
\quad U^{(1)}=\YT{0.16in}{}{
{,,},
{{},,1},
{1}
}\\
&\overset{\phi_{\cal{R}_2}=\phi_{122}}\longrightarrow T^{(2)}=\YT{0.16in}{}{
{,,,},
{{},,3},
{1,3},
{2}
}
\quad U^{(2)}=\YT{0.16in}{}{
{,,,2},
{{},,1},
{1,2},
{2}
}\\
&
\overset{\phi_{\cal{R}_3}=\phi_{3}}\longrightarrow P=T^{(3)}=\YT{0.16in}{}{
{,,,},
{{},,3},
{,3},
{1},
{2}
}
\quad Q=U^{(3)}=\YT{0.16in}{}{
{,,,2},
{{},,1},
{1,2},
{2},
{3}
}\\
&{\mbf{SS}}(T,U)=(P,Q)
\end{align*}
\end{ex}

The inverse Sagan-Stanley internal insertion correspondence $\mbf{SS^{-1}}$, $\alpha=(3,2,0)$, $\beta=(4,2,1)$:

\begin{eqnarray}\label{rskinverse}YT(\lambda/\beta,3)\times YT(\lambda/\alpha,3)&\longrightarrow &\bigcup_{\begin{smallmatrix}\lambda\nonumber\\
|\lambda|=11
\end{smallmatrix}}YT(\alpha/(1),3)\times YT(\beta/(1),3)\\
(P,{Q})\;\;\;\;\;&\overset{\mbf{SS^{-1}}}\longrightarrow& \;\;\;\;\;\;\;\;\;\;\;\;(T,U).
\end{eqnarray}

 The inverse Sagan-Stanley {\it  internal (row) insertion operator}, $\phi_i^{-1}$, or Sagan-Stanley {\it  deletion  operator}, is denoted  $\Delta_{i}$.
The \emph{reverse companion word} of $Q$, $rev\cal{R}(Q):={rev\cal{R}_1}{rev\cal{R}_2}{rev\cal{R}_3}=32\,431\,5$ defines the  \emph{deletion operator}
\begin{align}\label{deletionoperator}\Delta_{rev\cal{R}(Q)}&=\Delta_{{rev\cal{R}_1}{rev\cal{R}_2}{rev\cal{R}_3}}\\
&=\Delta_{rev\cal{R}_1}\circ \Delta_{rev\cal{R}_2}\circ \Delta_{rev\cal{R}_3}
\end{align}

\begin{align*}
P^{(0)}=
\YT{0.16in}{}{
{,,,},
{{},,3},
{,3},
{1},
{2}
}
\quad Q^{(0)}=\emptyset_\beta=\YT{0.16in}{}{
{,,,},
{{},},
{}
}
\overset{\Delta_{\cal{R}_3}=\Delta_{5}}\longrightarrow
P^{(1)}=
\YT{0.16in}{}{
{,,,},
{{},,3},
{1,3},
{2},
}
\quad Q^{(1)}=\YT{0.16in}{}{
{,,,},
{{},},
{3}
}\\
\overset{\Delta_{\cal{R}_2}=\Delta_{431}}\longrightarrow
P^{(2)}=
\YT{0.16in}{}{
{,,},
{{1},3,3},
{2},
}
\quad Q^{(2)}=\YT{0.16in}{}{
{,,,2},
{{2},2},
{3}
}
\overset{\Delta_{\cal{R}_1}=\Delta_{32}}\longrightarrow
P^{(3)}=
\YT{0.16in}{}{
{,1,3},
{{2},3},
}=T
\quad Q^{(3)}=\YT{0.16in}{}{
{,1,1,2},
{{2},2},
{3}
}=U
\\
\end{align*}
$$\mbf{SS^{-1}}(P,Q)=(T,U)$$
\end{ex}

\section{A preserver for the $P$-tableau in the skew RSK correspondence}\label{sec:preserver}


The Schensted row insertion  takes a SSYT  $T$ of partition shape,
and an element  $m$ in the $T$-alphabet, and constructs a new tableau, denoted $T\bigcdot m$ \cite{fulton}.
For skew tableaux there are two types of row insertion: \emph{external} and
\emph{internal} both of which based on Schensted insertion but with different procedures namely, in the former, the element in the $T$-alphabet to be inserted is  added to the multiset  $\{T\}$, and, in the latter,   is picked in $\{T\}$\cite[Section 2]{ss}.

\subsection{External row insertion on skew tableaux} \emph{Sagan-Stanley row external insertion} \cite[Section 2]{ss} on skew-tableaux is similar to Schensted's original procedure.
We start with a SSYT $T$ of shape $\lambda/\mu$ and an element $m$ in the $T$-alphabet to be  added to $T$.  To start, $m$ replaces the smallest entry in the first row of $T$ strictly larger than $ m$; in the
case where $m$ is bigger or equal than all entries in the first row, it is placed at that row's right end. If an entry was
displaced from the first row then it is inserted into the second using the
same rules as above. This process continues until some element comes to
rest at the end of a row.
The only need for caution in the skew case is when
something is to be inserted into row $i$ which is empty  and this can only happen
at the beginning or at the end of an insertion. It happens, when
$\mu_i=\lambda_i$  and in this case we put the element $m$ in cell $(i, \lambda_i+1)$. We obtain a new skew-tableau denoted $T.m$ with the same inner shape $\mu$ and $\{T\bigcdot m\}=\{T\}\cup \{m\}$.

If $T$ is a skew-tableau 
and $u$ is a word over the $T$-alphabet, $T.u$ denotes the skew tableau with the same inner shape as $T$ obtained by the \emph{Sagan-Stanley row external insertion} of $u$ in $T$. That is, if $u=u_1\cdots u_n$, $T\bigcdot u=(((T\bigcdot u_1)\bigcdot \cdots)\bigcdot u_n)$.
When $T$ is of straight shape $T\bigcdot u$  coincides with the usual Schensted insertion \cite{fulton}.
More generally, if $U$, $V$ and $W$ are tableaux of straight shape, $U.V$ means $U\bigcdot  w(V)$, and $(U\bigcdot V)\bigcdot W =U\bigcdot (V\bigcdot W)$.

We may look at Sagan-Stanley external insertion on a skew-tableau $T$ of shape $\lambda/\mu$ as an insertion on a straight tableau $Y_\mu\cup T$ of shape $\lambda$, where $\ell(\mu) $ is added to each entry of $T$, that is, $T$ is filled in the alphabet $\{\ell(\mu)+1,\ell(\mu)+2,\dots\}$.  If $ u$  is a  word over the alphabet $\{\ell(\mu)+1,\ell(\mu)+2,\dots\}$, then
\begin{equation}\label{straight}Y_\mu\cup (T\bigcdot u)=(Y_\mu\cup T)\bigcdot u.\end{equation}

 The following is a lemma on commutation
and cancellation in the plactic monoid of Lascoux and Sch\"tzenberger
\begin{lem}\label{lem:congr}  Let  $u$ and $u'$ be two words in  a same ordered alphabet $A=[n]$.
\begin{enumerate}
\item [(a)]\cite[Proposition 2.3]{ls88} 
Let   $B\subseteq A$ be  an interval and let  $u_{|B}$ and let $u'_{|B} $ be the restrictions of $u$ and $u'$  to $B$. Then one has
\begin{equation}\label{ls}u\equiv u'\Rightarrow u_{|B}\equiv u'_{|B}.\end{equation}

\item [(b)] \cite[Lemma 7.5]{rssw} Let $w$, $v$ be any words in the alphabet $[k]$ and $C=k\cdots 21$ be a column. Then

(i) $C\bigcdot w\equiv w.C$.

(ii) $C\bigcdot w\equiv C\bigcdot v$ if and only if $w \equiv v$.

\item [(c)]
If $ux\equiv u'x$ with $x=a_1\cdots a_k$ a row word in $A$  such that the smallest letter $a_1$ in  $x$ is bigger or equal than all letters in the words $u$ or $u'$ then

    (i) $ux\equiv u'x\equiv T x$  is  a tableau for some tableau $T$ with entries in $A$ less or equal than $a_1$ such that
 $ u\equiv u'\equiv T$.
    \end{enumerate}
 \end{lem}

\begin{lem}\label{skewknuth} (External insertion on skew-tableaux.)
Let $T$ be a skew-tableau  and $u$ and $ v$ two  words over the $T$-alphabet.
Then
\begin{enumerate}
\item[(a)] $w(T.u)\equiv w(T)u$, and
\item[(b)] $T\bigcdot u=T\bigcdot v$  whenever $u\equiv v$.
\end{enumerate}
\end{lem}
\begin{proof}
(a) 
If $T$ has straight shape then $w(T\bigcdot u)\equiv w(T)u$. In general, if $T$ has shape $\lambda/\mu$, using \eqref{straight}, one has $w[Y_\mu\cup (T\bigcdot u)]=w[(Y_\mu\cup T)\bigcdot u]\equiv w(Y_\mu\cup T)u$. From \eqref{ls}, it follows $w(T\bigcdot u)\equiv w(T)u$.

(b) Again using \eqref{straight}, as a consequence of the usual Schensted insertion on straight shapes, one has $Y_\mu\cup (T\bigcdot v)=(Y_\mu\cup T)\bigcdot v=(Y_\mu\cup T)\bigcdot u=Y_\mu\cup (T\bigcdot u)$ whenever $u\equiv v$ over the $T$-alphabet. Thus $Y_\mu\cup (T\bigcdot v)=Y_\mu\cup (T.u)$  whenever $u\equiv v$ over the $T$-alphabet implies $T\bigcdot u=T\bigcdot v$  whenever $u\equiv v$ over the $T$-alphabet.
\end{proof}
The reverse direction of Lemma \ref{skewknuth}, $(b)$, that is, the cancelation law, is trivially true when $T$ is just a Young diagram. Otherwise, unless additional conditions are satisfied as in Lemma \ref{lem:congr}, it is false in general. For instance, one has $yxz\equiv yzx$ for $x< y\le z$ and $xz\not\equiv zx$ with $x<z$. 

Indeed external insertion preserves Knuth equivalence: if $T$ and $Z$ are tableaux with the same inner shape, and $T\equiv Z$ and $u$ and $v$ are Knuth equivalent words, then from Lemma \ref{skewknuth}, $T\bigcdot u\equiv Z\bigcdot v$. In case $T$ and $Z$ are of the same straight shape  then we have $T=Z$ and  $T\bigcdot u=T\bigcdot v$.

\subsection{Internal   row insertion and internal bumping routes}



   \begin{defi}\label{def:restrictiontab} If $T\in YT(\lambda/\mu)$, with $\ell(\lambda)\le n$, $T^{[i]}\in YT((\lambda_1,\dots,\lambda_i)/(\mu_1,\dots,\mu_i))$ denotes the tableau  consisting of the first $i$ rows of $T$, for $i=0,1,\dots,n$, that is, the restriction of the tableau $T$ to the first $i$ rows. We put $T^{[0]}=\emptyset$ and $T^{[n]}=T$.
\end{defi}

\begin{defi}\label{def:tabfactor} Given the skew tableau $T\in YT(\lambda/\mu)$ and  $0\le i\le n$ where $\ell(\lambda)\le n$, we say that  $T$ is factorized across the row $i$ when we write   $T=T'\ast T^{[i]}$, $0\le i\le n$,
~with  $T'$  of shape $(\lambda_{i+1},\dots,\lambda_n)/(\mu_{i+1},\dots,\mu_n)$ is the restriction of $T$ to the last $n-i$ rows.
\end{defi}

   In particular, let  $T$ be a SSYT of straight shape  factored across the  $i$th row,  $T=T'\ast T^{[i]}$, with  $T'$ the tableau obtained by suppressing $T^{[i]}$ from $T$. Note  that $T=T'\ast T^{[i]}=T'.T^{[i]}$.

  Let $T\in YT(\lambda/\mu) $ with $\ell(\lambda)=n$, assume ${\color{red}1}<{\color{red}2}<\cdots< {\color{red}n}<1<2<\cdots<n$ and fill in the inner shape $\emptyset_\mu$   with the coloured alphabet  ${\color{red}1}<{\color{red}2}<\cdots< {\color{red}n}$ so that we get the Yamanouchi tableau $Y_\mu$. Let us factor  $Y_\mu\cup T$ across the $i$th row. Then one has,

  \begin{align}Y_\mu\cup T=(Y'\cup T')\ast(Y_\mu\cup T)^{[i]}=(Y'\cup T')\ast(Y_{(\mu_1,\dots,\mu_{i})}\cup T^{[i]})=(Y'\cup T')\bigcdot(Y_{(\mu_1,\dots,\mu_{i})}\cup T^{[i]})\label{starinsertion}\end{align}
  where $Y'$ and $T'$ denote
the restriction of $Y_\mu$ and $T$ to the last $n-i$ rows of $Y_\mu$ and $T$ respectively. Thereby, in the sense of \eqref{starinsertion}, we may write
$$T=T'\ast T^{[i]}=T'\bigcdot T^{[i]}.$$

Assume that $(i,\mu_i+1)$ is an inner corner of  $ T^{[i]}$ with entry $x$. It then follows that the action of the internal row insertion operator $\phi_i$ on $ T$ \eqref{def:operator} may be read as  an operation which bumps the entry $x$ in the cell $(i,\mu_i+1)$ of  $ T^{[i]}$ and then inserts externally the bumped element $x$ in the subtableau $ T'$.  That is,  $\phi_i$ on $ T^{[i]}$ bumps the entry $x$ and left justifies it in the $(i+1)$th row.
Then
 \begin{eqnarray}
 \phi_i (T^{[i]})&=&x\ast T^{[i]}_{-}\label{star1}\\
 \bar\phi_i(Y_{(\mu_1,\dots,\mu_{i})}\cup T^{[i]})&=&x\ast(Y_{(\mu_1,\dots,\mu_{i}+1)}\cup T^{[i]}_{-}),\label{star2}\end{eqnarray}
 where $T^{[i]}_{-}$ is $T^{[i]}$ with  the entry $x$ bumped out from the cell $(i,\mu_i+1)$. It means that $T^{[i]}_{-}$ $\in $ \linebreak $YT((\lambda_1,\dots, \lambda_i)/(\mu_1,\dots,\mu_i+1))$, that is,  $T^{[i]}$ with the left most entry of the $i$th row suppressed and the corresponding blank cell  added to the $i$th row of the inner shape of $T^{[i]}$.  Henceforth, from \eqref{starinsertion}, \eqref{star1} and \eqref{star2},
 \begin{eqnarray}\label{fact}\phi_i (T)&=&\phi_i[ T'\ast T^{[i]}] =T'.\phi_i( T^{[i]})=  T'\bigcdot[x\ast T^{[i]}_{-}]=
 (T'\bigcdot x)\ast T^{[i]}_{-},\\
 \bar\phi_i(Y_\mu\cup T)&=&\bar\phi_i[(Y'\cup T')\ast(Y_{(\mu_1,\dots,\mu_{i})}\cup T^{[i]})]=(Y'\cup T')\bigcdot \bar\phi_i(Y_{(\mu_1,\dots,\mu_{i})}\cup T^{[i]})
 \nonumber\\
&=&(Y'\cup T')\bigcdot [x\ast(Y_{(\mu_1,\dots,\mu_{i}+1)}\cup T^{[i]}_{-})]\nonumber\\
&=&[(Y'\cup T')\bigcdot x]\ast(Y_{(\mu_1,\dots,\mu_{i}+1)}\cup T^{[i]}_{-})\nonumber\\
&=&(Y'\cup T')\bigcdot x\ast(Y_{(\mu_1,\dots,\mu_{i}+1)}\cup T^{[i]}_{-}).\label{factorization}
\end{eqnarray}
If $\mu_i=\lambda_i$ and $(i,\lambda_i+1)$ is a blank  inner corner then $x=\emptyset$ and $\phi_i T^{[i]}=\emptyset\ast T^{[i]}_{-}=T^{[i]}_{-}$, with
$T^{[i]}_{-}\in YT((\lambda_1,\dots, \lambda_i+1)/(\mu_1,\dots,\mu_i+1))$ a. Also
$\phi_i T=T'\ast T^{[i]}_{-}$.

As an aside, observe that $w(\phi_i T^{[i]})=w(x\ast T^{[i]}_{-})=x\,w(T^{[i]}_{-})=w(T^{[i]})$, and obviously $\phi_i T^{[i]}\equiv  T^{[i]}$.
Therefore, using Lemma \ref{skewknuth},
\begin{eqnarray}\label{internalknuth}
w(\phi_i T)&=&w( T'\bigcdot x\ast T^{[i]}_{-})= w( T'\bigcdot x)w( T^{[i]}_{-})\equiv  w( T')x w( T^{[i]}_{-})= w(T')w(\phi_i T^{[i]})\nonumber\\&=&w(T')w(T^{[i]})=w(T'\ast T^{[i]})=w(T).
 \end{eqnarray}

 \begin{obs}\label{ob:factorLR} Let  $Y_{(\mu_1,\dots,\mu_{n+1})}\cup T\in\cal{LR}^{(n+1)}$ be a ballot tableau pair and $1\le m\le n$. Consider the factorisation through row $m$,
  $$Y_{(\mu_1,\dots,\mu_{n+1})}\cup T=[Y_{(\mu_{m+1},\dots,\mu_{n+1})}\cup \widehat T]\ast[Y_{(\mu_1,\dots,\mu_{m})}\cup T^{[m]}]\in{\cal{LR}}^{(n+1)}.$$ Then  ${T^{[m]}}$ is a ballot tableau on the alphabet $[m]$, and
 $\widehat T\in YT((\lambda_{m+1},\dots,\lambda_{n+1})/(\mu_{m+1},\dots,\mu_{n+1}))$ consisting of the last $n+1-m$ rows of $T$ is such that  the word  restricted to the alphabet $[m+1,n+1]$ satisfies the Yamanouchi condition.
 \end{obs}

\begin{ex}\label{ex:factor4} An illustration of \eqref{fact} and \eqref{factorization} is given below:
$$\begin{array}{ccccccccccccccc}
\bar\phi_3
\YT{0.16in}{}{
{{\color{red}1},{\color{red}1},{\color{red}1},{\color{red}1},1,1},
{{\color{red}2},{\color{red}2},{1},{2},2},
{{\color{red}3},{2},3,3},
{{\color{red}4},3,4,4},
{5,5},
}
=
\YT{0.16in}{}{
{{\color{red}4},3,4,4},
{5,5},
}
\,\bigcdot\,
\bar\phi_3
\YT{0.16in}{}{
{{\color{red}1},{\color{red}1},{\color{red}1},{\color{red}1},1,1},
{{\color{red}2},{\color{red}2},{1},{2},2},
{{\color{red}3},{2},3,3},
}
=
\YT{0.16in}{}{
{{\color{red}4},3,4,4},
{5,5},
}
\,\bigcdot\,
\YT{0.16in}{}{
{{\color{red}1},{\color{red}1},{\color{red}1},{\color{red}1},1,1},
{{\color{red}2},{\color{red}2},{1},{2},2},
{{\color{red}3},$\textcircled{\color{red}3}$,3,3},
{$\textcircled 2$},
}\\
\end{array}$$
$$
\begin{array}{cccccc}
=
\YT{0.16in}{}{
{{\color{red}4},3,4,4},
{5,5},
}
\,\bigcdot\,
\YT{0.16in}{}{
{$\textcircled 2$},}
\,\ast\,
\YT{0.16in}{}{
{{\color{red}1},{\color{red}1},{\color{red}1},{\color{red}1},1,1},
{{\color{red}2},{\color{red}2},{1},{2},2},
{{\color{red}3},{\color{red}$\textcircled {\color{red}3}$},3,3},
}\qquad
\end{array}$$
$$\begin{array}{ccccccccccccccc}
&&&&=
\YT{0.16in}{}{
{{\color{red}4},$\textcircled 2$,4,4},
{3,5},
{5},
}
\,\ast\,
\YT{0.16in}{}{
{{\color{red}1},{\color{red}1},{\color{red}1},{\color{red}1},1,1},
{{\color{red}2},{\color{red}2},{1},{2},2},
{{\color{red}3},$\textcircled {\color{red}3}$,3,3},
}
=\YT{0.16in}{}{
{{\color{red}4},$\textcircled 2$,4,4},
{3,5},
{5},
}
\,\bigcdot \,
\YT{0.16in}{}{
{{\color{red}1},{\color{red}1},{\color{red}1},{\color{red}1},1,1},
{{\color{red}2},{\color{red}2},{1},{2},2},
{{\color{red}3},$\textcircled {\color{red}3}$,3,3},
}\\
\\
&&&&=\YT{0.16in}{}{
{{\color{red}4},$\textcircled 2$,4,4},
{3,5},
{5},
}
\,\bigcdot \,\YT{0.16in}{}{
{{\color{red}3},$\textcircled {\color{red}3}$,3,3},
}\bigcdot
\YT{0.16in}{}{
{{\color{red}1},{\color{red}1},{\color{red}1},{\color{red}1},1,1},
{{\color{red}2},{\color{red}2},{1},{2},2},
}\\
\\
&&&&=
\YT{0.16in}{}{
{{\color{red}3},{\color{red}3},3,3},
{{\color{red}4},2,4,4},
{3,5},
{5},
}\bigcdot
\YT{0.16in}{}{
{{\color{red}1},{\color{red}1},{\color{red}1},{\color{red}1},1,1},
{{\color{red}2},{\color{red}2},{1},{2},2},
}=\YT{0.16in}{}{
{{\color{red}1},{\color{red}1},{\color{red}1},{\color{red}1},1,1},
{{\color{red}2},{\color{red}2},{1},{2},2},
{{\color{red}3},{\color{red}3},3,3},
{{\color{red}4},2,4,4},
{3,5},
{5},
}.
\end{array}$$
\end{ex}
An {\it internal row insertion operator} $\bar\phi_{i}$ ($\phi_i$) on $Y\cup T$ ($T$)
determines a collection $R_i$ of boxes, which are those where an element is
bumped from a row, together with the box where the last bumped element
lands and settles \cite{fulton}. Let us call to $R_i$ the  $\bar\phi_{i}$-{\em bumping route} of  $Y\cup T$. In particular, retain that whenever   $\bar\phi_{i}$ acts on $Y\cup T$ and the $\bar\phi_{i}$-{\em bumping route} terminates in some row $k> i$, this means that the $k$th row of $\bar\phi_{i}(Y\cup T)$  equals  the $k$th row  of $Y\cup T$ with  the last bumped entry in $\{T\}$  added at the end.
  If $\ell(\lambda)\le n$, the $\phi_i$-bumping path (route) terminates in some row $\le n+1$.

  For instance, below $Y\cup T$, $Y\cup H\in {\cal{LR}}^{(5)}$, and $\bar\phi_3$ acting on $Y\cup T$  bumps $2$ and fills the vacant cell with $\color{red}3$, then  the bumped $2$ is inserted in the $4$th row and bumps $3$, $3$ bumps $5$ which lands in the $6$th row, and $5$ is added to the $6$th row. In \eqref{route},  $R_3$ consists of the black boxes and the bumped numbers are $2$, $3$ and $5$, highlighted with circles,  the entries of the three last boxes of $R_3$. The $\bar\phi_1$-bumping routes $R_1$ in $Y\cup T$ and $Y\cup H$ are similarly displayed below in \eqref{route}:
\begin{eqnarray}\label{route1}
Y\cup T=\YT{0.145in}{}{
{{\color{red}1},{\color{red}1},{\color{red}1},{\color{red}1},1,1},
{{\color{red}2},{\color{red}2},{1},{2},2},
{{\color{red}3},{2},3,3},
{{1},3,4,4},
{5,5},
}&
\bar\phi_3(Y\cup T)=\YT{0.15in}{}{
{{\color{red}1},{\color{red}1},{\color{red}1},{\color{red}1},1,1},
{{\color{red}2},{\color{red}2},{1},{2},2},
{{\color{red}3},$\textcircled{\color{red}3}$,{3},3},
{1,$\textcircled 2$,4,4},
{$\textcircled 3$, 5},
{$\textcircled 5$},
}&
\bar\phi_1(Y\cup T)=\YT{0.15in}{}{
{{\color{red}1},{\color{red}1},{\color{red}1},{\color{red}1},$\textcircled{\color{red}1}$,1},
{{\color{red}2},{\color{red}2},1,$\textcircled 1$,2},
{{\color{red}3},{2},$\textcircled 2$,3},
{1,3,$\textcircled 3$,4 },
{$\textcircled 4$,5},
{$\textcircled 5$},
}
\end{eqnarray}
\begin{eqnarray}
Y\cup H=\YT{0.16in}{}{
{{\color{red}1},{\color{red}1},{\color{red}1},{\color{red}1},1,1},
{{\color{red}2},{\color{red}2},{1},{2},2},
{{\color{red}3},{2},3,3},
{{1},3},
{4,4},
}&\bar\phi_1(Y\cup H)=\YT{0.16in}{}{
{{\color{red}1},{\color{red}1},{\color{red}1},{\color{red}1},$\textcircled{\color{red}1}$,1},
{{\color{red}2},{\color{red}2},{1},$\textcircled{1}$,2},
{{\color{red}3},{2},$\textcircled{2}$,3},
{1, 3,$\textcircled 3$},
{4, 4},
}.\nonumber
\end{eqnarray}

The $\bar\phi_3$-bumping routes $R_3$, $\bar\phi_1$-bumping route $R_1$ of $Y\cup T$, and $\bar\phi_1$-bumping route of $Y\cup H$
\begin{align}\label{route}
\YT{0.16in}{}{
{{\color{red}1},{\color{red}1},{\color{red}1},{\color{red}1},1,1},
{{\color{red}2},{\color{red}2},{1},{2},2},
{{\color{red}3},\blacksquare,3,3},
{{1},\blacksquare,4,4},
{\blacksquare,5},
{\blacksquare},
}
&
\YT{0.16in}{}{
{{\color{red}1},{\color{red}1},{\color{red}1},{\color{red}1},{\blacksquare},1},
{{\color{red}2},{\color{red}2},1,\blacksquare,2},
{{\color{red}3},{2},\blacksquare,3},
{1,3,\blacksquare,4 },
{\blacksquare,5},
{\blacksquare},
}&
\YT{0.16in}{}{
{{\color{red}1},{\color{red}1},{\color{red}1},{\color{red}1},\blacksquare,1},
{{\color{red}2},{\color{red}2},{1},\blacksquare,2},
{{\color{red}3},{2},\blacksquare,3},
{1, 3,\blacksquare},
{4, 4},
}.
\end{align}

One often says  $a$ is a $\bar\phi_i$-bumped element of the skew-tableau $T$ to mean an  entry $a$ of $T$ that is bumped under the action of $\bar\phi_i$ on $T$.
 Thanks to \eqref{fact}, one has  as a consequence of the Row External Bumping Lemma in Section 1.1 of \cite{fulton}. This Lemma is instrumental in the proof of the Main Theorem and clearly shows that the internal insertion operators do not obey a naive commutation but instead a Knuth relation commutation as we shall see in the next section.

 \begin{lem} (\cite{fulton} Row internal bumping routes) \label{lem:1} Consider $1\le i\le j\le n$. Let $R_j$ and $R'_i$ be the pair of bumping routes of $\bar\phi_j$ on $Y\cup T$ and $\bar\phi_i$ on $\bar\phi_j(Y\cup T)$ respectively; and let $R_i$ and $R'_j$ be the pair of bumping routes of $\bar\phi_i$ on $Y\cup T$ and $\bar\phi_j$ on $\bar\phi_i(Y\cup T)$ respectively. Let $B$ and $B'$ be the corresponding pair of new boxes. Then it holds
\begin{enumerate}
 \item[(a)]  $R_j $ is strictly left of $R_i'$ and $B $ is strictly left of and
weakly below $B'$:
$$\begin{array}{cccc}\begin{smallmatrix}B&B'\end{smallmatrix}& \text{or} &
\begin{smallmatrix}
&B'\\
B&\\
\end{smallmatrix}\end{array}.$$

\item[ (b)] $R'_j$ is weakly left of $R_i$ and $B'$ is weakly left of and
strictly below $B$:

\begin{equation}\label{below}\begin{smallmatrix}
B\\
B'\\
\end{smallmatrix} \quad\text{or}\quad \begin{smallmatrix}
&B\\
B'&\\
\end{smallmatrix}.
\end{equation}
\end{enumerate}
In particular, $R'_j$ goes always strictly below the bottom box $B$ of $R_i$ by $\bar\phi_j$-bumping the element in $B$, this necessarily happens in  the case of the left hand side  of \eqref{below}, or by passing strictly to the left of $B$. Henceforth, if $x$ is $\bar\phi_i$-bumped  and $y$ is $\bar\phi_j$-bumped from the same row then $y<x$.
Moreover, if $B$ was created in the $(n+1)$th row then one has
$\begin{smallmatrix}
B\\
B'\\
\end{smallmatrix}$ and the last $\bar\phi_i$-bumped element resting in $B$ is $\bar\phi_j$-bumped out to be  settled in $B'$ and is strictly bigger than the element $\bar\phi_j$-inserted in $B$.
 \end{lem}
\begin{ex} \label{ex:route} An illustration of the internal insertion bumping routes lemma:
 \begin{align}\label{route23}
Y\cup T=\YT{0.16in}{}{
{{\color{red}1},{\color{red}1},{\color{red}1},1,1},
{{\color{red}2},{\color{red}2},{\color{red}2},2},
{{\color{red}3},{\color{red}3},2,3},
{{\color{red}4},1,3,4},
{1,4,5},
{2}
}&
\bar\phi_3\bar\phi_1(Y\cup T)=\bar\phi_3\YT{0.16in}{}{
{{\color{red}1},{\color{red}1},{\color{red}1},{\blacksquare},1},
{{\color{red}2},{\color{red}2},{\color{red}2},$\textcircled 1$},
{{\color{red}3},{\color{red}3},2,$\textcircled 2$},
{{\color{red}4},1,3,$\textcircled 3$},
{1,4,$\textcircled 4$},
{2,$\textcircled 5$},
}=
\YT{0.16in}{}{
{{\color{red}1},{\color{red}1},{\color{red}1},{\blacksquare},1},
{{\color{red}2},{\color{red}2},{\color{red}2},$\textcircled 1$},
{{\color{red}3},{\color{red}3},{\blacksquare},$\textcircled 2$},
{{\color{red}4},1,$\textcircled{ \color{blue}2}$,$\textcircled 3$},
{1,$\textcircled{\color{blue}3}$,$\textcircled 4$},
{2,$\textcircled{ \color{blue} 4}$},
{$\textcircled{ 5}$},
}
\end{align}
\begin{align}
\bar\phi_1\bar\phi_3(Y\cup T)=\bar\phi_1\YT{0.16in}{}{
{{\color{red}1},{\color{red}1},{\color{red}1},1,1},
{{\color{red}2},{\color{red}2},{\color{red}2}, 2},
{{\color{red}3},{\color{red}3},{\blacksquare}, 3},
{{\color{red}4},1,$\textcircled{ \color{blue} 2}$,4},
{1,$\textcircled{\color{blue}3}$, 5},
{2,$\textcircled{ \color{blue} 4}$},
}&=
\YT{0.16in}{}{
{{\color{red}1},{\color{red}1},{\color{red}1},{\blacksquare},1},
{{\color{red}2},{\color{red}2},{\color{red}2},$\textcircled 1$},
{{\color{red}3},{\color{red}3},{\blacksquare},$\textcircled 2$},
{{\color{red}4},1,$\textcircled{ \color{blue}2}$,$\textcircled 3$},
{1,$\textcircled{\color{blue}3}$,$\textcircled 4$},
{2,$\textcircled{ \color{blue} 4}$,$\textcircled 5$},
},
\end{align}
$$\bar\phi_1\bar\phi_3(Y\cup T)\neq \bar\phi_3\bar\phi_1(Y\cup T).$$
\end{ex}
\subsection{Proof of Theorem \ref{1}: 
Internal row insertion operators satisfy Knuth relations}
\label{sec: prooftheorem1}


In this section, for a fixed  $T$,  a sufficient condition on $U$ for the coincidence of $P(T,U)$ in the Sagan-Stanley internal row insertion bijection in Theorem \ref{th:internal} is provided. The sufficient condition does not involve the tableau $U$ directly but rather its companion word.  The companion word of $U$, Definition \ref{comp},
 encodes the inner corners for the action of the sequence of internal insertion operators acting on $T$.
Before giving the proof of Theorem \ref{1} (Theorem \ref{th:U}) we start with some warmup results.

When   the shape $\beta$ of ${U}$   is  a rectangle of height $\ge \ell(\mu)$ and width $\ge \alpha_1$, then the part of
 $P$ occupying the rows $>\ell(\beta)$ is a rectification of $T$. Then the internal insertion procedure is independent of a particular sequence of inner corners in $T$ chosen \cite{rssw}.  We may therefore have $P(T,{U})=P(T,{U}')$ with ${U}\neq {U}'$and as we have seen this happens  for example when  the shape $\beta$ of ${U}$ and   ${U}'$  is  a rectangle of height $\ge \ell(\mu)$ and width $\ge \mu_1$.
 We have the following characterization which only takes into account the shape of $U$ and is independent of $T$.
\begin{lem} \label{U} Let $T\in YT(\alpha/\mu)$ and  $U\in YT(\beta/\mu)$ where
 $\beta$ is of rectangle shape.  Then  $\cal{R}(U)\equiv$  $ {\ell(\beta)}^{\beta_1-\mu_{\ell(\beta)}}$ $ \cdots 2^{\beta_1-\mu_2}$ $1^{\beta_1-\mu_1},$ is a reverse Yamanouchi word of content $(\beta_1^{\ell(\beta)})-\mu=(\beta_1-\mu_1, \beta_1-\mu_2,\dots,\beta_1- \mu_{\ell(\beta)})$, and, therefore, $P(T,U)=P(T,U')$ is the skew tableau with inner shape the rectangle diagram $\emptyset_\beta$ where   below it is  the rectification of the subtableau of $T$ consisting of the first $\beta_1$ columns, and to the right of it is the subtableau of $T$ consisting of the last $max\{\alpha_1-\beta_1,0\}$ columns of $T$..
\end{lem}

 Knuth equivalence on $\cal{R}({U})$ and $\cal{R}({U}')$ provide only sufficient  conditions for the equality of $P(T,{U})=P(T,{U}')$.  In fact, we may have $P(T,{U})=P(T,{U}')$ and $\cal{R}(U)\not\equiv \cal{R}(U)'$. The example below illustrates this fact.
\begin{ex}\label{ex:sufficientcond}
\begin{enumerate}
\item  Companion  words of $U$ and $U'$:
\begin{align}
 U=\YT{0.16in}{}{
{,1,1},
{{2},2},
}\quad
std U=\YT{0.16in}{}{
{,1,2},
{{3},4},
}
\quad U'=\YT{0.16in}{}{
{,1,2},
{{1},2},
}\quad std U'=\YT{0.16in}{}{
{,2,4},
{{1},3},
}\nonumber
\end{align}
 \begin{eqnarray}\label{stand2}
\cal{R}(U)=2211=\cal{R}(std U),\;\cal{R}(U')=\cal{R}(std U')=1212.\nonumber
\end{eqnarray}

$(b)$ $\cal{ R}(U)$ and $\cal{ R}(U')$ are internal insertion order words of $T=\YT{0.16in}{}{
{,1,3},
{{2},3},
}$, and
\begin{align}
\phi_{221\bf 1}T
=\YT{0.16in}{}{
{,,},
{,,3},
{{1},3},
{2},
}=
\phi_{{\bf 1}212}T\nonumber
 \end{align}
 but $R(U)\not\equiv R(U')$.
 However, this property does not hold for every given $T$ having $\cal R(U)$ and $\cal R(U')$ as insertion words which shows that necessary and sufficient conditions for the equality $P(T,U)=P(T,U')$
 also depend on $T$. 

For  $T=\YT{0.16in}{}{
{,1,2,3},
{{3},3},
}$, one has
\begin{align}
\phi_{221\bf 1}T
=\YT{0.16in}{}{
{,,,3},
{,},
{{1},2},
{3,3},
}
\neq
\phi_{{\bf 1}212}T=\YT{0.16in}{}{
{,,,3},
{,,2},
{{1},3},
{3},
}\nonumber
\end{align}

 \item Companion words of $V$ and $V'$ below are internal insertion order words of

 $T'=\YT{0.16in}{}{
{,,1},
{,{2}},
{3},
}$,
 $
V=\YT{0.16in}{}{
{,,2},
{,{1}},
3}\equiv V'=\YT{0.16in}{}{
{,,2},
{,{3}},
1},\quad
 R(V)=312\not\equiv R(V')=213,$

$
\phi_{312}T'=\YT{0.16in}{}{
{,,},
{,,{1}},
{{}},
{2},
{3},}
\neq
\phi_{213}T'=\YT{0.16in}{}{
{,,},
{,},
{,1},
{2},
{3},}
$


\item
$
W_1=\YT{0.16in}{}{
{,1},
{{1},2},
},W_2=\YT{0.16in}{}{
{,1},
{{2},2},
},R(W_1)=212\equiv R(W_2)=221$,
  $
221\equiv 212\;\text{but}\;221{\bf{1}}\not\equiv{\bf 1} 212
$
\begin{align}
\phi_{R(W_2)}T=\phi_{221}T
&=\YT{0.16in}{}{
{,,3},
{,},
{{1},3},
{2},
}=
\phi_{212}T=\phi_{R(W_1)}T\nonumber
 \end{align}

 but as we have seen  above
$\phi_{221\bf 1}T=
\phi_{{\bf 1}212}T.
 $
On the other hand, if $T=\YT{0.16in}{}{
{,1,3},
{2,3},
}$,
$U=\YT{0.16in}{}{
{,1},
{2},
}
$ and $V=\YT{0.16in}{}{
{,3},
{1},
}
$,
$
\phi_2\phi_1 T=\YT{0.16in}{}{
{,,3},
{,3},
{1},
{2},
}
\neq\phi_1\phi_2 T=\YT{0.16in}{}{
{,,3},
{,1},
{2,3},
}
.$
However with
$U'=\YT{0.16in}{}{
{,2},
{1},
{3},
}
$ and $V'=\YT{0.16in}{}{
{,3},
{1},
{2},}
$, $R(U')=312\equiv R(V')=132$, one has
$
\phi_{312} T=
\phi_{132} T=\YT{0.16in}{}{
{,,3},
{,1},
{,3},
{2},}
.$

\end{enumerate}
\end{ex}

The following lemmata and proposition 
{ show that Sagan-Stanley internal row insertion operators satisfy Knuth relations.}
 Let $T$ be a SSYT   of shape $\lambda/\mu$ with $\ell(\lambda)=n$. When we write $\phi_u( T)$, for some word $u$,  it is assumed that there exists $U\in YT(\beta/\mu)$ such that $\cal{R}(U)=u$, that is, $u$ is an internal insertion order word for $T$.

Observe that if $m$ is the largest entry of $U\in YT(\beta/\mu)$ with $\gamma$ its content, and  $1\le d\le m$, we may decompose $ U= U_{|[d]}\cup U_{|[d+1,m]}$  and
$\std U
=(\std U)_{|[q]}\cup (\std U)_{|[q+1,|\gamma|]}$ with $q=\gamma_1+\cdots+\gamma_d$, and $\cal{R} (U)=\cal{R}( U_{|[d+1,m]}) \cal{R}( U_{|[d]})$. On the other hand, if for some $V\in YT(\gamma/\epsilon)$, $\cal{R}(\std V)=u_3u_2u_1=ijk$ and $1\le i\le k< j$ (respectively $\cal{R}(\std V)=u_3u_2u_1=kji$ and $1\le i<k\le j$) then $\gamma$ is not a row nor a column and it is easily checked that there exists $V'\in YT(\gamma/\epsilon)$ such that $\cal{R} (\std V')=jik$ (respectively $\cal{R} (\std V')=kij$). 
Without loss of generality, we may consider $i<k=i+1< j=i+2 \mbox{ or } i= k<j=i+1 $,
$V=\YT{0.16in}{}{
{,2},
{1},
{3},
}
 \mbox{ or } \tilde V=\YT{0.16in}{}{
{1,3},
{2},
} $ and $V'=\YT{0.16in}{}{
{,3},
{1},
{2},
}$  or  $ \tilde V'=\YT{0.16in}{}{
{1,2},
{3},}
$, $\cal{R}(V)=jik\equiv \cal {R}(V')=ijk$, $\cal{R}(\tilde  V)=iji\equiv \cal{R}(\tilde V')=jii$; and $i<k=i+1<j= i+2 $, 
$V=\YT{0.16in}{}{
{,1},
{,3},
{2},
}$,
$V'=\YT{0.16in}{}{
{,2},
{,3},
{1},
}$,
$\cal{R}(V)=kji\equiv \cal{R}(V')=kij$.

Reciprocally if $\cal{R}(\std U)=wz$ then $ U= Z\cup  W$ with $\cal{R}(\std W)=w$ and $\cal{R}(\std Z)=z$.
 Let $\cal{R}(\std U)=u_2jik u_1$ with $1\le i\le k< j$ and $u_1$, $u_2$ words. Decompose $\std U=U_1\cup V\cup U_2$ such that $\cal {R}(U_1)=u_1$, $\cal{R}(V)=ijk$ and $\cal{R}(U_2)=u_2$. Then $U':=U_1\cup V'\cup U_2$ has $\cal{R}(U')=u_2jiku_1$.
Therefore if $u$ is an internal insertion order word of $T$ and $u'\equiv u$ then $u'$ is also an internal insertion order word of $T$. Lemmata and proposition below show that the plactic class of  an internal  insertion order word of $T$ gives rise to the same $P$-tableau of $T$ in the  Sagan-Stanley Theorem \ref{th:internal}.

\begin{lem}\label{lem:knuth-intern3} Let $T\in YT(\lambda/\mu)$ with $\ell(\lambda)=n$. Let $n-1\le \ell(\mu)\le n$. Then

\begin{align}\bar\phi_i\bar\phi_n\bar\phi_k T &=\bar\phi_n\phi_i\bar\phi_k T ,\;\; \mbox{ $1\le i\le k< n$}.\label{eq:knuthn}
\end{align}

In addition, if $\ell(\mu)= n$,

\begin{align}
\bar\phi_i\bar\phi_{n+1}\bar\phi_k T&=\bar\phi_{n+1}\bar\phi_i\bar\phi_k T, \;\;\mbox{ $1\le i\le k< n+1$.}\label{eq:knuthn+1}
\end{align}

\end{lem}
\begin{proof} Recall that we are assuming that $(k,\mu_k+1)$ is an inner corner of $T$.

We start with identity \eqref{eq:knuthn}.
Let $T=U\ast V$ be the factorization of $T$ across  the $k$th row, $1\le i\le k< n$. Let $\phi_k V=\alpha\ast V'$ where $|\alpha|=0,1$ such that $\alpha$ is the empty word if and only if $\mu_k=\lambda_k$ and in this case $V'$ is obtained from $V$ by adding an empty box at the end of row $k$, otherwise the call $(k,\mu_k+1)$ is vacated from the entry $\alpha$ which migrates to the cell $(k+1,1)$. Let $\phi_i V'=\beta \ast V''$. Since $i\le k$, by Lemma \ref{lem:1}, $(a)$, if $|\alpha|=0$ then $|\beta|=0$ and $\phi_i(V')=V''$, and if $|\alpha|=1$, then  either $|\beta|=0$  or  $\alpha\le \beta$. Then, using the factorization \eqref{fact}, one has to analyse  the  cases

$\bf{I}.$   $|\alpha|=0$.
One has $\bar\phi_k V=V'$ and $\bar\phi_i V'=V''$, and it follows
\begin{align*}
\bar\phi_i\bar\phi_n\phi_k(U\ast V)&=\bar\phi_i\bar\phi_n(U\ast V')=\bar\phi_n U \ast\bar \phi_i V'=\bar\phi_n U \ast  V''=\bar\phi_n\bar\phi_i(U\ast V' )
\\
&=\bar\phi_n\bar\phi_i\bar\phi_k(U\ast V).
                                   \end{align*}

  $\bf{II}.$  $|\alpha|=1$.

$(a)$ If $|\beta|=0$,
one has $\phi_k V=\alpha\ast V'$ and $\phi_i V'=V''$, and
\begin{align*}
\bar\phi_i\phi_n\bar\phi_k(U\ast V)&=\bar\phi_i(\bar\phi_n(U\bigcdot\alpha)\ast V')=\bar\phi_n ( U.\alpha) \ast  V''=\phi_n \bar\phi_i(U\bigcdot\alpha)\ast  V'
\\
&=\phi_n\phi_i\phi_k(U\ast V).
                                   \end{align*}

$(b)$ If $|\beta|=1$,
one has  $\bar\phi_k V=\alpha\ast V'$ and $\bar\phi_i V'=\beta\ast V''$, with $\alpha\le \beta$,  that is, $\phi_i \phi_k V=\alpha\beta\ast V''$, $\alpha\le \beta$, and
                                   \begin{align*}
                                   \bar\phi_i\bar\phi_n\phi_k(U\ast V)&=\bar\phi_i\bar\phi_n(U.\bar\phi_k V)= \bar\phi_i\bar\phi_n[(U.\alpha)\ast V']\\
                                   &=\bar\phi_n (U\bigcdot\alpha).\bar\phi_i V'=[\bar\phi_n( U\bigcdot\alpha)\bigcdot\beta]\ast V''\\
                                   &=[\bar\phi_n( U\bigcdot\alpha\beta)]\ast V'',\;\mbox{ $\alpha\le \beta$;   Lemma \ref{lem:1}, $(a)$, external bumping version,}\\
                                   &=\bar\phi_n( U\bigcdot\alpha\beta\ast V'')=\bar\phi_n( U\bigcdot\bar\phi_i \bar\phi_k V)=\bar\phi_n\bar\phi_i( U.\bar\phi_kV)=\bar\phi_n\bar\phi_i\bar\phi_k( U \ast V).
\end{align*}
It is similarly checked that \eqref{eq:knuthn+1} holds.
\end{proof}

\begin{lem} \label{lem:knuth-intern2} Let $T\in YT(\lambda/\mu)$ with $\ell(\lambda)=n$. Let $n-1\le \ell(\mu)\le n$. Then

\begin{align}\bar\phi_k\bar\phi_i\bar\phi_n T =\bar\phi_k\bar\phi_n\bar\phi_i  T,\;\mbox{ $1\le i< k\le n$}.\label{eq:reverseknuth}
\end{align}

In addition if $\ell(\mu)=n$,

\begin{align}\label{n+1-2}
\bar\phi_k\phi_i\bar\phi_{n+1}T=\bar\phi_k\bar\phi_{n+1}\bar\phi_i  T, \;\mbox{ $1\le i< k\le n+1$.}
\end{align}
 \end{lem}
 \begin{proof}
 We  first  consider \eqref{eq:reverseknuth} with  $1\le i<k<n$. Let $T=U\ast W\ast V$ be the factorization of $ T$ across  the $(k-1)$th and $(n-1)$th rows.  Since $\ell(\lambda)=n$ and  $n-1\le \ell(\mu)\le n$, let $U=\emptyset_{\mu_n}. u$ with  $u=u_1\dots u_r$, $r=\lambda_n-\mu_n$, be the $n$th row of $T$ with $|u|\ge 0$ and $\mu=(\mu_1\dots, \mu_n)$, $\mu_n\ge 0$, and $\lambda=(\lambda_1,\dots,\lambda_n)$, $\lambda_n>0$. Using the factorization \eqref{fact}, one has

\begin{align}
\bar\phi_k\bar\phi_i\bar\phi_n(U\ast W\ast V)&=\bar\phi_k\phi_i[(\bar\phi_nU)\ast W\ast V)\nonumber \\
&=\bar\phi_k[\bar\phi_n U\bigcdot\bar\phi_i(W\ast V)], \mbox{ $i<n$},\nonumber\\
&=\bar\phi_k[(\phi_n U\bigcdot\beta)\ast W'\ast V'],\;\;\phi_i(W\ast V)=\beta\ast W'\ast V',\;|\beta|\ge 0,\nonumber\\
&=(\phi_n U\bigcdot\beta)\bigcdot\phi_k W'\ast V'\label{ident1}
                                   \end{align}

                                   and
                                   \begin{align}
                                  \bar\phi_k\bar\phi_n\bar\phi_i(U\ast W\ast V) &=\bar\phi_k\phi_n [U\bigcdot\phi_i(W\ast V)]=\phi_k\phi_n( U\bigcdot\beta\ast W'\ast V')\nonumber\\
                                  &=\bar\phi_n( U\bigcdot\beta)\bigcdot\bar\phi_k W'\ast V'\label{ident2}.
                                 \end{align}
                                 We want to show that \eqref{ident1} and \eqref{ident2} are equal, that is, $\bar\phi_n( U.\beta)=\bar\phi_n U.\beta$.
We have two main cases, either the $\phi_i$-bumping route reaches the $n$th row or not.

\text{I.} The $\phi_i$-bumping  route does not reach the $n$-th row, that is,  $|\beta|=0$. Then $\bar\phi_n( U\bigcdot\beta)=\bar\phi_nU=\phi_n U.\beta$, and $\eqref{ident1}=\eqref{ident2}=\bar\phi_n U\bigcdot\bar\phi_k W'\ast V'$.

\text{II.} The $\phi_i$-bumping  route  reaches the $n$-th row. That is  $|\beta|=1$, and $\beta$ is the $\phi_i$-bumped  element from the $(n-1)$th row to the $n$th row. Since $i<k$, by Lemma \ref{lem:1}, $(b)$,
$\phi_k W'=@\ast W''$ with $\beta>@$. Then
\begin{align*}
\eqref{ident1}&=(\bar\phi_n\bigcdot\beta)\bigcdot\bar\phi_kW'\ast V'=(\bar\phi_nU\bigcdot \beta\bigcdot @)\ast W''\ast V',\\
\eqref{ident2}&=\bar\phi_n(U\bigcdot\beta)\bigcdot \bar\phi_k W'\ast V'=\bar\phi_n(U\bigcdot\beta)\bigcdot @\ast W''\ast V'.
\end{align*}
 Next we consider the following cases according to the length of of $u$ in the $n$th row $U$ of $T$.

$(a)$  $|u|=0$.
One has $U=\emptyset_{\mu_n}$,  $ \bar\phi_n (U\bigcdot\beta)=\bar\phi_n (\emptyset_{\mu_n}.\beta)=\beta\ast \emptyset_{\mu_n+1}$. 
Then, since $\beta>@$, it follows

\begin{eqnarray}
\bar\phi_n U\bigcdot \beta\bigcdot @&=&\emptyset_{\mu_n+1}\bigcdot\beta\bigcdot @=\beta\ast \emptyset_{\mu_n+1}@,\nonumber\\
\phi_n( U\bigcdot\beta)\bigcdot @&=&\phi_n( \emptyset_{\mu_n}\beta)\bigcdot @=\beta\ast \emptyset_{\mu_n+1}@,\;\text{and}\nonumber\\
\eqref{ident1}=\eqref{ident2}&=&\beta\ast \emptyset_{\mu_n+1}@\ast W''\ast V'.
\end{eqnarray}

$(b)$ $|u|=1$.
One has $U=\emptyset_{\mu_n}u_1$,  $\bar\phi_n U.\beta=u_1\ast \emptyset_{\mu_n+1}\beta$ and
$$\bar\phi_n( U\bigcdot \beta)=\bar\phi_n (\emptyset_{\mu_n}u_1\bigcdot \beta)=\begin{cases}
\bar\phi_n (\emptyset_{\mu_n}u_1\beta),\;u_1\le \beta,\\
\phi_n (\emptyset_{\mu_n}u_1\ast \beta),\;u_1>\beta.\\
\end{cases}
$$

Henceforth
\begin{eqnarray*}
\phi_n U.\beta.@&=&(u_1\ast \emptyset_{\mu_n+1}\beta).@=\begin{cases}
u_1\beta \ast\emptyset_{\mu_n+1} @,\;u_1\le \beta,\\
\begin{smallmatrix}\beta\\
u_1\end{smallmatrix} \ast \emptyset_{\mu_n+1}@,\;u_1>\beta,\\
\end{cases}\\
\bar\phi_n (U.\beta).@&=&\begin{cases}
\phi_n (\emptyset_{\mu_n}u_1\beta)\bigcdot @=u_1\ast \emptyset_{\mu_n+1}\beta\bigcdot @=u_1\beta\ast \emptyset_{\mu_n+1}@,\;u_1\le \beta,\nonumber\\
\phi_n (u_1\ast \emptyset_{\mu_n+1}\beta)@=\begin{smallmatrix}\beta\\
u_1
\end{smallmatrix}\ast\emptyset_{\mu_n+1} @,\;u_1>\beta.\nonumber\\
\end{cases}
\end{eqnarray*}
Thus $\eqref{ident1}=\eqref{ident2}$ are equal to $u_1\beta\ast\emptyset_{\mu_n+1} @\ast W''\ast V'$, if $u_1\le \beta$,  or $\begin{smallmatrix}\beta\\
u_1
\end{smallmatrix}\ast \emptyset_{\mu_n+1}@\ast W''\ast V'$, otherwise.

$(c)$ $|u|\ge 2$.
One has $U=\emptyset_{\mu_n}u_1\cdots u_r$, with $r\ge 2$, and $\phi_n U= u_1\ast \emptyset_{\mu_n+1}u_2\cdots u_r$.

Either $u_r\le\beta$ or $u_r>\beta$.

$(i)$ $u_r\le \beta$.

Since $\beta>@$, let $x:=\min\{z\in\{u_2,\dots,u_r,\beta\}:z>@\}$ and $(u_2\dots u_r \beta).@=:x\ast u'$. Note $x\ge u_2\ge u_1$. Henceforth
\begin{eqnarray*}\bar\phi_n(U.\beta).@&=&\phi_n(\emptyset_{\mu_n}u_1\cdots u_r\beta).@
=u_1\ast (\emptyset_{\mu_n+1}u_2\cdots u_r\beta).@
= u_1x\ast \emptyset_{\mu_n+1}u',
\end{eqnarray*} and
\begin{eqnarray*}\phi_n U\bigcdot \beta\bigcdot @&=&u_1\ast (\emptyset_{\mu_n+1}u_2\cdots u_r)\bigcdot \beta.@
=u_1\ast (\emptyset_{\mu_n+1}u_2\cdots u_r\beta).@
= u_1x\ast \emptyset_{\mu_n+1}u'.
\end{eqnarray*}
Therefore, $\eqref{ident1}=\eqref{ident2}=u_1x\ast \emptyset_{\mu_n+1}u'\ast W''\ast V'.$

$(ii)$ It remains to study when $u_r>\beta$.

$(ii.1)$ $u_1>\beta$.
One has $@<\beta<u_1\le u_2\le \cdots\le u_r$,
\begin{eqnarray*}\bar\phi_n (U.\beta)\bigcdot @&=&\bar\phi_n(u_1\ast \emptyset_{\mu_n}\beta u_2\cdots u_r).@
=\begin{smallmatrix}
\beta\\
u_1
\end{smallmatrix}\ast(\emptyset_{\mu_n+1}u_2\cdots u_r).@\\
\\
&=&\begin{smallmatrix}
\beta&u_2\\
u_1&\\
\end{smallmatrix}\ast \emptyset_{\mu_n+1}@u_3\cdots u_r,
\end{eqnarray*} and
\begin{eqnarray*}\bar\phi_n U\bigcdot \beta.@&=&(u_1\ast \emptyset_{\mu_n+1} u_2\cdots u_r)\bigcdot \beta\bigcdot @=u_1u_2\ast (\emptyset_{\mu_n+1}\beta u_3\cdots u_r).@\\
\\&=&
\begin{smallmatrix}
\beta&u_2\\
u_1&\\
\end{smallmatrix}\ast \emptyset_{\mu_n+1}@u_3\cdots u_r.
\end{eqnarray*}
Thus, $\eqref{ident1}=\eqref{ident2}=\begin{smallmatrix}
\beta&u_2\\
u_1&\\
\end{smallmatrix}\ast \emptyset_{\mu_n+1}@u_3\cdots u_r\ast W''\ast V',$  where  $u_1u_2\beta\equiv u_1\beta u_2\equiv \begin{smallmatrix}
\beta&u_2\\
u_1&\\
\end{smallmatrix},$ with $\beta<u_1\le u_2$.

$(ii.2)$ $u_r\ge\cdots\ge u_i>\beta\ge u_{i-1}\ge\cdots\ge u_1$, for some $i\in\{2,\dots,r\}$.
One has $$U.\beta=\emptyset_{\mu_n}  u_1\cdots u_{i-1}u_i\cdots u_r.\beta=u_i\ast \emptyset_{\mu_n}u_1\cdots u_{i-1}\beta u_{i+1}\cdots u_r.$$
Let $$x:=\min\{z\in\{u_2,\dots,u_{i-1},\beta\}:z>@\}$$ and $(u_2\cdots u_{i-1} \beta u_{i+1}\cdots u_r).@=:x\ast u'.$ Note that $u_r\ge \cdots \ge u_{i+1}\ge u_i>\beta$ $\ge x$ $\ge u_2$ $\ge u_1$. Henceforth
\begin{eqnarray*}
\bar\phi_n(U.\beta).@&=&\bar\phi_n(u_i\ast \emptyset_{\mu_n}u_1\cdots u_{i-1}\beta u_{i+1}\cdots u_r).@=\begin{smallmatrix}u_1\\
u_i\\
\end{smallmatrix}\ast(\emptyset_{\mu_n+1}u_2\cdots u_{i-1}\beta u_{i+1}\cdots u_r).@\\
&=&\begin{smallmatrix}u_1&x\\
u_i&\\
\end{smallmatrix}\ast \emptyset_{\mu_n+1}u',
\end{eqnarray*} and
\begin{eqnarray*}
\bar\phi_n U.\beta.@&=&(u_1\ast \emptyset_{\mu_n+1}u_2\cdots  u_r).\beta.@=u_1u_i\ast(\emptyset_{\mu_n+1}u_2\cdots u_{i-1}\beta u_{i+1}\cdots u_r).@\\
&=&u_1u_i.x\ast \emptyset_{\mu_n+1}u',\;u_i>x\\
&=&\begin{smallmatrix}u_1&x\\
u_i&\\
\end{smallmatrix}\ast \emptyset_{\mu_n+1}u'.
\end{eqnarray*}

Again, $\eqref{ident1}=\eqref{ident2}= \begin{smallmatrix}u_1&x\\
u_i&\\
\end{smallmatrix}\ast \emptyset_{\mu_n+1}u'\ast W''\ast V'$ and
$u_iu_1x\equiv u_1u_ix\equiv\begin{smallmatrix}u_1&x\\
u_i&\\
\end{smallmatrix}$, with $u_1\le x<u_i$.

\medskip
We now study for $1\le i<k=n$, which is similarly checked.
\begin{eqnarray}\label{k=n1}
\bar\phi_n\phi_i\phi_n(U\ast V)&=&\phi_n(\bar\phi_n U.\phi_iV), \;\bar\phi_iV=\beta\ast V',\; \ell(\beta)=0,1\nonumber\\
&=&\phi_n(\phi_n U\bigcdot \beta)\ast V'),
\end{eqnarray}
\begin{eqnarray}\label{k=n2}
\bar\phi_n\phi_n\phi_i(U\ast V)&=&\bar\phi_n\phi_n (U\ast\phi_iV)=\phi_n \phi_n[( U\bigcdot \beta)\ast V'].
\end{eqnarray}

I. $\ell(\beta)=0$.
\begin{eqnarray*}\eqref{k=n1}=\phi_n(\phi_n U)\ast V'
=\phi_n^2 U\ast V'=\eqref{k=n2}.
\end{eqnarray*}

II. $\ell(\beta)=1$.

$(a)$ $\ell(u)=0$.
One has $U=\emptyset_{\mu_n}\Rightarrow \phi_n U=\emptyset_{\mu_n+1}\Rightarrow  \phi_n (U.\beta)=\beta\ast \emptyset_{\mu_n+1}$. Therefore
\begin{eqnarray*}\eqref{k=n1}=\bar\phi_n(\bar\phi_n U\bigcdot \beta)\ast V'=\bar\phi_n(\beta\ast \emptyset_{\mu_n+1})\ast V'=\beta\ast \emptyset_{\mu_n+2}\ast V'.
\end{eqnarray*} and
\begin{eqnarray*}\eqref{k=n2}=\bar\phi_n\phi_n (U.\beta)\ast V'=\bar\phi_n^2( y\beta)\ast V'=\bar\phi_n( \beta\ast \emptyset_{\mu_n+1})\ast V'=\beta\ast \emptyset_{\mu_n+2}\ast V'.
\end{eqnarray*}
$(b)$ $\ell(u)=1$.

One has $U=\emptyset_{\mu_n} u_1$, $\phi_n U=u_1\ast \emptyset_{\mu_n+1}$ and $\phi_n(U.\beta)=\phi_n(\emptyset_{\mu_n} u_1.\beta)$. Then
\begin{eqnarray*}\eqref{k=n1}=\phi_n\phi_n (U.\beta)\ast V'&=&\phi_n^2( \emptyset_{\mu_n}\beta)\ast V'\\
&=&\phi_n( \beta\ast \emptyset_{\mu_n+1})\ast V'=\begin{cases}
u_1\beta\ast \emptyset_{\mu_n+2}\ast V',\; u_1\le \beta,\\
\begin{smallmatrix}
\beta\\
u_1
\end{smallmatrix}\ast V',\;u_1>\beta,
\end{cases}
\end{eqnarray*}
\begin{eqnarray*}\eqref{k=n2}=\bar\phi_n\bar\phi_n (U\bigcdot \beta)\ast V'&=&\bar\phi_n^2(yu_1\bigcdot \beta)\ast V'
=\begin{cases}
\bar\phi_n^2(\emptyset_{\mu_n}u_1\beta)\ast  V',\; u_1\le \beta,\\
\bar\phi_n^2(u_1\ast y\beta)\ast  V'\;u_1>\beta,
\end{cases}\\
&=&
\begin{cases}
\bar\phi_n(u_1\ast \emptyset_{\mu_n+1}\beta)\ast  V',\; u_1\le \beta,\\
\bar\phi_n(\begin{smallmatrix}
\beta\\
u_1
\end{smallmatrix}u_1\ast \emptyset_{\mu_n+1})\ast  V',\;u_1>\beta,
\end{cases}\\
&=&
\begin{cases}
u_1\beta\ast \emptyset_{\mu_n+2}\ast  V'',\; u_1\le \beta,\\
\begin{smallmatrix}
\beta\\
u_1
\end{smallmatrix}\ast \emptyset_{\mu_n+2}\ast  V',\;u_1>\beta.
\end{cases}
\end{eqnarray*}

$(c)$ $\ell(u)\ge 2$.

Let $U=yu_1\cdots u_r$, $r\ge 2$.

$(i)$ $u_r\le \beta$
\begin{eqnarray*}
\bar\phi_n(\phi_n U\bigcdot \beta)=\bar\phi_n(u_1\ast \emptyset_{\mu_n+1}u_2\cdots u_r\beta)=u_1u_2\ast \emptyset_{\mu_n+1}u_3\cdots u_r\beta.
\end{eqnarray*}
\begin{eqnarray*}
\bar\phi_n\bar\phi_n( U\bigcdot \beta)=\bar\phi_n(u_1\ast \emptyset_{\mu_n+1}u_2\cdots u_r\beta)=u_1u_2\ast \emptyset_{\mu_n+1}u_3\cdots u_r\beta.
\end{eqnarray*}

$(ii)$ $u_1> \beta$
\begin{eqnarray*}
\phi_n(\phi_n U.\beta)&=&\phi_n(u_1\ast \emptyset_{\mu_n+1}u_2\cdots u_r.\beta)
=\phi_n(u_1u_2\ast \emptyset_{\mu_n+1}\beta u_3\cdots u_r)\\
&=&\begin{smallmatrix}
\beta&u_2 \\
u_1&
\end{smallmatrix}\ast \emptyset_{\mu_n+2}u_3\cdots u_r\ast V'.
\end{eqnarray*}
\begin{eqnarray*}
\bar\phi_n \phi_n (U\bigcdot \beta)\ast V'&=&\phi_n\phi_n(u_1\ast y\beta u_2\cdots u_r)\ast V'=
\phi_n(\begin{smallmatrix}
\beta \\
u_1
\end{smallmatrix}\ast \emptyset_{\mu_n+1}u_2 u_3\cdots u_r)\ast V'\\
&=&\begin{smallmatrix}
\beta&u_2 \\
u_1&
\end{smallmatrix}\ast \emptyset_{\mu_n+2}u_3\cdots u_r\ast V'.
\end{eqnarray*}
$(iii)$ $u_r\ge\cdots\ge u_i>\beta\ge u_{i-1}\ge\cdots\ge u_1$, for some $i\in\{2,\dots,r\}$.

One has $U.\beta=\emptyset_{\mu_n}  u_1\cdots u_{i-1}u_i\cdots u_r\bigcdot \beta=u_i\ast yu_1\cdots u_{i-1}\beta u_{i+1}\cdots u_r$.
\begin{align*}
\eqref{k=n1}&= \bar\phi_n[(u_1\ast yu_2\cdots u_r)\bigcdot\beta]\ast V'=
\bar\phi_n(u_1u_i\ast \emptyset_{\mu_n+1} u_2\cdots u_{i-1}\beta u_{i+1}\cdots u_r)\ast V'\\
&=\begin{smallmatrix}
u_1&u_2 \\
u_i&
\end{smallmatrix}\ast \emptyset_{\mu_n+2}u_3\cdots u_{i-1}\beta u_{i+1}\cdots u_r\ast V'
\end{align*}
\begin{align*}
\eqref{k=n2}&= \bar\phi_n \bar\phi_n (U.\beta)\ast V'=\bar\phi_n^2[(u_i\ast yu_1\cdots u_{i-1}\beta u_{i+1}\cdots u_r)]\ast V'\\
&=
\bar\phi_n(\begin{smallmatrix}
u_1 \\
u_i&
\end{smallmatrix}\ast \emptyset_{\mu_n+1} u_2\cdots u_{i-1}\beta u_{i+1}\cdots u_r)\ast V'\\
&=\begin{smallmatrix}
u_1&u_2 \\
u_i&
\end{smallmatrix}\ast \emptyset_{\mu_n+2}u_3\cdots u_{i-1}\beta u_{i+1}\cdots u_r\ast V'.
\end{align*}
It is similarly checked that \eqref{n+1-2} holds.
\end{proof}

\begin{obs} \label{obsF} Let $F$ be the tableau restricted to the   rows strictly below the $n$th row of $\bar\phi_k\bar\phi_i\bar\phi_n T=\bar\phi_k\bar\phi_n\bar\phi_i T$,for $1\le i< k\le  n$. If $w$ and $w'$ are the words consisting  of the elements of $T$ successively bumped out from the $n$-th row under the action of  $\bar\phi_k\bar\phi_i\bar\phi_n$ and  $\bar\phi_k\bar\phi_n\bar\phi_i$, for $1\le i< k\le  n$, on $ T$ respectively, then $w\equiv w'\equiv F$. This easily follows from the fact that $F$ is the external insertion of the elements of $T$ successively bumped out of the $n$-th row. This also applies to the action of $\bar\phi_i\bar\phi_n\bar\phi_k$ and $ \bar\phi_n\bar\phi_i\bar\phi_k$, for $1\le i\le k< n$, on $ T$.
\end{obs}

\begin{obs} \label{obsbar} Thanks to \eqref{def:operator} and \eqref{def:operatorbar}, lemmata \ref{lem:knuth-intern3} and \ref{lem:knuth-intern2} are generalized to $Y_\mu\cup T$.
\end{obs}


\begin{ex}\label{ex:knuth} Illustration of Lemma \ref{lem:knuth-intern3}, 
and of previous observations with $n=3$,
$$\begin{array}{ccccccccc}
\bar\phi_3 \bar\phi_1\bar\phi_2 T&=&\bar\phi_3 \bar\phi_1\bar\phi_2\YT{0.16in}{}{
{{\color{red}1},{\color{red}1},{\color{red}1},{\color{red}1},1,1},
{{\color{red}2},{\color{red}2},{1},{2},2},
{{\color{red}3},{2},3,3},
}&=&\bar\phi_3 \bar\phi_1\YT{0.16in}{}{
{{\color{red}1},{\color{red}1},{\color{red}1},{\color{red}1},1,1},
{{\color{red}2},{\color{red}2},{\color{red}2},{2},2},
{{\color{red}3},{1},3,3},
{2},
}\\
\\
&=&\bar\phi_3 \YT{0.16in}{}{
{{\color{red}1},{\color{red}1},{\color{red}1},{\color{red}1},{\color{red}1},1},
{{\color{red}2},{\color{red}2},{\color{red}2},{1},2},
{{\color{red}3},{1},2,3},
{2,3},
}
&=&\YT{0.16in}{}{
{{\color{red}1},{\color{red}1},{\color{red}1},{\color{red}1},{\color{red}1},1},
{{\color{red}2},{\color{red}2},{\color{red}2},{1},2},
{{\color{red}3},{\color{red}3},2,3},
{1,3},
{2},
}& w=231\equiv \YT{0.16in}{}{
{1,3},
{2},
}
\end{array}
$$
$$
\begin{array}{cccccccccccc}
\bar\phi_1 \bar\phi_3\bar\phi_2 T=\bar\phi_1 \bar\phi_3\bar\phi_2\YT{0.16in}{}{
{{\color{red}1},{\color{red}1},{\color{red}1},{\color{red}1},1,1},
{{\color{red}2},{\color{red}2},{1},{2},2},
{{\color{red}3},{2},3,3},
}=
\bar\phi_1 \bar\phi_3\YT{0.16in}{}{
{{\color{red}1},{\color{red}1},{\color{red}1},{\color{red}1},1,1},
{{\color{red}2},{\color{red}2},{\color{red}2},{2},2},
{{\color{red}3},{1},3,3},
{2},
}
=\bar\phi_1 \YT{0.16in}{}{
{{\color{red}1},{\color{red}1},{\color{red}1},{\color{red}1},1,1},
{{\color{red}2},{\color{red}2},{\color{red}2},{2},2},
{{\color{red}3},{\color{red}3},3,3},
{1},
{2},
}\quad
\\
=\YT{0.16in}{}{
{{\color{red}1},{\color{red}1},{\color{red}1},{\color{red}1},{\color{red}1},1},
{{\color{red}2},{\color{red}2},{\color{red}2},{1},2},
{{\color{red}3},{\color{red}3},2,3},
{1,3},
{2},
}
=\bar\phi_3 \bar\phi_1\bar\phi_2 T,\quad
w=231\equiv w'=213\equiv\YT{0.16in}{}{
{1,3},
{2},
}.
\end{array}$$
\end{ex}
 \begin{prop}\label{propp:knuth} ( Knuth relations of internal row insertion operators.) Let $Y\cup T$ be a  tableau pair  with $Y=Y_\mu$, $T$ a SSYT  of shape $\lambda/\mu$ and $\ell (\lambda)=n$. Suppose that $kij$ with $1\le i< k\le j\le n+1$, or $ijk$ with $1\le i\le k< j\le n+1$, are  internal insertion  order words of $T$. Then,  it holds
 \begin{eqnarray}
 \bar\phi_k\bar\phi_i\bar\phi_j(Y\cup T)&=&\bar\phi_k\bar\phi_j\bar\phi_i(Y\cup T), \;1\le i< k\le j\le n+1,\label{gen1}\\
  \bar\phi_i\bar\phi_j\bar\phi_k(Y\cup T)&=&\bar\phi_j\bar\phi_i\bar\phi_k(Y\cup T),\; 1\le i\le k< j\le n+1.\label{gen2}
  \end{eqnarray}
   More generally, if $u\equiv u'$ are internal insertion  order words of $T$, then
\begin{eqnarray}
\bar\phi_{u}(Y\cup T)=\bar\phi_{u'}(Y\cup T).\label{gen3}
\end{eqnarray}
 \end{prop}
 \begin{proof} Let $1<j\le n$ and $Y_\mu\cup T=(Y'\cup T')\ast(Y_{(\mu_1,\dots,\mu_{j})}\cup T^{[j]})$ with $Y'\cup T'$ the restriction of $Y\cup T$ to the last $n-j$ rows of $Y\cup T$. Lemmas  \ref{lem:knuth-intern3} and \ref{lem:knuth-intern2}  guarantee that
 $$\bar\phi_k\bar\phi_i\bar\phi_j(Y_{(\mu_1,\dots,\mu_{j})}\cup T^{[j]})=\bar\phi_k\bar\phi_j\bar\phi_i(Y_{(\mu_1,\dots,\mu_{j})}\cup T^{[j]}).$$
 Let $F$ be the restriction of $\bar\phi_k\bar\phi_i\bar\phi_j(Y_{(\mu_1,\dots,\mu_{j})}\cup T^{[j]})$ $=\bar\phi_k\bar\phi_j\bar\phi_i(Y_{(\mu_1,\dots,\mu_{j})}\cup T^{[j]})$ to the rows strictly below row $j$.
Recalling Remark \ref{obsF}, the action of   $\bar\phi_k\bar\phi_i\bar\phi_j$ and $\bar\phi_k\bar\phi_j\bar\phi_i$ on $Y\cup T$ inserts    words $w$ and $w'$ respectively, which are Knuth equivalent to $ F$, into  $Y'\cup T'$,
$$\bar\phi_k\bar\phi_i\bar\phi_j(Y\cup T)=[(Y'\cup T').w]\ast\bar\phi_k\bar\phi_i\bar\phi_j(Y_{(\mu_1,\dots,\mu_{j})}\cup T^{[j]}),$$
and
$$\bar\phi_i\bar\phi_j\bar\phi_k(Y\cup T)=[(Y'\cup T').w']\ast\bar\phi_k\bar\phi_j\bar\phi_i(Y_{(\mu_1,\dots,\mu_{j})}\cup T^{[j]}).$$
It follows from Lemma \ref{skewknuth}, $(b)$, that $ (Y'\cup T').w= (Y'\cup T'). w'=(Y'\cup T').F$, and thus \eqref{gen1} holds,
 \begin{eqnarray}\bar\phi_k\bar\phi_i\bar\phi_j(Y\cup T)&=&[(Y'\cup T').F]\ast\bar\phi_k\bar\phi_i\bar\phi_j(Y_{(\mu_1,\dots,\mu_{j})}\cup T^{[j]})\nonumber\\
&=&[(Y'\cup T').F]\ast\bar\phi_k\bar\phi_i\bar\phi_j(Y_{(\mu_1,\dots,\mu_{j})}\cup T^{[j]})=\bar\phi_i\bar\phi_j\bar\phi_k(Y\cup T).\end{eqnarray}
 Equality \eqref{gen2} follows similarly. Equality \eqref{gen3} follows from \eqref{gen1} and \eqref{gen2} and the definition of Knuth equivalent words.
 \end{proof}
\begin{ex} Illustration of this proposition using \eqref{route23} in Example \ref{ex:route}

$(a)$ Consider $Y\cup T$  as in
\eqref{route23} and $213\equiv 231$,
 $$\begin{array}{ccccccccc}
\bar\phi_2\bar\phi_1\bar\phi_3( \eqref{route23})&=&
\bar\phi_2\bar\phi_3\bar\phi_1( \eqref{route23})
&=&\YT{0.17in}{}{
{{\color{red}1},{\color{red}1},{\color{red}1},\blacksquare,1},
{{\color{red}2},{\color{red}2},{\color{red}2},\blacksquare},
{{\color{red}3},{\color{red}3},\blacksquare,$\textcircled 1$},
{{\color{red}4},1,$\textcircled{ \color{blue}2}$,$\textcircled 2$},
{1,$\textcircled{\color{blue}3}$,$\textcircled 3$},
{2,$\textcircled{ \color{blue} 4}$,$\textcircled 4$},
{$\textcircled 5$},
}
\end{array}$$
$(b)$  $312\equiv 132$,
\begin{eqnarray*}
\bar\phi_3 \bar\phi_1\bar\phi_2
\YT{0.16in}{}{
{{\color{red}1},{\color{red}1},{\color{red}1},{\color{red}1},1,1},
{{\color{red}2},{\color{red}2},{1},{2},2},
{{\color{red}3},{2},3,3},
{{1},3,4,4},
{5,5},
}
=
\YT{0.16in}{}{
{{1},3,4,4},
{5,5},
}
\,\bigcdot\,
\bar\phi_3 \bar\phi_1\bar\phi_2
\YT{0.16in}{}{
{{\color{red}1},{\color{red}1},{\color{red}1},{\color{red}1},1,1},
{{\color{red}2},{\color{red}2},{1},{2},2},
{{\color{red}3},{2},3,3},
}
\end{eqnarray*}
\begin{eqnarray*}
=
\YT{0.16in}{}{
{{1},3,4,4},
{5,5},
}
\,\bigcdot\,
\bar\phi_3 \bar\phi_1\YT{0.16in}{}{
{{\color{red}1},{\color{red}1},{\color{red}1},{\color{red}1},1,1},
{{\color{red}2},{\color{red}2},{\color{red}2},{2},2},
{{\color{red}3},{1},3,3},
{2},
}
=
\YT{0.16in}{}{
{{1},3,4,4},
{5,5},
}
\,\bigcdot\,
\bar\phi_3 \YT{0.16in}{}{
{{\color{red}1},{\color{red}1},{\color{red}1},{\color{red}1},{\color{red}1},1},
{{\color{red}2},{\color{red}2},{\color{red}2},{1},2},
{{\color{red}3},{1},2,3},
{2,3},
}
\end{eqnarray*}
\begin{eqnarray*}
=
\YT{0.16in}{}{
{{1},3,4,4},
{5,5},
}
\,\bigcdot\,
\YT{0.16in}{}{
{2,3},
}
\,\bigcdot\,
\bar\phi_3 \YT{0.16in}{}{
{{\color{red}1},{\color{red}1},{\color{red}1},{\color{red}1},{\color{red}1},1},
{{\color{red}2},{\color{red}2},{\color{red}2},{1},2},
{{\color{red}3},{1},2,3},
}
=
\YT{0.16in}{}{
{{1},3,4,4},
{5,5},
}
\,\bigcdot\,
\bar\phi_1 \bar\phi_3\bar\phi_2
\YT{0.16in}{}{
{{\color{red}1},{\color{red}1},{\color{red}1},{\color{red}1},1,1},
{{\color{red}2},{\color{red}2},{1},{2},2},
{{\color{red}3},{2},3,3},
}
\end{eqnarray*}
\begin{eqnarray*}
=\YT{0.16in}{}{
{{1},3,4,4},
{5,5},
}
\,\bigcdot\,
\YT{0.16in}{}{
{{\color{red}1},{\color{red}1},{\color{red}1},{\color{red}1},{\color{red}1},1},
{{\color{red}2},{\color{red}2},{\color{red}2},{1},2},
{{\color{red}3},{\color{red}3},2,3},
{1,3},
{2},
}
=
\YT{0.16in}{}{
{{1},3,4,4},
{5,5},
}
\,\bigcdot\,
\YT{0.16in}{}{{1,3},
{2},
}
\,\bigcdot\,
\YT{0.16in}{}{
{{\color{red}1},{\color{red}1},{\color{red}1},{\color{red}1},{\color{red}1},1},
{{\color{red}2},{\color{red}2},{\color{red}2},{1},2},
{{\color{red}3},{\color{red}3},2,3},
}\\
=\YT{0.16in}{}{
{{1},1,3,4},
{2,4},
{3,5},
{5},
}
\,\bigcdot\,
\YT{0.16in}{}{
{{\color{red}1},{\color{red}1},{\color{red}1},{\color{red}1},{\color{red}1},1},
{{\color{red}2},{\color{red}2},{\color{red}2},{1},2},
{{\color{red}3},{\color{red}3},2,3},
}
=
\YT{0.16in}{}{
{{\color{red}1},{\color{red}1},{\color{red}1},{\color{red}1},{\color{red}1},1},
{{\color{red}2},{\color{red}2},{\color{red}2},{1},2},
{{\color{red}3},{\color{red}3},2,3},
{{1},1,3,4},
{2,4},
{3,5},
{5},
}.
\end{eqnarray*}
\end{ex}

\section{
Switching  on ballot tableau pairs as a recursive Sagan-Stanley internal row insertion
}\label{switchingllkt}

\subsection{Littewood-Richardson rules and commutation symmetry}
Let $n$ be a fixed positive integer and let $x=(x_1,x_2,\ldots,x_n)$ be a sequence of
indeterminates. Then, for each partition $\lambda$ of length $\leq n$,
there exists a Schur function $s_\lambda(x)$ which is a homogeneous symmetric polynomial in the $x_k$
of total degree $|\lambda|$. The product of two Schur functions is explicitly given by the Littlewood-Richardson rule which amounts to finding how many SSYT's satisfy certain conditions.
\begin{thm} \cite{LR,gpthomas,sch}  The Littlewood-Richardson (LR) rule.  The coefficients appearing in the
expansion  of a product of Schur polynomials $s_\mu$ and $s_\nu$
\begin{equation}s_\mu(x)\ s_\nu(x)  = \sum_\lambda\ c_{\mu\nu}^\lambda\ s_\lambda(x)\label{schur}\end{equation}
are given by
$
  c_{\mu\nu}^\lambda = \#\{\text{ballot SSYT of shape $\lambda/\mu$ and content $\nu$}\}.
$
The coefficients $c_{\mu\nu}^\lambda$ are known as Littlewood--Richardson (LR) coefficients, and the ballot SSYT's are also known as Littlewood-Richardson tableaux.
\end{thm}

 The Schubert structure coefficients of the product in $H^*(G(d,n))$, the cohomology of the Grassmannian $G(d,n)$, (as a $\mathbb{Z}$-module),   are also given by the LR rule.
 The connection with $H^*(G(d,n))$, the cohomology of the Grassmannian $G(d,n)$ is due to L. Lesieur, \cite{lesieur}.
\begin{equation}\sigma_\mu \sigma_\nu=\sum_{\lambda\subseteq d\times(n-d)}{c_{\mu\;\nu}^\lambda} \sigma_\lambda.\label{schubert}\end{equation}
The Schur structure coefficients \eqref{schur} are not only
Schubert structure coefficients  \eqref{schubert}. They are also multiplicities in tensor products of $GL_n(\mathbb{C})$-
representations and in induction products of $\mathfrak{S}_n$-representations.

 The \emph{rectification} of a SSYT $T$ is the unique SSYT of normal shape whose reading word is Knuth equivalent to that of $T$. Using the notion of rectification of a SSYT,  the LR rule may also be formulated in the language of M.-P. Sch\"utzenberger's {\it jeu de taquin} \cite{stanley}. The SSYT's  $U$ and $V$ are said to be {\em jeu de taquin} equivalent if one can be obtained from another by a sequence of {\em jeu de taquin} slides \cite{stanley,fulton}.  Recall that each stage of {\em jeu de taquin slide} converts the reading word of a semistandard Young tableau  into a Knuth equivalent one, and  {\em jeu de taquin} commutes with standardisation.
\begin{thm}   Littlewood-Richardson rule's {\em jeu de taquin} version (\cite{stanley}, Appendix 1.) Fix a  standard tableau $S$ of shape $\nu$. Then
\begin{eqnarray}
  c_{\mu\nu}^\lambda &=& \#\{\text{ SYT of shape $\lambda/\mu$  whose rectification is  $S$}\}.\nonumber
\end{eqnarray}
(
The special choice of $S=\std Y_\nu$ relates this version with the LR tableau version above.)
\end{thm}
The rectification of a SSYT $T$ does not depend on the order of {\em jeu de taquin slides}.
 We can then consider $Y_\mu$ to be the inner shape of a SSYT of shape $\lambda/\mu$ and content $\nu$   and then look at $Y_\mu$ as a set of instructions to  tell where  {\em jeu de taquin} contracting slides start to rectify $T$. (Standardise $Y_\mu$ and $T$. The   {\em jeu de taquin} slides start with the biggest entry of $\std Y_{\mu}$, seen as a hole. The "hole" will slide until it reaches the outer boundary \cite{stanley}.  Then proceed similarly with {\em jeu de taquin} slides into the remaining  entries of $\std Y_\mu$, in decreasing numerical order. When $\std\, T$ is rectified to some $S$ the elements of $\std Y_\mu$ are the entries of the skew shape $\lambda/\nu$ and  encode the order in which the boxes were vacated in the {\em jeu de taquin} sliding process.) Such tableau sliding process correspond to a particular presentation of the  tableau switching procedure \cite{bss96} on  tableau pairs of partition shape. (When $S=\std Y_\nu$ one has the tableau sliding presentation for ballot tableau pairs of normal shape.)

 \emph{Tableau switching} process is  outlined in the next subsection.
Let $\rho_1$ be the  involution map  that  the  tableau switching procedure calculates on tableau pairs. We call it the \emph{switching involution}. The tableau sliding presentation of $\rho_1$  is called {\em infusion} by Thomas and Yong in \cite{tyong1, infusion}, and it can be translated to the language of Fomin's {\em jeu de taquin growths} (see Fomin's Appendix 1 in \cite{stanley}) which shows that  {\em infusion} is an involution. This approach is realised by  Thomas and Yong in \cite{tyong2} to exhibit an involution for the commutation of LR coefficients.

\subsection{The tableau switching map}
Switching ~\cite{bss96} is  an operation that takes a pair of tableaux $U\cup V$ and   moves them through each other  giving another
such pair $V'\cup U'$ of the same shape, in a way
that preserves Knuth equivalence, $V\equiv V'$ and $U\equiv U'$, and the shape of their union.  Loosely speaking, the switching algorithm may be realised as a mixture of Sch\"utzenberger's {\em jeu de taquin} and its reverse process in the sense that it calculates an involutive  map on pairs of tableaux and if $U\cup V$ has normal shape then $V'$ is the rectification of $V$ and $U$ is the rectification of $U'$.
The \emph{switching map} on $U\cup V$  is processed through  local moves to interchanging two  vertically or horizontally adjacent letters
 $\bf u$ and $\bf v$ from $U$ and $V$ respectively. The intermediate objects have to be defined when the switching procedure moves  $U$ and $V$ through each other.
A {\em  perforated tableau} $U$ of shape $\lambda$ is a filling of
some of the boxes in $\lambda$ with integers satisfying some restrictions: whenever
$\bf x$ and $\bf x'$ are entries of $U$ where $\bf x'$ is to the north-west of $\bf x$, then
\begin{align}&\mbox{ if not in the same column, }
\begin{smallmatrix}\bf{x'}&\\
\\
 &\bf{x}\end{smallmatrix},\quad
 \begin{smallmatrix}\bf{x'}&&\bf x \end{smallmatrix}\Rightarrow\;    \mathbf{x }\ge \mathbf{x}', \mbox{ and }\\
 &\mbox{  if in the same column }
 \begin{smallmatrix}\bf{x'}\\
 \\
 \bf{x}\end{smallmatrix}\Rightarrow\;  \mathbf{x}> \mathbf{x}'. \label{switch}
 \end{align}

We may switch an integer with the neighbour empty box $\blacksquare$ to the south, east, north or west, in a perforated tableau, so that the result is still a perforated tableau:
 $$\text{contracts $U$}\quad\begin{smallmatrix}\blacksquare&\\
 \bf{x}\end{smallmatrix}  {\underset s\rightarrow} \begin{smallmatrix}\mathbf{x}&\\
 \blacksquare&\end{smallmatrix}\qquad \begin{smallmatrix}\blacksquare&\mathbf{x}\end{smallmatrix}  {\underset s\rightarrow} \begin{smallmatrix}\mathbf{x}&\blacksquare\\
\end{smallmatrix}$$
$$\text{expands $U$}\quad\begin{smallmatrix}\blacksquare&\\
 \bf{x}\end{smallmatrix}  {\underset s\rightarrow} \begin{smallmatrix}\mathbf{x}&\\
 \blacksquare&\end{smallmatrix}\qquad \begin{smallmatrix}\blacksquare&\mathbf{x}\end{smallmatrix}  {\underset s\rightarrow} \begin{smallmatrix}\mathbf{x}&\blacksquare\\
\end{smallmatrix}$$
A {\em perforated tableau pair} $U\cup V$ of shape $\lambda$ is the superimposing of two perforated tableaux $U$ and $V$ of shape $\lambda$,
so that together they
completely fill $\lambda$.
If $\bf u$ and $\bf v$ are vertically or horizontally adjacent letters from $U$ and $V$ respectively, then an interchanging  of $\bf u$ with $\bf v$ is a {\em switch}, written ${\bf u } \underset{s}\leftrightarrow{\bf v}$, provided  it produces  a new perforated tableau pair,
   $$ \begin{smallmatrix}\bf u&\\
 \bf{v}\end{smallmatrix}  {\underset s\leftrightarrow} \begin{smallmatrix}\bf v&\\
 \bf u&\end{smallmatrix}\qquad \begin{smallmatrix}\bf u&\bf v\end{smallmatrix}  {\underset s\leftrightarrow} \begin{smallmatrix}\bf v&\bf u\\
\end{smallmatrix}.$$

  We collect the following elementary perforated tableau pairs with the corresponding elementary moves. For the sake of clarity, the entries of $U$ are drawn in red:
$$\begin{array}{cccccccccccccccccccc}
{\color{red} \rm u}&\rm v \;\;{\underset{\bf s}\leftrightarrow}&\rm v&{\color{red} \rm u},&\quad\quad
&&&&&{\color{red} \rm u}&{\underset{\bf s}\leftrightarrow}&\rm v\\
&&&&&\quad \quad &&&&\rm v&&{\color{red} \rm u}\\
\\
{\color{red} \rm u}&\rm v& {\underset{\bf s}\leftrightarrow} \;&\rm v&{\color{red}\rm u}&,\quad a> \rm v,&&&&{\color{red} \rm u}&a& {\underset{\bf s}\leftrightarrow } \;\;\; \rm v&{a}&,\quad a\ge \rm v.\\
a&&&{a}&&&&&&\rm v&&\;\;\;\;{\color{red} \rm u}
\end{array}.
$$
The above switches recover the {\em jeu de taquin } and {\em reverse jeu de taquin } switches when the ${\color{red} \rm u}$-entry  is seen as an empty entry.
An example of a sequence of switches
$$\YT{0.15in}{}{
 {{\color{red}1},{\color{red}1},{\color{red}1},{\color{red}1},{1},{1}},
 {{\color{red}2},{\color{red}2},{\color{red}2},{\color{red}2},2},
 {{\color{red}3},{1},{2},{3}},
}
{\underset{\bf s}\leftrightarrow}
\YT{0.15in}{}{
 {{\color{red}1},{1},{1},{\color{red}1},{\color{red}1},{\color{red}1}},
 {{\color{red}2},{\color{red}2},2,{\color{red}2},{\color{red}2}},
 {{1},{2},{3},{\color{red}3}},
}
{\underset{\bf s}\leftrightarrow}
\YT{0.15in}{}{
 {{\color{red}1},1,1,{\color{red}1},{\color{red}1},{\color{red}1}},
 {{1},{\color{red}2},2,{\color{red}2},{\color{red}2}},
 {{\color{red}2},{2},{3},{\color{red}3}},
}
{\underset{\bf s}\leftrightarrow}
\YT{0.15in}{}{
 {1,1,1,{\color{red}1},{\color{red}1},{\color{red}1}},
 {{\color{red}1},{\color{red}2},2,{\color{red}2},{\color{red}2}},
 {{\color{red}2},{2},{3},{\color{red}3}},
}.
$$

\begin{thm}{\em\cite[Theorem 2.3]{bss96}} The Switching Procedure  on  tableau pairs, calculates \emph{switching }  the unique map on tableau pairs with the following properties

$(1)$ Start with a  tableau pair $U\cup V$;

$(2)$ Switch integers from $U$ with integers from $V$ until it is no longer
possible to do so.
This produces a new  tableau pair $S\cup H$.
$$U\cup V\underset{\bf s}\longleftrightarrow S\cup H,$$
where $U\equiv H$ and $ V\equiv S$. (If $U\cup V$ has normal shape, $U$ is the rectification of $H$ and $S$  the rectification of $V$). If subtableaux decompose $U$, $V$ can switch with $U$ in stages. Similarly if $V$ decomposes, $U$ can switch with $V$ in stages. Switching  is an \emph{involution  on tableau pairs} denoted by  $\rho_1$.
\end{thm}
\begin{defi}\label{rho1}We  also  write $\rho_1^{(n)}$  when  the switching  $\rho_1$ is acting on  tableau pairs with shape length $\le n$.
\end{defi}

\begin{cor} The map $\rho_1^{(n)}$ is an involution in ${\cal{LR}}^{(n)}$, for all $n\ge 1$.
Moreover $c_{\mu,\nu}^\lambda=c_{\nu,\mu}^\lambda$.
 \end{cor}

For example, for $n=3$, $\mu=(4,4,1)$ and $\nu=(3,2,1)$,
\begin{align*}
Y_{(4,4,1)}\cup V&=\YT{0.14in}{}{
 {{\color{red}1},{\color{red}1},{\color{red}1},{\color{red}1},{1},{1}},
 {{\color{red}2},{\color{red}2},{\color{red}2},{\color{red}2},{2}},
 {{\color{red}3},{1},{2},{3}},
}
\underset{s}\leftrightarrow
\YT{0.14in}{}{
 {{\color{red}1},{\color{red}1},1,{\color{red}1},{\color{red}1},{1}},
 {{\color{red}2},{\color{red}2},{\color{red}2},2,{\color{red}2}},
 {{1},{2},{3},{\color{red}3}},
}
\underset{s}\leftrightarrow
\YT{0.14in}{}{
 {1,{\color{red}1},1,{\color{red}1},{\color{red}1},{1}},
 {{\color{red}1},{\color{red}2},{\color{red}2},2,{\color{red}2}},
 {{\color{red}2},{2},{3},{\color{red}3}},
}\\
&\underset{s}\leftrightarrow
Y_{(3,2,1)}\cup H=\YT{0.14in}{}{
 {{1},{1},{1},{\color{red}1},{\color{red}1},{\color{red}1}},
 {{2},{2},{\color{red}1},{\color{red}2},{\color{red}2}},
 {3,{\color{red}2}, {\color{red}2},{\color{red}3}},
}=\rho_1^{(3)}(Y_{(4,4,1)}\cup V)\in {\cal{LR}}^{(3)}.
\end{align*}

\subsection{
 Henriques-Kamnitzer $\mathfrak{gl}_n$-crystal commuter and  Sagan-Stanley internal row insertion}\label{subsec:skew}
For definitions, in this section, we refer the reader to \cite[Section 12]{akt16}, \cite[Section 3.4]{azkoma25} and \cite{sathishtorres}. { One way to conclude that a ballot tableau pair commuter $\rho$ coincides with the  Henriques-Kamnitzer $\mathfrak{gl}_n$--crystal commuter $Com_{HK}$ \cite{HenKam,HK2} is to show that $\rho(Y_\mu\cup T)=Y_\nu \cup H$ with  $H\equiv Y_\mu$ and $T\equiv Y_\nu$   satisfy the following: the corresponding left and right  Gelfand-Tsetlin (GT) pattern pairs $(G_\mu(T),G_\nu(T))$ respectively $(G_\nu(H),G_\mu(H))$ are related through  the Sch\"utzenberger involution $\xi$,  $G_\mu(H)=\xi(G_\mu(T))$  and $G_\nu(H)=\xi(G_\nu(T))$.  We shall show that our commuter ($\rho_3$ in \cite{pakvallejo}) based on the Sagan-Stanley internal insertion will produce such GT patterns in Theorem \ref{2} thanks to the coincidence of the commuter $\rho_1$ with $Com_{HK}$.
We will work on an  illustration as a motivation for Theorem \ref{2}.

Given $Y\cup T\in {\cal{LR}}^{(n)}$ with $Y=Y_\mu$  and $T$ a ballot tableau of shape $\lambda/\mu$ and weight $\nu$, we define for $i=1,\dots,n$, the partition $\nu^{(i)}$ to be the content of the ballot tableau $T^{[i]}$, and the partition $\widehat\nu=$ $(\widehat\nu_1,\dots,\widehat\nu_{n-1},\widehat\nu_n=\nu_n)\subseteq \nu$ where $\widehat\nu_i$ is the number of $i$'s in row $i$ of $T$, for  $i=1,\dots,n$.

The Sagan-Stanley internal row insertion correspondence \eqref{P} applied to the pair
$(\emptyset_{\widehat \nu},$ $ U)$  and $U\in YT(\nu/\hat\nu)$ with $\hat\nu\subseteq \nu$,
  produces $P$ equals to $\emptyset_\nu$ the  Young diagram of shape $\nu$.
Let  $G_\nu$  be the \emph{Gelfand-Tsetlin (GT) pattern} of type $\nu$ or the \emph{companion tableau} of $T$.
Note $G_\nu$ is defined by the nested sequence $ \nu^{(1)}\subseteq \nu^{(2)}\subseteq\cdots \subseteq\nu^{(n)}=\nu$. We use the skew GT pattern $G_{\nu/\widehat\nu}$, as internal insertion order tableau, and the parts of $\mu$ to be added properly, to construct an LR commuter as we explain next. { The skew GT pattern $G_{\nu/\widehat\nu}$ is the skew tableau of shape $\nu/\widehat\nu$ obtained by vacating the cells in the Young diagram of shape $\widehat\nu$ in the GT pattern $G_\nu$. See Example \ref{ex:skewgt}}

 The internal insertion  order word $\cal{R}(G_{\nu/\widehat\nu})=V_n\cdots V_3V_2V_1$ is decomposed into row words $V_i$, $i=1,\cdots, n$. The row word $V_i$ is  precisely the $i$th row word  of $T$ restricted to the alphabet $[i-1]$, and $|V_i|=\lambda_i-\mu_i-\hat\nu_i$, for $i=1,\dots,n$. Thus $V_1$ is the empty word,  $\phi_{V_1}=1$, $\phi_{\cal{R}(G_{\nu/\widehat\nu})}=\phi_{V_n}\cdots\phi_{V_3}\phi_{V_2}$ and $$\bar\phi_{\cal{R}(G_{\nu/\widehat\nu})}(\emptyset_{\hat\nu})=\bar\phi_{V_n}\cdots\bar\phi_{V_3}\bar\phi_{V_2}(\emptyset_{\hat\nu})=\emptyset_{\nu},$$ the Young diagram of shape $\nu$. For $i=1,\dots,n$, one needs the operator $\bar\chi_i$, to be  iterated $\mu_i$ times, written $\bar\chi_i^{\mu_i}=\underbrace{\bar\chi_i\circ\cdots\circ\bar\chi_i}_{\mu_i}$,
over $\bar\phi_{V_{i}}(\bar\chi_{i-1}^{\mu_{i-1}}\bar\phi_{V_{i-1}})\cdots(\bar\chi_{2}^{\mu_{2}}\bar\phi_{V_{2}})\bar\chi_1^{\mu_1}(\emptyset_{\hat\nu})$, to recursively adding  in each iteration one $i$  at the end of $i$th row. Given $Y_\mu\cup T\in {\cal{LR}}^{(n)}$ this procedure encodes $Y_\mu\cup T$ in the form $(\emptyset_{\widehat\nu}, G_{\nu/\widehat \nu},\mu)$ and gives the \emph{$\mu$-augmented internal insertion operator} an involution in ${\cal{LR}}^{(n)}$

\begin{align}\label{recursion0}\bar\phi^\mu_{\cal{R}( G_{\nu/\widehat \nu})}(\emptyset_{\widehat\nu})=(\bar\chi_n^{\mu_n}\bar\phi_{V_n})\cdots(\bar\chi_3^{\mu_3}\bar\phi_{V_3})(\bar\chi_2^{\mu_2}
\bar\phi_{V_{2}})
\bar\chi_1^{\mu_1}(\emptyset_{\widehat\nu})&=\emptyset_\nu\cup H\in {\cal{LR}}^{(n)},
\end{align}
where $H$ is a ballot tableau of shape $\lambda/\nu$ and content $\mu$.

\begin{ex}\label{ex:skewgt} Let $\lambda=(6,5,5,4,3)$, $\nu=(4,4,3,2,0) $, $\widehat\nu=(2,1,1,1,0)$ and $\mu=(4,3,2,1,0)$
\begin{eqnarray}\label{gtnu}Y_{\mu}\cup T&=&\YT{0.15in}{}{
 {{\color{red}1},{\color{red}1},{\color{red}1},{\color{red}1},{1},{1}},
 {{\color{red}2},{\color{red}2},{\color{red}2},{1},{2}},
 {{\color{red}3},{\color{red}3},{1},{2},{3}},
 {{\color{red}4},{2},{3},{4}},
 {{2},{3},{4}},
}\quad  G_{\nu}(T)=
\YT{0.15in}{}{
 {1,1,{2},{3}},
 {2,{3},{4},{5}},
 {3,{4},{5}},
 {4,5},
}
\quad
 G_{\nu/\widehat \nu}(T)=
\YT{0.15in}{}{
 {,,{2},{3}},
 {,{3},{4},{5}},
 {,{4},{5}},
 {,5},
}
\quad \std G_{\nu/\widehat \nu}=
\YT{0.15in}{}{
 {,,{1},{3}},
 {,{2},{5},{8}},
 {,{4},{7}},
 {,6},
}\nonumber
\\
\\
Y_{\mu}&=&\YT{0.15in}{}{
 {{\color{red}1},{\color{red}1},{\color{red}1},{\color{red}1}},
 {{\color{red}2},{\color{red}2},{\color{red}2}},
 {{\color{red}3},{\color{red}3}},
 {{\color{red}4}},
}\quad 
\cal{R}( G_{\nu/\widehat \nu})={\cal{R}}_5{\cal{R}}_4{\cal{R}}_3\,{\cal{R}}_2\,\emptyset={234}\,{23}\,{12}\,{1}=V_5V_4V_3V_2.\nonumber
\end{eqnarray}
Then for $n=5$,

 $$\bar\phi_{\cal{R}( G_{\nu/\widehat \nu})}(\emptyset_{\widehat\nu})=\bar\phi_{{\cal{R}}_5{\cal{R}}_4{\cal{R}}_3\,{\cal{R}}_2\,\emptyset}(\emptyset_{\widehat \nu})
=\bar\phi_{234}\bar\phi_{23}\bar\phi_{12}\bar\phi_{1}(\emptyset_{\widehat\nu})=\emptyset_\nu.$$

The $\mu$-augmentation of $\bar\phi_{\cal{R}( G_{\nu/\widehat \nu})}$:
\begin{align*}\bar\phi^\mu_{\cal{R}( G_{\nu/\widehat \nu})}(\emptyset_{\widehat \nu})&=\bar\chi_5^{\mu_5}\bar\phi_{{\cal{R}}_5}\bar\chi_4^{\mu_4}\bar\phi_{\cal{R}_4}\bar\chi_3^{\mu_3}\bar\phi_{\cal{R}_3}\bar\chi_2^{\mu_2}
\bar\phi_{\cal{R}_2}
\bar\chi_1^{\mu_1}\bar\phi_\emptyset(\emptyset_{\widehat \nu})\\
&=\bar\phi_{234}\bar\chi_4\bar\phi_{23}\bar\chi_3^{2}\bar\phi_{12}\bar\chi_2^{3}\bar\phi_1
\bar\chi_1^{4}(\emptyset_{\widehat\nu})\\
&=\emptyset_{\nu}\cup H
\end{align*}

gives
\begin{eqnarray}\emptyset_{\widehat\nu}=\YT{0.15in}{}{
 {,},
{{}},
{{}},
{{}},}
{\underset{\bar\chi_1^{4}}\rightarrow}\YT{0.15in}{}{
 {,,\color{red}{1},\color{red}{1},\color{red}1,\color{red}1},
{{}},
{{}},
{{}},
}
{\underset{\bar\chi_2^{3}\bar\phi_1}\rightarrow}
\YT{0.15in}{}{
 {,,,\color{red}1,\color{red}1,\color{red}1},
{{},\color{red}1,\color{red}2,\color{red}2,\color{red}2},
{{}},
{{}},
}
{\underset{\bar\chi_3^{2}\bar\phi_1\bar\phi_2}\rightarrow}
\YT{0.15in}{}{
 {,,,,\color{red}1,\color{red}1},
{{},,\color{red}1,\color{red}2,\color{red}2},
{{},\color{red}1,\color{red}2,\color{red}3,\color{red}3},
{{}},
}
\nonumber\\
\nonumber\\
{\underset{\bar\chi_4\bar\phi_2\bar\phi_3}\rightarrow }
\YT{0.15in}{}{
 {,,,,\color{red}1,\color{red}1},
{{},,,\color{red}2,\color{red}2},
{{},,\color{red}1,\color{red}3,\color{red}3},
{{},\color{red}1,\color{red}2,\color{red}4},
}
{\underset{\bar\phi_2\bar\phi_3\bar\phi_4}\rightarrow }
\YT{0.15in}{}{
 {,,,,\color{red}1,\color{red}1},
{{},,,,\color{red}2},
{{},,,\color{red}2,\color{red}3},
{{},,\color{red}1,\color{red}3},
{\color{red}1,\color{red}2,\color{red}4},
}=\emptyset_{\nu}\cup H.\label{operator}
\end{eqnarray}
where $H$ is a ballot tableau of content $\mu$ and shape $\lambda/\nu$. This augmented operator creates the skew GT pattern $G_{\nu/\widehat \nu}$ defined by the sequence of inner shapes read in \eqref{operator} onwards and displayed in \eqref{nuhat}.

The inverse of the $\mu$-augmented operator is the $\mu$-deletion operator
$$\Delta^\mu(\emptyset_\nu\cup H)=(\bar\chi_1^{\mu_1})^{-1}\Delta_2^{\lambda_2-\nu_2-\widehat\mu_2}(\bar\chi_2^{\mu_2})^{-1}\Delta_3^{\lambda_3-\nu_3-\widehat\mu_3}(\chi_3^{\mu_3})^{-1}\Delta_4^{\lambda_4-\nu_4-\widehat\mu_4}
(\bar\chi_4^{\mu_4})^{-1}\Delta_5^{\lambda_5-\nu_5-\widehat\mu_5}  (\emptyset_\nu\cup H)$$

\noindent which  reads \eqref{operator} backwards and thereby creates $G_{\nu/\widehat\nu}$ the  nested sequence of inner shapes
\begin{align}\label{nuhat}\nu=(4,4,3,2,0)=\widehat\nu^{(5)}\supseteq (4,3,2,1,0)=\widehat\nu^{(4)}\supseteq (4,2,1,1,0)=\widehat\nu^{(3)}\supseteq(3,1,1,1,0)=\widehat\nu^{(2)}\supseteq(2,1,1,1,0)=\widehat\nu
\end{align}

\end{ex}

\begin{ex}\label{ex:skewgt2} Let $\lambda=(6,5,5,4,3)$, $\nu=(4,4,3,2,0) $, $\widehat\mu=(2,1,1,0,0)$ and $\mu=(4,3,2,1)$
\begin{eqnarray}\label{gtmu}Y_{\nu}\cup H&=&\YT{0.15in}{}{
 {{\color{red}1},{\color{red}1},{\color{red}1},{\color{red}1},{1},{1}},
 {{\color{red}2},{\color{red}2},{\color{red}2},{\color{red}2},{2}},
 {{\color{red}3},{\color{red}3},{\color{red}3},{2},{3}},
 {{\color{red}4},{\color{red}4},{1},{3}},
 {{1},{2},{4}},
}\quad
G_{\mu}(H)=
\YT{0.15in}{}{
 {1,1,{4},{5}},
 {2,{3},{5}},
 {3,{4}},
 {{5}},
}
\quad
 G_{\mu/\widehat \mu}(H)=
\YT{0.15in}{}{
 {,,{4},{5}},
 {,{3},{5}},
 {,{4}},
 {{5}},
}
\\
Y_{\nu}&=&\YT{0.15in}{}{
 {{\color{red}1},{\color{red}1},{\color{red}1},{\color{red}1}},
 {{\color{red}2},{\color{red}2},{\color{red}2},{\color{red}2}},
 {{\color{red}3},{\color{red}3},{\color{red}3}},
 {{\color{red}4},{\color{red}4}},
}
\quad \cal{R}( G_{\mu/\widehat \mu})={\cal{R}}_5{\cal{R}}_4{\cal{R}}_3\emptyset\,\emptyset={124}\,{13}\,{2}\emptyset\emptyset.\nonumber
\end{eqnarray}
Then $\bar\phi_{\cal{R}( G_{\mu/\widehat \mu})}(\emptyset_{\widehat\mu})=\bar\phi_{124}\bar\phi_{13}\bar\phi_{2}(\emptyset_{\widehat\mu})=\emptyset_\mu,$ and
$\bar\phi^\nu_{\cal{R}( G_{\mu/\widehat \mu})}(\emptyset_{\widehat\mu})=\bar\phi_{124}\bar\chi_4^2\bar\phi_{13}\bar\chi_3^{3}\bar\phi_{2}\bar\chi_2^{4}
\bar\chi_1^{4}(\emptyset_{\widehat\mu})$ $=$ $\emptyset_{\mu}\cup T$ are displayed below
\begin{eqnarray}\emptyset_{\widehat\mu}=\YT{0.15in}{}{
 {,},
{{}},
{{}}
}
{\underset{\bar\chi_1^{4}}\rightarrow}\YT{0.15in}{}{
 {,,\color{red}{1},\color{red}{1},\color{red}1,\color{red}1},
{{}},
{{}},
}
{\underset{\bar\chi_2^{4}}\rightarrow}
\YT{0.15in}{}{
 {,,\color{red}1,\color{red}1,\color{red}1,\color{red}1},
{{},{\color{red}2},\color{red}2,\color{red}2,\color{red}2},
{{}},
}
{\underset{\bar\chi_3^{3}\bar\phi_2}\rightarrow}
\YT{0.15in}{}{
 {,,\color{red}1,\color{red}1,\color{red}1,\color{red}1},
{{},,\color{red}2,\color{red}2,\color{red}2},
{{},\color{red}2,\color{red}3,\color{red}3,\color{red}3},
}\nonumber
\\
\nonumber\\
{\underset{\bar\chi_4^2\bar\phi_1\bar\phi_3}\rightarrow }
\YT{0.15in}{}{
 {,,,\color{red}1,\color{red}1,\color{red}1},
{{},,\color{red}1,\color{red}2,\color{red}2},
{{},,\color{red}2,\color{red}3,\color{red}3},
{\color{red}2,\color{red}3,\color{red}4,\color{red}4},
}
{\underset{\bar\phi_1\bar\phi_2\bar\phi_4}\rightarrow }
\YT{0.15in}{}{
 {,,,,\color{red}1,\color{red}1},
{{},,,\color{red}1,\color{red}2},
{{},,\color{red}1,\color{red}2,\color{red}3},
{{},\color{red}2,\color{red}3,\color{red}4},
{\color{red}2,\color{red}3,\color{red}4},
}=\emptyset_{\mu}\cup T.\label{operator2}
\end{eqnarray}

The deletion operator is obtained from \eqref{operator2} backwards while recording the inner shapes to define $G_{\mu/\widehat\mu}$ the skew companion tableau of $H$
$$\mu=(4,3,2,1,0)\supseteq \widehat\mu^{(4)}=(3,2,2,0,0)\supseteq \widehat\mu^{(3)}=(2,2,1,0,0)\supseteq \widehat\mu^{(2)}=(2,1,1,0,0)\supseteq\widehat\mu=(2,1,1,0,0).$$


$$\Delta^\nu(\emptyset_\mu\cup T)=(\bar\chi_1^{\nu_1})^{-1}\Delta_2^{\lambda_2-\mu_2-\widehat\nu_2}(\bar\chi_2^{\nu_2})^{-1}\Delta_3^{\lambda_3-\mu_3-\widehat\nu_3}(\chi_3^{\nu_3})^{-1}
\Delta_4^{\lambda_4-\mu_4-\widehat\nu_4}
(\bar\chi_4^{\nu_4})^{-1}\Delta_5^{\lambda_5-\mu_5-\widehat\nu_5}  (\emptyset_\mu\cup T)=\emptyset_\nu\cup H$$
$$H\equiv Y_\mu$$
\end{ex}


\begin{obs}\cite{fultonbuch,pakvallejo,HenKam, akt16, tka18,azkoma25} Recall $\rm{evac} G_\mu(T)=G_\mu(H)$ where $\rm{evac}=\xi$ denotes the Sch\"utzenberger evacuation

$$\rm{evac} G_\mu(T)=\rm{evac}  \YT{0.15in}{}{
 {1,1,{1},{3}},
 {2,{2},{4}},
 {3,{5}},
 {{5}},
}=G_{\mu}(H)=
\YT{0.15in}{}{
 {1,1,{4},{5}},
 {2,{3},{5}},
 {3,{4}},
 {{5}},
}$$

\end{obs}

To avoid the  skew GT patterns $G_{\nu/\widehat\nu}$  and $G_{\mu/\widehat\mu}$ in the previous examples, and since $G_\nu=Y_{\widehat\nu}\cup G_{\nu/\widehat\nu}$ and similarly $G_\mu=Y_{\widehat\mu}\cup G_{\mu/\widehat\mu}$, we get supplied  with another adding operator $\bar\omega_i$ regarding the partitions $\widehat\nu=(\widehat\nu_1,\dots,\widehat\nu_n)$ or $\widehat\mu=(\widehat\mu_1,\dots,\widehat\mu_n)$  as follows.

Put $(Y\cup T)^{[0]}:=\emptyset$. For $i=0,1,\dots,n$, let $(Y_\mu\cup T)^{[i]}=Y_{(\mu_1,\dots,\mu_i)}\cup T^{[i]}\in {\cal LR}^{(i)}$  where $T^{[i]}$ has content $\nu^{(i)}$, and  row $i$ consists of  $V_i$,   the  row subword  restricted to the  entries in $[i-1]$, followed with   $\widehat\nu_i$ $i$'s.
For each $i=1,\dots,n$, we now consider the \emph{operator $\bar\omega_i$}, to be  iterated $\hat\nu_i$ times,
to  contributing,  in each iteration over
\begin{align*} Y_{\nu^{(i-1)}}\cup H^{(i-1)}=(\bar\chi_{i-1}^{\mu_{i-1}}\bar\phi_{V_{i-1}}\bar\omega_{i-1}^{\hat\nu_{i-1}})\cdots(\bar\chi_2^{\mu_2}
\bar\phi_{V_2}\omega_2^{\hat\nu_2}
)({\bar \chi_1^{\mu_1}\bar\omega_1^{\hat\nu_1}})(\emptyset),\\
 \mbox{ $H^{(i-1)}\equiv Y_{(\mu_1,\dots,\mu_{i-1})}$, skew shape $(\lambda_1,\dots,\lambda_{i-1})/\nu^{(i-1)}$ ,\; $H^{(0)}=\emptyset $,}
 \end{align*}
\noindent  with one $i$ to the $i$th row of the inner shape $Y_{\nu^{(i-1)}}$.
This allows to give the following  recursive presentation of the switching map $\rho_1^{(n)}$ on ballot tableau pairs in ${\cal LR}^{(n)}$, for $n\ge 1$,

\begin{align}\rho_1^{(n)}(Y_\mu\cup T)&=(\bar\chi_n^{\mu_n}\bar\phi_{V_n}\bar\omega_n^{\widehat\nu_n})
\cdots(\bar\chi_3^{\mu_3}\bar\phi_{V_3}\bar\omega_3^{\widehat\nu_3})(\bar\chi_2^{\mu_2}
\bar\phi_{V_2}\bar\omega_2^{\widehat\nu_2})
(\bar\chi_1^{\mu_1}\bar\omega_1^{\widehat\nu_1})(\emptyset)\label{lrnu1}\\
&=(\bar\chi_n^{\mu_n}\bar\phi_{V_n}\omega_n^{\widehat\nu_n})\rho^{(n-1)}[(Y_\mu\cup T)^{[n-1]}]\label{lrnu2}\\
&=Y_\nu\cup H, \mbox{ with $H\equiv Y_\mu$} .\label{lrnu3}
\end{align}
The Yamanouchi tableau $Y_{\hat\nu}$ is constructed recursively by means of $\hat\nu_i$ iterations of the operator $\omega_i$, for $i=1,\dots,n$.
In particular, if $V$ is the empty word and $\mu$ is the zero partition, then $\hat\nu=\nu$ and $Y_\nu=\bar\omega_n^{\nu_n}\cdots \bar\omega_2^{\nu_2}\bar\omega_1^{\nu_1}(\emptyset)$. Observe that $\bar\omega_i$ commutes with $\bar\phi_{V_i}$.



\begin{ex}\label{ex:main} Let us resume to  Example \ref{ex:skewgt} where  $\nu=(4,4,3,2,0) $, $\widehat\nu=(2,1,1,1,0)$ and $\mu=(4,3,2,1)$. We illustrate
\begin{align}
\rho_1^{(5)}(Y_\mu\cup T)=\bar\phi_2\bar\phi_3\bar\phi_4\bar\chi_4^{\mu_4}\bar\phi_2\bar\phi_3 \bar\omega_4\bar\chi_3^{\mu_3}\bar\phi_1\bar\phi_2\bar\omega_3\bar\chi_2^{\mu_2}\bar\phi_1\bar\omega_2
\bar\chi_1^{\mu_1}\bar\omega_1^2(\emptyset).\nonumber
\end{align}
 Recall $\bf{s}$ denotes the switching operator.
\begin{align*}
(Y\cup T)^{(0)}&=\emptyset\\(Y\cup T)^{[1]}&=
\YT{0.15in}{}{
 {{\color{red}1},{\color{red}1},{\color{red}1},{\color{red}1},{1},{1}},
 }
\\
(Y\cup T)^{[2]}&=\YT{0.15in}{}{
 {{\color{red}1},{\color{red}1},{\color{red}1},{\color{red}1},{1},{1}},
 {{\color{red}2},{\color{red}2},{\color{red}2},{1},{2}},
 }
 \\
(Y\cup T)^{[3]}&=\YT{0.15in}{}{
 {{\color{red}1},{\color{red}1},{\color{red}1},{\color{red}1},{1},{1}},
 {{\color{red}2},{\color{red}2},{\color{red}2},{1},{2}},
 {{\color{red}3},{\color{red}3},{1},{2},{3}},
 }\\
(Y\cup T)^{[4]}&= \YT{0.15in}{}{
 {{\color{red}1},{\color{red}1},{\color{red}1},{\color{red}1},{1},{1}},
 {{\color{red}2},{\color{red}2},{\color{red}2},{1},{2}},
 {{\color{red}3},{\color{red}3},{1},{2},{3}},
 {{\color{red}4},{2},{3},{4}},
 }\\
 (Y\cup T)^{[5]}&=Y\cup T.
\end{align*}
We check $\rho_1^{(i)}(Y_\mu\cup T)^{[i]}=(\bar\chi_i^{\mu_i}\bar\phi_{V_i}\bar\omega_i^{\widehat\nu_i})\rho^{(i-1)}[(Y_\mu\cup T)^{[i-1]}],$ for $i=1,\dots,5$.

 \begin{align*}
 (Y\cup T)^{[1]}&=
\YT{0.15in}{}{
 {{\color{red}1},{\color{red}1},{\color{red}1},{\color{red}1},{1},{1}},
 } {\underset s\rightarrow}\YT{0.15in}{}{
 {{1},{1},{\color{red}1},{\color{red}1},{\color{red}1},{\color{red}1}},
}\\
&=\rho_1^{(1)}[(Y\cup T)^{[1]}]\\
&=Y_{\widehat\nu_1}\cup H^{(1)},\quad H^{(1)}\equiv Y_{(\mu_1)}
 \end{align*}

 \begin{align*}
(Y\cup T)^{[0]}=Y_0\cup H^{(0)}&=\emptyset {\underset{\omega_1^{2}}\longrightarrow}\YT{0.15in}{}{
 {{1},{1}},
}{\underset{\chi_1^{\mu_1}}\longrightarrow}\YT{0.15in}{}{
 {{1},{1},{\color{red}1},{\color{red}1},{\color{red}1},{\color{red}1}},
}\\
&=\chi_1^{\mu_1}\omega_1^{\widehat\nu_1}(Y\cup T)^{[0]}\\
&=\rho_1^{(1)}[(Y\cup T)^{[1]}]
\\
&= Y_{\widehat\nu_1}\cup H^{(1)},\quad H^{(1)}\equiv Y_{(\mu_1)}
\end{align*}
\begin{align*}
(Y\cup T)^{[2]}&=\YT{0.15in}{}{
 {{\color{red}1},{\color{red}1},{\color{red}1},{\color{red}1},{1},{1}},
 {{\color{red}2},{\color{red}2},{\color{red}2},{1},{2}},
 }
 {\underset s\rightarrow}
 \YT{0.15in}{}{
 {{\color{red}1},{1},{1},{\color{red}1},{\color{red}1},{\color{red}1}},
 {{1},{2},{\color{red}2},{\color{red}2},{\color{red}2}},
 }
 {\underset s\rightarrow}
 \YT{0.15in}{}{
 {{1},{1},{1},{\color{red}1},{\color{red}1},{\color{red}1}},
 {{\color{red}1},{2},{\color{red}2},{\color{red}2},{\color{red}2}},
 }
 {\underset s\rightarrow}
 \YT{0.15in}{}{
 {{1},{1},{1},{\color{red}1},{\color{red}1},{\color{red}1}},
 {2,{\color{red}1},{\color{red}2},{\color{red}2},{\color{red}2}},
 }
\\
&=\rho_1^{(2)}[(Y\cup T)^{[2]}]\\
&=Y_{\nu^{(2)}}\cup H^{(2)},\quad H^{(2)}\equiv Y_{(\mu_1,\mu_2)}
\end{align*}
 \begin{align*}
 \rho_1^{(1)}[(Y\cup T)^{[1]}]&=\YT{0.15in}{}{
 {{1},{1},{\color{red}1},{\color{red}1},{\color{red}1},{\color{red}1}},
}
{\underset{\omega_2^{}}\rightarrow}
\YT{0.15in}{}{
 {{1},{1},{\color{red}1},{\color{red}1},{\color{red}1},{\color{red}1}},
 {{2}},
}
{\underset{\bar\phi_1}\rightarrow}
\YT{0.15in}{}{
 {{1},{1},1,{\color{red}1},{\color{red}1},{\color{red}1}},
 {{2},{\color{red}1}},
}
{\underset{\chi_2^{\mu_2}}\rightarrow}
\YT{0.15in}{}{
 {{1},{1},1,{\color{red}1},{\color{red}1},{\color{red}1}},
 {{2},{\color{red}1},{\color{red}2},{\color{red}2},{\color{red}2}},
}\\
&=\bar\chi_2^{\mu_2}\bar\phi_1\bar\omega_2^{}\rho_1^{(1)}[(Y\cup T)^{[1]}]\\
&=\rho_1^{(2)}[(Y\cup T)^{[2]}]
\end{align*}
\begin{align*}
(Y\cup T)^{[3]}&=\YT{0.15in}{}{
 {{\color{red}1},{\color{red}1},{\color{red}1},{\color{red}1},{1},{1}},
 {{\color{red}2},{\color{red}2},{\color{red}2},{1},{2}},
 {{\color{red}3},{\color{red}3},{1},{2},{3}},
 }
 {\underset s\rightarrow}
 \YT{0.15in}{}{
 {{\color{red}1},{\color{red}1},{1},{1},{\color{red}1},{\color{red}1}},
 {{\color{red}2},{1},{2},{\color{red}2},{\color{red}2}},
 {{1},{2},{3},{\color{red}3},{\color{red}3}},
 }{\underset s\rightarrow}
 \YT{0.15in}{}{
 {{\color{red}1},{\color{red}1},{1},{1},{\color{red}1},{\color{red}1}},
 {1,{1},{2},{\color{red}2},{\color{red}2}},
 {{2},{3},\color{red}2,{\color{red}3},{\color{red}3}},
 }{\underset s\rightarrow}
 \YT{0.15in}{}{
 {{1},{1},{1},{1},{\color{red}1},{\color{red}1}},
 {{2},{2},{\color{red}1},{\color{red}2},{\color{red}2}},
 {3,{\color{red}1}, {\color{red}2},{\color{red}3},{\color{red}3}},
}\\
&=\rho_1^{(3)}[(Y\cup T)^{[3]}]\\
&= Y_{\nu^{(3)}}\cup H^{(3)}\quad H^{(3)}\equiv Y_{(\mu_1,\mu_2,\mu_3)}
 \end{align*}
\begin{align*}
\rho_1^{(2)}[(Y\cup T)^{[2]}]{\underset{\omega_3^{}}\rightarrow}&\YT{0.15in}{}{
 {{1},{1},1,{\color{red}1},{\color{red}1},{\color{red}1}},
 {{2},{\color{red}1},{\color{red}2},{\color{red}2},{\color{red}2}},
 {3},
}
{\underset{\bar\phi_2}\rightarrow}
\YT{0.15in}{}{
 {{1},{1},1,{\color{red}1},{\color{red}1},{\color{red}1}},
 {{2},{2},{\color{red}2},{\color{red}2},{\color{red}2}},
 {3, {\color{red}1}},
}
{\underset{\bar\phi_1}\rightarrow}
\YT{0.15in}{}{
 {{1},{1},{1},1,{\color{red}1},{\color{red}1}},
 {{2},{2},{\color{red}1},{\color{red}2},{\color{red}2}},
 {3,{\color{red}1}, {\color{red}2}},
}
{\underset{\chi_3^{\mu_3}}\rightarrow}
\YT{0.15in}{}{
 {{1},{1},{1},{1},{\color{red}1},{\color{red}1}},
 {{2},{2},{\color{red}1},{\color{red}2},{\color{red}2}},
 {3,{\color{red}1}, {\color{red}2},{\color{red}3},{\color{red}3}},
}
\\
&=\bar\chi_3^{\mu_3}\bar\phi_1\bar\phi_2\bar\omega_3\rho_1^{(2)}(Y\cup T)^{[2]}
=\bar\chi_3^{\mu_3}\bar\phi_1\bar\phi_2\bar\omega_3\bar\chi_2^{\mu_2}\bar\phi_1\bar\omega_2\rho_1^{(1)}[(Y\cup T)^{[1]}]
\\
&=\rho_1^{(3)}[(Y\cup T)^{[3]}]
\end{align*}

\begin{align*}
(Y\cup T)^{[4]}&= \YT{0.15in}{}{
 {{\color{red}1},{\color{red}1},{\color{red}1},{\color{red}1},{1},{1}},
 {{\color{red}2},{\color{red}2},{\color{red}2},{1},{2}},
 {{\color{red}3},{\color{red}3},{1},{2},{3}},
 {{\color{red}4},{2},{3},{4}},
 }
 {\underset s\rightarrow}
 \YT{0.15in}{}{
 {{\color{red}1},{\color{red}1},1,1,{\color{red}1},{\color{red}1}},
 {{\color{red}2},1,2,{\color{red}2},{\color{red}2}},
 {1,2,{3},{\color{red}3},{\color{red}3}},
 {{2},{3},{4},{\color{red}4}},
 }
 {\underset s\rightarrow}
  \YT{0.15in}{}{
 {{\color{red}1},{\color{red}1},1,1,{\color{red}1},{\color{red}1}},
 {1,1,2,{\color{red}2},{\color{red}2}},
 {2,{2},{3},{\color{red}3},{\color{red}3}},
 {{\color{red}2},{3},{4},{\color{red}4}},
 }
  {\underset s\rightarrow}
  \YT{0.15in}{}{
 {1,1,1,1,{\color{red}1},{\color{red}1}},
 {2,{2},2,{\color{red}2},{\color{red}2}},
 {{\color{red}1},{\color{red}1},{3},{\color{red}3},{\color{red}3}},
 {{\color{red}2},{3},{4},{\color{red}4}},
 }\\
 & {\underset s\rightarrow}
  \YT{0.15in}{}{
 {1,1,1,1,{\color{red}1},{\color{red}1}},
 {2,{2},2,{\color{red}2},{\color{red}2}},
 {3,{\color{red}1},{3},{\color{red}3},{\color{red}3}},
 {{\color{red}1},{\color{red}2},{4},{\color{red}4}},
 }
  {\underset s\rightarrow}\YT{0.15in}{}{
 {1,1,1,1,{\color{red}1},{\color{red}1}},
 {2,{2},2,{\color{red}2},{\color{red}2}},
 {3,3,{\color{red}1},{\color{red}3},{\color{red}3}},
 {4,{\color{red}1},{\color{red}2},{\color{red}4}},
 }\\
 &=\rho_1^{(4)}[(Y\cup T)^{[4]}]\\
 &= Y_{\nu^{(4)}}\cup H^{(4)}\quad H^{(4)}\equiv Y_{(\mu_1,\mu_2,\mu_3,\mu_4)}
 \end{align*}

\begin{align*}\rho_1^{(3)}[(Y\cup T)^{[3]}]&{\underset{\omega_4^{}}\rightarrow}
 \YT{0.15in}{}{
{{1},{1},{1},{1},{\color{red}1},{\color{red}1}},
 {{2},{2},{\color{red}1},{\color{red}2},{\color{red}2}},
 {{3},{\color{red}1},{\color{red}2},{\color{red}3},{\color{red}3}},
 {4},
}{\underset{\bar\phi_3}\rightarrow}
\YT{0.15in}{}{
{{1},{1},{1},{1},{\color{red}1},{\color{red}1}},
 {{2},{2},{\color{red}1},{\color{red}2},{\color{red}2}},
 {{3},{3},{\color{red}2},{\color{red}3},{\color{red}3}},
 {4,{\color{red}1}},
}
{\underset{\bar\phi_2}\rightarrow}
\YT{0.15in}{}{
{{1},{1},{1},{1},{\color{red}1},{\color{red}1}},
 {{2},{2},{2},{\color{red}2},{\color{red}2}},
{{3},{3},{\color{red}1},{\color{red}3},{\color{red}3}},
{4,{\color{red}1},{\color{red}2}},
}
{\underset{\chi_4^{\mu_4}}\rightarrow}
\YT{0.15in}{}{
{{1},{1},{1},{1},{\color{red}1},{\color{red}1}},
 {{2},{2},{2},{\color{red}2},{\color{red}2}},
{{3},{3},{\color{red}1},{\color{red}3},{\color{red}3}},
{4,{\color{red}1},{\color{red}2},{\color{red}4}},
}\\
&=\chi_4^{\mu_4}\bar\phi_2\bar\phi_3\bar\omega_4\rho_1^{(3)}[(Y\cup T)^{[3]}]=\bar\chi_4^{\mu_4}\bar\phi_2\bar\phi_3\bar\omega_4\bar\chi_3^{\mu_3}\bar\phi_1\bar\phi_2\bar\omega_3\rho_1^{(2)}(Y\cup T)^{[2]}\\
&=\rho_1^{(4)}[(Y\cup T)^{[4]}]
\end{align*}

\begin{align*}
[Y_{\mu}\cup T]^{[5]}&=\YT{0.15in}{}{
 {{\color{red}1},{\color{red}1},{\color{red}1},{\color{red}1},{1},{1}},
 {{\color{red}2},{\color{red}2},{\color{red}2},{1},{2}},
 {{\color{red}3},{\color{red}3},{1},{2},{3}},
 {{\color{red}4},{2},{3},{4}},
 {{2},{3},{4}},
}
{\underset s\rightarrow} \YT{0.15in}{}{
 {{\color{red}1},{\color{red}1},{\color{red}1},{\color{red}1},{1},{1}},
 {{\color{red}2},{\color{red}2},{\color{red}2},{1},{2}},
 {{\color{red}3},{\color{red}3},{1},{2},{3}},
 {{2},{2},{3},{4}},
 {{\color{red}4},{3},{4}},
}
{\underset s\rightarrow} \YT{0.15in}{}{
 {{\color{red}1},{\color{red}1},{\color{red}1},{\color{red}1},{1},{1}},
 {{\color{red}2},{\color{red}2},{\color{red}2},{1},{2}},
 {{\color{red}3},{\color{red}3},{1},{2},{3}},
 {{2},{2},{3},{4}},
 {{3},{4},{\color{red}4}},
}
{\underset s\rightarrow} \YT{0.15in}{}{
 {{\color{red}1},{\color{red}1},{\color{red}1},{\color{red}1},{1},{1}},
 {{\color{red}2},{\color{red}2},{\color{red}2},{1},{2}},
 {{1},{2},{2},{3},{\color{red}3}},
 {{2},{3},{4},{\color{red}3}},
 {{3},{4},{\color{red}4}},
}\\
&{\underset s\rightarrow} \YT{0.15in}{}{
 {{\color{red}1},{\color{red}1},{\color{red}1},{1},{1},{\color{red}1}},
 {{\color{red}2},{\color{red}2},1,
 {2},{\color{red}2}},
 {{1},{2},{2},{3},{\color{red}3}},
 {{2},{3},{4},{\color{red}3}},
 {{3},{4},{\color{red}4}},
}
{\underset s\rightarrow} \YT{0.15in}{}{
 {{\color{red}1},{\color{red}1},{1},{1},{\color{red}1},{\color{red}1}},
 {{\color{red}2},1,{\color{red}2},
 {2},{\color{red}2}},
 {{1},{2},{2},{3},{\color{red}3}},
 {{2},{3},{4},{\color{red}3}},
 {{3},{4},{\color{red}4}},
 }
 {\underset s\rightarrow} \YT{0.15in}{}{
 {{\color{red}1},{\color{red}1},{1},{1},{\color{red}1},{\color{red}1}},
 {{\color{red}2},1,
 {2},2,{\color{red}2}},
 {{1},{2},3,{\color{red}2},{\color{red}3}},
 {{2},{3},{4},{\color{red}3}},
 {{3},{4},{\color{red}4}},
 }
 {\underset s\rightarrow} \YT{0.15in}{}{
 {{\color{red}1},{\color{red}1},{1},{1},{\color{red}1},{\color{red}1}},
 {1,1,{2},2,{\color{red}2}},
 {{2},{2},3,{\color{red}2},{\color{red}3}},
 {{3},{3},{4},{\color{red}3}},
 {4,{\color{red}2},{\color{red}4}},
 }{\underset s\rightarrow}
 \YT{0.15in}{}{
{{1},{1},{1},{1},{\color{red}1},{\color{red}1}},
 {{2},{2},{2},{2},{\color{red}2}},
{{3},{3},3,{\color{red}2},{\color{red}3}},
{4,4,{\color{red}1},{\color{red}3}},
{{\color{red}1},{\color{red}2},{\color{red}4}},
}
\\
 &=\rho_1^{(5)}[(Y\cup T)^{[5]}]\\
 &= Y_{\nu}\cup H\quad H\equiv Y_{\mu}
\end{align*}

\begin{align*}\rho_1^{(4)}[(Y\cup T)^{[4]}]&{\underset{\bar\phi_4}\rightarrow}
\YT{0.15in}{}{
{{1},{1},{1},{1},{\color{red}1},{\color{red}1}},
 {{2},{2},{2},{\color{red}2},{\color{red}2}},
{{3},{3},{\color{red}1},{\color{red}3},{\color{red}3}},
{4,4,{\color{red}2},{\color{red}4}},
{{\color{red}1}},
}{\underset{\bar\phi_3}\rightarrow}
\YT{0.15in}{}{
{{1},{1},{1},{1},{\color{red}1},{\color{red}1}},
 {{2},{2},{2},{\color{red}2},{\color{red}2}},
{{3},{3},3,{\color{red}3},{\color{red}3}},
{4,4,{\color{red}1},{\color{red}4}},
{{\color{red}1},{\color{red}2}},
}
{\underset{\bar\phi_2}\rightarrow}
\YT{0.15in}{}{
{{1},{1},{1},{1},{\color{red}1},{\color{red}1}},
 {{2},{2},{2},{2},{\color{red}2}},
{{3},{3},3,{\color{red}2},{\color{red}3}},
{4,4,{\color{red}1},{\color{red}3}},
{{\color{red}1},{\color{red}2},{\color{red}4}},
}\\
&=\bar\phi_2\bar\phi_3\bar\phi_4\rho_1^{(4)}[(Y\cup T)^{[4]}]\\
&=\rho_1^{(5)}[Y\cup T]\\
&=Y_\nu\cup {H} \quad H\equiv Y_{\mu}.
\end{align*}

Observe that $\widehat \nu_1=\nu^{(1)}=(2)\subseteq \nu^{(2)}=(3,1)\subseteq \nu^{(3)}=(4,2,1)\subseteq \nu^{(4)}=(4,3,2,1)\subseteq \nu^{(5)}=(4,4,3,2,0)$ defines $G_\nu(T)$.
\end{ex}

The switching map on ballot tableau pairs of normal shape can now be translated to the language of internal row insertion operations.

\begin{thm} (Theorem \ref{2}) \label{th:main} (Main Theorem)
Let $n\ge 1$ and $Y\cup T\in {\cal{LR}}^{(n)}$ with $Y=Y_\mu$  and $T$ a ballot tableau of shape $\lambda/\mu$ and weight $\nu$.
For $1\le i\le n$, let $(Y_\mu\cup T)^{[i]}=Y_{(\mu_1,\dots,\mu_i)}\cup T^{[i]}\in {\cal{LR}}^{(i)}$ 
with $T^{[i]}$ of weight $\nu^{(i)}$. Consider the $i$th row word of $T^{[i]}$ where $V_i$ is the  row subword  restricted to the  entries in $[i-1]$, and  $\widehat \nu_i=\lambda_i-\mu_i-|V_i|$ is the number of entries equal to $i$. Put  $(Y_\mu\cup T)^{[0]}=Y_{\nu^{0}}\cup H^{(0)}:=\emptyset$, $\nu^{0}:=0$,
$\bar\phi_{\emptyset}=id$  and  $\rho_1^{(0)}(\emptyset):=\emptyset$.
Then, for $i=1,\dots,n$, it holds
\begin{align}
\rho_1^{(i)}[(Y_\mu\cup T)^{[i]}]&=(\bar\chi_i^{\mu_i}\circ\bar\phi_{V_i}\circ\bar\omega_i^{\widehat\nu_i})\circ\rho^{(i-1)}[(Y_\mu\cup T)^{[i-1]}]
\label{introd:mainrecursionxx}
\\
&=\bar\chi_i^{\mu_i}\circ\bar\omega_i^{\widehat\nu_i}\circ\bar\phi_{V_i}(Y_{\nu^{(i-1)}}\cup H^{(i-1)})
=Y_{\nu^{(i)}}\cup H^{(i)}\in {\cal{LR}}^{(i)},\label{introd:mainrecursionx}
\end{align}
where $\bar\omega_i^{\widehat\nu_i}$ adds the $i$th row word $i^{\widehat\nu_i}$ to $ Y_{\nu^{(i-1)}}$, $\bar\chi_i^{\mu_i}$ adds the row word $i^{\mu_i}$ at the end  of the $i$th row of $\bar\phi_{V_i}\circ\bar\omega_i^{\widehat\nu_i}(Y_{\nu^{(i-1)}}\cup H^{(i-1)})$  and  $H^{(i)}\equiv Y_{(\mu_1,\dots,\mu_i)}$.
In particular, all   bumping routes of $\bar\phi_{V_i}$ are pairwise disjoint and terminate in the $i$th row.
\end{thm}

\subsection {Lecouvey-Lenart and Kumar-Torres bijections between  Kwon and Sundaram branching models}\label{conjecture}
For detailed definitions pertaining this section we refer to \cite{leclen,krv21,kwon18, sathishtorres}.
A fundamental fact of Kumar-Torres bijection  is that  it restricts to  tableaux satisfying the Sundaram condition and those whose evacuation satisfy the Kwon condition by considering and recognizing that they can be embedded in  the 
Kushwaha–Raghavan–Viswanath \cite{krv21} bijection on flagged hives. To settle the Lecouvey-Lenart conjecture \cite{leclen}, it remains to know  the coincidence of  LR commuters. Note that in our notation the role of $\lambda$ and  $\nu$ are swapped in \cite{sathishtorres}.

 We  fix a positive integer
$n$, and assume, unless otherwise stated, that $\ell(\lambda) \le 2n -1$ and $\ell(\mu)\le n$. Let $m=2n$.
A Littlewood–Richardson tableau of shape $\lambda/\mu$ and weight $\nu$ satisfies the
\emph{Sundaram property} if for each $i = 0,\dots , \ell(\nu)/2$, the entry $2i +1$ appears in row $n+i$
or above in the Young diagram of $\lambda$. The set of $T \in  LR(\lambda/\mu, \nu)$ satisfying
the Sundaram property is denoted by $LRS(\lambda/\mu, \nu)$ and called the set of \emph{ Sundaram LR tableaux}.

 Denote by $LRS_{\mu,\nu}^\lambda$  the subset  of $LR_{\mu,\nu}^\lambda$ consisting of the right companions of $LRS(\lambda/\mu, \nu)$. Then those tableaux consist of the tableaux in $LR_{\mu,\nu}^\lambda$ satisfying the flag property that the entries in the $kth$ row are bounded above by $n+ \left \lfloor k/2\right \rfloor$, for $k=1,\dots,2n$, \cite[Proposition 4.7]{sathishtorres}. We call them   \emph{Sundaram right companions}. Denote by ${}^-LRS_{\mu,\nu}^\lambda$  the subset  of ${}^-LR_{\mu,\nu}^\lambda$ consisting of the left companions of $LRS(\lambda/\mu,\nu)$.

A semistandard tableau of shape $\mu$ and $\ell(\mu)\le n $ is said to satisfy the \emph{Kwon
property} if the entries in row $i$ are at least $2i - 1$, for $i = 1,\dots, n$. Denote by $LRK_{\nu,\mu}^\lambda$ the subset
of $LR^\lambda_{\nu,\mu}$ consisting of tableaux $G$
such that their {Sch\"utzenberger evacuation} $\xi=evac_{2n}$ within the crystal $B(\mu,2n)$, $\xi(G)=evac_{2n}(G )$, satisfy the Kwon
property.

The Sundaram branching rule (see \cite{sathishtorres} and \cite{sundaram} for the definition)  states the following.

\begin{thm}\cite{sundaram} The branching coefficient
$c_{\mu}^\lambda$
equals the cardinality
of the set
 $$LRS(\lambda,\mu):=\bigcup LRS(\lambda/\mu,\nu),$$ where the union is taken over all even partitions $\nu$, that is, $\nu_{2i-1}=\nu_{2i}$, $i\ge 1$.
\end{thm}

The Kwon’s branching rule as reformulated by Lecouvey–
Lenart \cite[ Lemma 6.11]{leclen} says the following.

\begin{thm}\cite{kwon18,leclen} The branching coefficient
$c_{\mu}^\lambda$
equals the cardinality
of the set
 $$LRK(\lambda,\mu):=\bigcup LRK^\lambda_{\nu,\mu}$$ where the union is taken over all even partitions $\nu$, that is, $\nu_{2i-1}=\nu_{2i}$, $i\ge 1$.
\end{thm}

\begin{thm}[Kumar-Torres  Theorem 3.6 \cite{sathishtorres}] The composition
$$LR(\lambda/\mu, \nu)\overset{\sim}\longrightarrow  LR_{\mu,\nu}^\lambda
\overset{U}\longrightarrow LR_{\nu,\mu}^\lambda$$
restricts to a bijection
\begin{align}LRS(\lambda/\mu, \nu)\overset{\sim}\longrightarrow  LRK^\lambda_{\nu, \mu}\end{align}
where $U : LR_{\mu,\nu}^\lambda
\overset{U}\longrightarrow LR_{\nu,\mu}^\lambda$ is the LR commuter by Kushwaha–Raghavan–Viswanath\cite{krv21}.
(See \cite{sathishtorres} for definitions). Therefore, the above composition
induces a bijection between $ LRS(\lambda,\mu)$ and $LRK(\lambda,\mu)$.
\end{thm}

From \cite[Section 12]{akt16} we know the LR commuter $U$ coincides with the Henriques-Kamnitzer $\mathfrak{gl}_n$-crystal commuter because it produces the same GT pattern pair. The Kumar-Torres bijection together with the coincidence of LR commuters  gives the corollary:
\begin{cor}\label{cor:llkt} The Kumar-Torres bijection
\begin{align*}LR(\lambda/\mu, \nu)\overset{\sim}\longrightarrow  LR_{\mu,\nu}^\lambda
\overset{U}\longrightarrow LR_{\nu,\mu}^\lambda\end{align*}
and the Lecouvey-Lenart bijection
\begin{align*}LR(\lambda/\mu, \nu)\overset{\sim}\longrightarrow  LR_{\mu,\nu}^\lambda
\overset{U'}\longrightarrow LR_{\nu,\mu}^\lambda\end{align*}
where $U'$ is the LR commuter defined by Henriques–Kamnitzer, coincide. Thereby both restrict to a bijection $$LRS(\lambda, \mu)\overset{\sim}\longrightarrow  LRK(\lambda, \mu).$$
\end{cor}

\begin{cor} \label{kwon-sundaram-pair} The Henriques-Kamnitzer commuter \eqref{comhk} or \eqref{comhk-} restrict to LRS tableaux and gives
\begin{align}&LRS_{\mu,\nu}^\lambda\overset{\sim}\longrightarrow LRS(\lambda/\mu,\nu)\overset{\sim}\longrightarrow {}^-LRS_{\mu,\nu}^\lambda \overset{\xi}\longrightarrow LRK_{\nu,\mu}^\lambda,\; G_\nu\mapsto T\mapsto G_\mu\mapsto \xi(G_\mu),
\end{align}

Therefore, ${}^-LRS_{\mu,\nu}^\lambda=\xi(LRK_{\nu,\mu}^\lambda)$ which are precisely  those  tableaux satisfying the Kwon condition in ${}^-LR_{\mu,\nu}^\lambda$.
\end{cor}


\begin{ex}\label{llkt}Motivated by  the question raised by the authors in \cite[Remark 3.7]{sathishtorres}, we now illustrate the Lecouvey-Lenart and Kumar-Torres bijections with our LR commuter based  on the Sagan-Stanley internal insertion. We consider \cite[Example 4.11]{sathishtorres} where $n=3$, $m=6$, $\ell(\lambda)=4$, $\ell(\nu)=4$,  $\ell(\mu)=3$, and $\mu=(2,1,1,0,0,0),$ $\nu=(4,4,2,1,0,0)$ $\subseteq $ $\lambda=(5,4,3,3,0,0)$:

\begin{align}\label{sator}Y_{\mu}\cup T=\YT{0.15in}{}{
 {{\color{red}1},{\color{red}1},{1},{1},1},
 {{\color{red}2},1,2,{2}},
 {{\color{red}3},{2},{3}},
 {{2},{3},{4}},
}, \quad T\in LRS(\lambda/\mu,\nu)\subseteq LR(\lambda/\mu,\nu).\end{align}

The right and left companions of the Sundaram LR tableau $T$ are respectively
\begin{align}
G_\nu(T)=\YT{0.15in}{}{
 {{1},{1},1,2},
 {2,{2},3,4},
 {{3},4},
 {{4}},
}\in LRS_{\mu,\nu}^\lambda,
\quad  G_\mu(T)=\YT{0.15in}{}{
 {{3},{3}},
 {4},
 {{6}},
}\in {}^-LRS_{\mu,\nu}^\lambda\subseteq {}^-LR_{\mu,\nu}^\lambda
\end{align}
where $G_\mu(T)$ of weight $rev(\lambda-\nu)=(0,0,2,1,0,1))$ is defined by the nested sequence $$(2,1,1,0,0,0)\supseteq (21000)\supseteq (2100)\supseteq(200)\supseteq(00)\supseteq (00).$$
We may check that
$G_\mu(T)$ satisfy the Kwon condition:  the entries in row $i$ are at least $2i - 1$, for $i = 1,2,3$: entries in row 1 are $3\ge 1$, entries in row 2 are $4\ge 3$, entries in row 3 are $6\ge 6-1$.

\begin{align} evac_{6}\, G_\mu(T)=\YT{0.15in}{}{
 {{1},{4}},
 {3},
 {{4}},
}\in LRK_{\nu,\mu}^\lambda\subseteq  LR_{\nu,\mu}^\lambda  \subseteq B(\mu,6),\quad \mbox{weight $\lambda-\nu=(1,0,1,2,0,0)$}
\end{align}

 The Henriques-Kamnitzer LR commuter, the LR commuter by Kushwaha–Raghavan–Viswanath and our commuter all of them send $(G_\mu(T), G_\nu(T))$ to $(evac_6 \,G_\nu(T),evac_6\, G_\mu(T))$:

 \begin{align}\label{azesator}Y_{\mu}\cup T=\YT{0.15in}{}{
 {{\color{red}1},{\color{red}1},{1},{1},1},
 {{\color{red}2},1,2,{2}},
 {{\color{red}3},{2},{3}},
 {{2},{3},{4}},
}\rightarrow
\YT{0.15in}{}{
 {{\color{red}1},{\color{red}1},{1},{1},1},
 {{\color{red}2},1,2,{2}},
 {{\color{red}3},{2},{3}},
 {{2},{3}},
}\rightarrow
\YT{0.15in}{}{
 {{\color{red}1},{1},{1},{1},1},
 {{\color{red}2},2,2,{2}},
 {{\color{red}3},{3},{3}},
 {{2}},
}\rightarrow
\YT{0.15in}{}{
 {{\color{red}1},{1},{1},{1},1},
 {{\color{red}2},2,2,{2}},
 {2,{3},{3}},
}\end{align}
\begin{align}
\rightarrow
\YT{0.15in}{}{
 {{\color{red}1},{1},{1},{1},1},
 {{\color{red}2},2,2,{2}},
 {2},
}
\rightarrow \YT{0.15in}{}{
 {{\color{red}1},{1},{1},{1},1},
 {{2},2,2,{2}},
}\rightarrow
\YT{0.15in}{}{
 {{\color{red}1},{1},{1},{1},1},
}\rightarrow
\YT{0.15in}{}{
 {{\color{red}1}},
}\rightarrow \emptyset
\end{align}

One then has the GT pattern of type $\mu$ defined by the red nested sequence of partitions
$$(1)\subseteq (1)\subseteq (1,1)\subseteq (2,1,1),$$
 that  gives the tableau $evac_6\,G_\mu(T)$ as in the bijection by Kumar-Torres based on the LR commuter by Kushwaha–Raghavan–Viswanath,
\begin{align*}evac_6\,G_\mu(T)=\YT{0.15in}{}{
 {{1},{4}},
 {3},
 {{4}},
}\overset{c^{-1}}\longrightarrow  Y_\nu\cup H=\YT{0.15in}{}{
 {{\color{red}1},{\color{red}1},{\color{red}1},{\color{red}1},{1}},
 {{\color{red}2},{\color{red}2},{\color{red}2},{\color{red}2}},
 {{\color{red}3},{\color{red}3},{2}},
{{\color{red}4},1,3},
} \\
 G_\nu(H)=\YT{0.15in}{}{
 {{3},3,3,{5}},
 {4,4,6,6},
 {5,5},
{6 },
}=evac_6(G_\nu(T))\qquad G_\mu(H)=evac_{6}\,G_\mu(T) \in LRK^\lambda_{\mu,\nu}
\end{align*}
\end{ex}


\section{Recursion of switching on ballot tableau pairs} \label{recursionswitch}
 Switching can be performed in stages whenever a decomposition of  $Y$ or $T$ in  the tableau pair $Y\cup T$ is considered.
  This property allows to exhibit a recursion of the switching map on ballot tableau pairs $Y_\mu\cup T$ by reducing the size of the partition $\mu$. { Due to the switching condition \eqref{switch}, switches in a tableau pair $Y\cup T$  where $Y$ or $T$ is a ballot tableau are such that the entries $i$ in the $i$th row of a ballot tableau can not be switched upwards.  Thereby switching  $T$  with the rows of $Y_\mu$ in stages incurs that  the length of the word in the $n$th row of $T$ restricted to the alphabet $[n-1]$ eventually reduces. Because an entry $i$ in row $i$ of $Y_\mu$ either switches horizontally with an entry of $T$ or vertically  with an  entry of $T$ below row $i$ and further moves of that $i$ will never occur with entries of $T$ above row $i$.} Recall Remark \ref{ob:factorLR}.

\begin{lem}\label{lem:2-0}
Let $2\le l\le n$ and  $T\in YT((\lambda_l,\dots,\lambda_{n+1})/(\mu_l))$, $\mu_l>0$, of weight $(\alpha_1,\dots,\alpha_{l-1},\nu_l,\dots,\nu_{n+1})$, with $(\nu_l,\dots,\nu_{n+1})$ a partition. For some $\beta\subseteq (\lambda_l,\dots,\lambda_{n+1})$, assume the decomposition $Y_{(\mu_l)}\cup T=Y_{(\mu_l)}\cup A\cup B$ where $A=T_{|[l-1]}\in YT(\beta/(\mu_l))$ has weight $(\alpha_1,\dots,\alpha_{l-1})$, and $B=T_{|[l,n+1]} $
is a ballot tableau of shape $(\lambda_l,\dots,\lambda_{n+1})/\beta$ and weight $(\nu_l,\dots,\nu_{n+1})$. Then if $F(n+1)^{\nu_{n+1}}$ is the $(n-l+2)$th row of $T$ with $F$ a word in $[n]$, one has
$$\rho_1^{(n-l+2)}(Y_{(\mu_l)}\cup T)=S\cup Q,$$
such that $S=S_1\cup S_2\equiv T$, with $S_1=S_{|[l-1]}\equiv A$ and $S_2=S_{|[l,n+1]}\equiv B\equiv Y_{(\nu_l,\dots,\nu_{n+1})}$, and $Q\equiv Y_{(\mu_l)}$.
Moreover, if $\hat F(n+1)^{\nu_{n+1}}$ is the $(n-l+2)$th row of $S$
 and  $D$ is the $(n-l+2)$th row of  $Q$ then  $\hat F$ is a subword of $F$ and $|\hat F|+|D|=|F|$.
\end{lem}
\begin{proof} Let $S$ be the rectification of $T$. The rectification can proceed in stages by switching $Y_{(\mu_l)}$ with $A$ and $B$ in stages.   Observe  that  since $T$ restricted to the alphabet $[l,n+1]$ is the ballot tableau $B$ of skew shape $(\lambda_l,\dots,\lambda_{n+1})/\beta$ and weight  $(\nu_l,\dots,\nu_{n+1})$, the  $\nu_{n+1}$ entries $n+1$ in   the last row of $T$ stay there until the end of the rectification of $T$.  Therefore $\rho_1^{(n-l+2)}(Y_{(\mu_l)}\cup T)=S\cup Q,$ such that $S=S_1\cup S_2\equiv T$, with $S_1=S_{|[l-1]}\equiv A$ and $S_2=S_{|[l,n+1]}\equiv B\equiv Y_{(\nu_l,\dots,\nu_{n+1})}$, and $Q\equiv Y_{(\mu_l)}$. In particular,  the last row of $S$ has $\nu_{n+1}$ entries $n+1$.

 If $F=\emptyset$ then   $Y_{(\mu_l)}\cup T=Y_{(\mu_l)}\cup [(n+1)^{\nu_{n+1}}\ast T^{[n-l+1)]}]=Y_{(\mu_l)}\cup [A\cup\big((n+1)^{\nu_{n+1}}\ast  B^{[n-l+1)]}\big)]$. Hence  $\rho_1^{(n-l+2)}(Y_{(\mu_l)}\cup T)=[S_1\cup \big((n+1)^{\nu_{n+1}}\ast S_2^{[n-l+1]}\big)]\cup (\emptyset\ast Q^{[n-l+1]}) $ with $S_1\equiv A$, $S_2=(n+1)^{\nu_{n+1}}\ast S_2^{[n-l+1]}\equiv Y_{(\nu_l,\dots,\nu_{n+1})}$ and $S=S_1\cup S_2\equiv T$.   Thus
 $$\begin{array}{cccccc}\rho_1^{(n-l+2)}(Y_{(\mu_l)}\cup T)
 &=&(n+1)^{\nu_{n+1}}\ast\rho_1^{(n-l+1)}(Y_{(\mu_l)}\cup T^{[n-l+1]})\nonumber\\
 &=&((n+1)^{\nu_{n+1}}\ast S^{[n-l+1]})\cup (\emptyset\ast Q^{[n-l+1]}),
 \end{array}$$
where $S=((n+1)^{\nu_{n+1}}\ast S^{[n-l+1]})\equiv T$ and $Q=\emptyset\ast Q^{[n-l+1]}\equiv Y_{(\mu_l)}$.  The last rows, $\hat F$ of  $S$, and $D$ of $Q$ are both the empty word.

Let $F\neq\emptyset$ and $Y_{(\mu_l)}=l^{\mu_l}$. Either an entry $l$ of $Y_{(\mu_l)}$ reaches  in  the switching process the row next to the last or not. In the later case $D=\emptyset$ and $\widehat F=F$. In the former case,   there exists  a perforated tableau pair where the row next to the last has one entry of $Y_{(\mu_l)}$. Then
 the two last rows  are either of the form:

 $(a)$
$$\begin{matrix}\cdots& {\color{red}l}&a&\cdots&n^{\hat\nu_n}\\
\cdots&b&\cdots&(n+1)^{\nu_{n+1}}&\\
\end{matrix}
,\;\text{with $1\le a, b\le n$ and $\hat\nu_n\le \nu_n$},$$  and  the next switches either are
$$\begin{matrix}\cdots& {\color{red}l}&a&\cdots&n^{\hat\nu_n}\\
\cdots&b&\cdots&(n+1)^{\nu_{n+1}}\\
\end{matrix}\quad{\underset{\bf s}\rightarrow}\quad\begin{matrix} \cdots&b&a&\cdots&n^{\hat\nu_n}\\
 \cdots&{\color{red}l}&\cdots&(n+1)^{\nu_{n+1}}\\
\end{matrix}\quad{\underset{\bf s}\rightarrow}\quad\\
$$

$$
\qquad\qquad\qquad\qquad\qquad\qquad\qquad\qquad\qquad\qquad\qquad{\underset{\bf s} \rightarrow}\quad \begin{matrix}\cdots&b&a&\cdots&n^{\hat\nu_n}\\
&\cdots &(n+1)^{\nu_{n+l}}&\color{red}l\\
\end{matrix},\quad  \text{if $b\le a\le n$},$$

or

$$\qquad\qquad\qquad\quad\begin{matrix}\cdots& \color{red}l&a&\cdots&n^{\hat\nu_n}\\
\cdots&b&\cdots &(n+1)^{\nu_{n+1}}\\
\end{matrix} \qquad{\underset{\bf s}\rightarrow}\qquad \begin{matrix}\cdots& a&\color{red}l&\cdots&n^{\hat\nu_n}\\
\cdots& b&\cdots &(n+1)^{\nu_{n+1}}\\
\end{matrix},\quad  \text{if $a<b \le n$};$$

\noindent or $(b)$   $$\begin{matrix}\cdots& &{\color{red}l}&\cdots &n &n&\cdots n\\
\cdots&n+1\cdots& n+1&\cdots &n+1&\\
\end{matrix}\quad{\underset {\text{horizontal switch {\bf s}}}\rightarrow}\quad
\begin{matrix}\cdots&& n^{\hat\nu_n}&{\color{red}l}\\
\cdots&(n+1)^{\nu_{n+1}}&\\
\end{matrix}.$$

In any case, $S=[\hat F(n+1)^{\nu_{n+1}}]\ast S^{[n-l+1]}$ with $\hat F$ a subword of $F$, and $ Q=D\ast Q^{[n-l+1]}\equiv Y_{(\mu_l)}$ such that $D=l^{|D|}$ and  $|\hat F|+|D|=|F|$.
\end{proof}
\begin{ex}Below one 
illustrates the previous lemma.
\begin{enumerate}
\item $l=3<n=4$, $\mu_3=3$, $T\in YT((7,6,4)/(3))$ of content $(2,1; 5,3,3)$ such that $T_{|\{1,2\}}\in YT((4,2)/(3))$ has content $(\alpha_1=2,\alpha_2=1)$  and  $T_{|\{3,4,5\}}$ a ballot tableau of skew-shape $(7,6,4)/(4,2)$ content $\nu=(\nu_3=5,\nu_4=3,\nu_5=3)$
and $F=3$.
$$Y_{(3)}\cup T=\YT{0.14in}{}{
 {\color{red}3,\color{red}3,\color{red}3,1,3,3,3},
 {1,2,3,4,4,4},
 {3,5,5,5},
 }
 \underset {s}\rightarrow
 \YT{0.14in}{}{
 {\color{red}3,\color{red}3,1,3,3,3,\color{red}3},
 {1,2,3,4,4,4},
 {3,5,5,5},
 }
 \underset {s}\rightarrow
 \YT{0.14in}{}{
 {\color{red}3,1,\color{red}3,3,3,3,\color{red}3},
 {1,2,3,4,4,4},
 {3,5,5,5},
 }
 $$
 $$
 \underset {s}\rightarrow
 \YT{0.14in}{}{
 {\color{red}3,1,3,3,3,3,\color{red}3},
 {1,2,\color{red}3,4,4,4},
 {3,5,5,5},
 }
\underset {s}\rightarrow
 \YT{0.14in}{}{
 {\color{red}3,1,3,3,3,3,\color{red}3},
 {1,2,4,4,4,\color{red}3},
 {3,5,5,5},
 }
 \underset {s}\rightarrow
 \YT{0.14in}{}{
 {1,1,3,3,3,3,\color{red}3},
 {\color{red}3,2,4,4,4,\color{red}3},
 {3,5,5,5},
 }
 $$
 $$
 \underset {s}\rightarrow
 \YT{0.14in}{}{
 {1,1,3,3,3,3,\color{red}3},
 {2,4,4,4,\color{red}3,\color{red}3},
 {3,5,5,5},
 }=S\cup Q,\; S\equiv T, \; Q=\emptyset\ast Q^{[2]}\equiv Y_{(3)},\, D=\emptyset,\, F=\widehat F,
 \;\;\widehat F5^{\nu_5}=F\,5^{\nu_5}=3\,5^3.$$

\item
$l=3<n=4$, $\mu_3=3$, $H\in YT((7,6,4)/(3))$ of content $(2,1; 5,4,2)$ such that $H_{|\{1,2\}}\in YT((4,2)/(3))$ has content $(\alpha_1=2,\alpha_2=1)$  and  $H_{|\{3,4,5\}}$ a ballot tableau of skew-shape $(7,6,4)/(4,2)$ content $\nu=(\nu_3=5,\nu_4=4,\nu_5=2)$
and $F=34$.
 $$Y_{(3)}\cup H=\YT{0.14in}{}{
 {\color{red}3,\color{red}3,\color{red}3,1,3,3,3},
 {1,2,3,4,4,4},
 {3,4,5,5},
 }
 \underset {s}\rightarrow
 \YT{0.14in}{}{
 {\color{red}3,\color{red}3,1,3,3,3,\color{red}3},
 {1,2,3,4,4,4},
 {3,4,5,5},
 }
 \underset {s}\rightarrow
 \YT{0.14in}{}{
 {\color{red}3,1,\color{red}3,3,3,3,\color{red}3},
 {1,2,3,4,4,4},
 {3,4,5,5},
 }$$
 $$
 \underset {s}\rightarrow
 \YT{0.14in}{}{
 {\color{red}3,1,3,3,3,3,\color{red}3},
 {1,2,\color{red}3,4,4,4},
 {3,4,5,5},
 }
 \underset {s}\rightarrow
 \YT{0.14in}{}{
 {\color{red}3,1,3,3,3,3,\color{red}3},
 {1,2,4,4,4,\color{red}3},
 {3,4,5,5},
 }
 \underset {s}\rightarrow
 \YT{0.14in}{}{
 {1,1,3,3,3,3,\color{red}3},
 {\color{red}3,2,4,4,4,\color{red}3},
 {3,4,5,5},
 }
 \underset {s}\rightarrow
 \YT{0.14in}{}{
 {1,1,3,3,3,3,\color{red}3},
 {2,\color{red}3,4,4,4,\color{red}3},
 {3,4,5,5},
 }$$
 $$
 \underset {s}\rightarrow
 \YT{0.14in}{}{
 {1,1,3,3,3,3,\color{red}3},
 {2,4,4,4,4,\color{red}3},
 {3,\color{red}3,5,5},
 }
 \underset {s}\rightarrow
 \YT{0.14in}{}{
 {1,1,3,3,3,3,\color{red}3},
 {2,4,4,4,4,\color{red}3},
 {3,5,5,\color{red}3},
 }=S\cup Q,\; S\equiv H,\; Q=3\ast Q^{[2]},\; D=3, \, F=34,\,\widehat F=3$$
$$\text{with $\widehat F5^{\nu_5}=3\,5^2$ a subword of $ F5^{\nu_5}=34\,5^2$}.$$
 \end{enumerate}
\end{ex}
\begin{lem}\label{lem:2-1} Let $n\ge 1$ and $Y_\mu\cup T\in {\cal{LR}}^{(n+1)}$,  with  $\mu=(\mu_1,\dots,\mu_n,0)$  a non zero partition and  $T\equiv Y_\nu$. Suppose  $\lambda_{n+1}=\nu_{n+1}\ge 0$, that is, the $(n+1)$th row of $T$ is the word $(n+1)^{\nu_{n+1}}$. Then $$\rho_1^{(n+1)}(Y_\mu\cup T)=(n+1)^{\nu_{n+1}}\ast\rho_1^{(n)}(Y_\mu\cup T^{[n]})= Y_\nu\cup Q,$$ where $Q=\emptyset\ast Q^{[n]}\equiv Y_\mu$.

\end{lem}
\begin{proof} One has  $Y\cup T=(n+1)^{\nu_{n+1}}\ast(Y\cup T^{[n]})$ a ballot tableau pair. Hence  the switching procedure on $Y\cup T$  only comprises  the entries of $T^{[n]}$, all $\le n$, and the entries of $Y_\mu$. Thus $\rho_1^{(n+1)}(Y_\mu\cup T)=(n+1)^{\nu_{n+1}}\ast\rho_1^{(n)}(Y_\mu\cup T^{[n]})=Y_\nu\cup Q$ and  the $(n+1)$th row of $Q$ is empty.
\end{proof}
Let $n\ge 1$ and $Y_\mu\cup T\in {\cal{LR}}^{(n+1)}$,  with   $\mu=(\mu_1,\dots,\mu_n,0)$  a non zero partition, and  $T\equiv Y_\nu$. Next theorem  uses switching into stages. For  some  $1\le d\le n$ with $\mu_d>0$, we decompose $Y_\mu=Y_{(\mu_1,\dots,\mu_{d-1})}\cup Y_{(\mu_d,\dots,\mu_n)}$, and  thereby decomposing $Y_\mu\cup T= Y_{(\mu_1,\dots,\mu_{d-1})}\cup Y_{(\mu_d,\dots,\mu_n)}\cup T$. Then switch $T$ with $Y_{(\mu_d,\dots,\mu_n)}$  to get $Y_{(\mu_1,\dots,\mu_{d-1})}\cup S\cup Q$, with $S\equiv T$ and  $Q\equiv Y_{(\mu_d,\dots,\mu_n)}$ consisting of the entries of $Y_\mu$ moved  out to the skew shape of $T$. The choice of $d$ is made with the purpose to reduce   the length of the $(n+1)$th row word $F$ of $T$, restricted to the alphabet $[n]$.
We have the guarantee that this happens with $d=1$ but at this point, when $T$ is rectified,  $ F$ is empty. (Note that when $d=1$ we  have the full switch of $Y$ with $Y_\mu$ which gives
$S\cup Q$, with $S=Y_\nu\equiv T$, therefore the $(n+1)$th row of $S$ is empty, and  $Q\equiv Y_{(\mu_1,\dots,\mu_n)}$.) The choice of $d$ is made when for the first time an entry $d$ of $Y_\mu$ reaches the $(n+1)$th row.   When this happens and $T$ is full switched with $Y_{(\mu_d,\dots,\mu_n)}$, we stop the switching.  At this stage the   $(n+1)$th row  $D$ of $Q$ comprises only  entries equal to $d$ and the   word $F$ is reduced to a subword $\widehat F$ of length $|F|-|D|$ with $|D|>0$.

 \begin{thm} \label{recursion}Let $n\ge 1$ and $Y_\mu\cup T\in {\cal{LR}}^{(n+1)}$,  with   $\mu=(\mu_1,\dots,\mu_n,0)$  a non zero partition, and  $T\equiv Y_\nu$. Suppose    $\lambda_{n+1}-\nu_{n+1}\ge 1$, that is, the $(n+1)$th-row  of $T$ is the word $F(n+1)^{\nu_{n+1}}$ with $F $  a non empty word   in the alphabet $[n]$. Then, there exists
 $1\le d\le n$, with $\mu_d> 0$, such that
\begin{equation}\label{eq:switch0}
Y_\mu\cup T\underset{s}\rightarrow Y_{(\mu_1,\dots,\mu_{d-1})}\cup S\cup Q=[\widehat F(n+1)^{\nu_{n+1}}\ast(Y_{(\mu_1,\dots,\mu_{d-1})}\cup S^{[n]})]\cup [D\ast Q^{[n]}]
 \end{equation}
 where  $S\equiv Y_\nu$  and  $Q={ D\ast Q^{[n]}}\equiv Y_{(\mu_d,\dots,\mu_n)}$ is over the alphabet $[{d},{n}]$, with $Q^{[d-1]}=\emptyset$, and $D=d^{|D|}$.
 In addition, the $(n+1)$th row  $\widehat F(n+1)^{\nu_{n+1}}$ of $S$ is such that  $\widehat F$ is a strict subword of $F$ whose length  satisfies $ |D|+|\widehat F|=|F|>|\widehat F|\ge 0$.
Also,\begin{equation}\label{eq:x1recur}
 \rho_1^{(n+1)}(Y_\mu\cup T)=\begin{cases}
Y_\nu\cup Q,\;\text{with}\;|D|=\lambda_{n+1}-\nu_{n+1},\;\text{if $d=1$},\\
\rho_1^{(n+1)}[Y_{(\mu_1,\dots,\mu_{d-1})}\cup S]\cup Q,\;\;\text{if $d>1$.}
\end{cases}
\end{equation}
 \end{thm}

 \begin{proof}
 We handle the proof  by induction on the length  $l:=\ell(\mu)\ge 1$. Let $n\ge l$ and   $v:=\lambda_{n+1}-\nu_{n+1}=|F|>0$  the number of entries $\le n$ in the $(n+1)$th row of $T$.

 For $l=1$, $Y_{(\mu_1)}\cup T\in {\cal{LR}}^{(n+1)}$ with $\mu_1>0$, and, therefore, the $(n+1)$th row of $T$  is of the form $n^v(n+1)^{\nu_{n+1}}$
 with $1\le v:=\lambda_{n+1}-\nu_{n+1}\le \mu_1$.
 Rectify $T$ with {\em jeu de taquin} slides looking at the entries of $Y_{(\mu_1)}$ as holes. Then $$Y_{\mu_1}\cup T\underset{s}\rightarrow Y_\nu\cup Q,$$
where $Q=D\ast Q^{[n]}\equiv Y_{\mu_1}$. Since the shape of $Y_{(\mu_1)}\cup T$ is preserved in the switching procedure, $|D|=|F|=v$ and $D={ 1}^v$. In the case $l=1$, \eqref{eq:switch0} and \eqref{eq:x1recur} hold with $d=1$, and $\widehat F=\emptyset$.

 Let $l>1$, and assume the statement true for $1,\dots, l-1$. Then $n+1\ge l+1\ge 3$, and  consider the factorisation $$Y_{(\mu_1,\dots,\mu_l)}\cup T=[Y_{(\mu_l)}\cup \widehat T]\ast[Y_{(\mu_1,\dots,\mu_{l-1})}\cup T^{[l-1]}]\in{\cal{LR}}^{(n+1)},$$ where ${T^{[l-1]}}$ is a ballot tableau on the alphabet $[l-1]$, and
 $\widehat T\in YT((\lambda_l,\dots,\lambda_{n+1})/(\mu_l))$ consists of the last $n-l+2$ rows of $T$ whose word  restricted to the alphabet $[l,n+1]$ satisfies the Yamanouchi condition.
Then, by Lemma \ref{lem:2-0}, the switching procedure   gives
\begin{equation}Y_{(\mu_l)}\cup \widehat T\underset{s}\rightarrow S\cup Q,\label{eq:induction}\end{equation}
  with $S\equiv \widehat T$ (the rectification of $\widehat T$)  in the alphabet $\{1,\dots,n+1\}$,  and $Q=D\ast Q^{[n-l+1]}\equiv Y_{(\mu_l)}$ in the alphabet $\{ l\}$.  If $|D|>0$ the last row of $S$ is the word $\widehat F(n+1)^{\nu_{n+1}}$ with $\widehat F$ a subword of $F$. Observing that $T=\widehat T\ast  T^{[l-1]}\equiv S\ast  T^{[l-1]}$ is a ballot tableau, we have
in addition
$$Y_{(\mu_1,\dots,\mu_l)}\cup T\underset{s}\rightarrow [Y_{(\mu_1,\dots,\mu_{l-1})}\cup (S\ast{T^{[l-1])}})]\cup(Q\ast\emptyset^{l-1}),$$
where $ Y_{(\mu_1,\dots,\mu_{l-1})}\cup (S\ast{T^{[l-1]}})$ is a ballot tableau pair. If $|D|>0$, \eqref{eq:switch0} and  \eqref{eq:x1recur} holds with $d=l$. Otherwise, $|D|=0$ and henceforth
the  $(n+1)$th row of $S$ \eqref{eq:induction} is the  $(n+1)$th row of $T$ with  length $\lambda_{n+1}$. Thereby
$$Y_{(\mu_1,\dots,\mu_l)}\cup T\underset{s}\rightarrow [Y_{(\mu_1,\dots,\mu_{l-1})}\cup (S\ast{T^{[l-1]}} )]\cup (\emptyset\ast Q^{[n-l+1]}\ast\emptyset^{l-1}),$$
 where $Y_{(\mu_1,\dots,\mu_{l-1})}\cup (S\ast{T^{[l-1]}})\in {\cal{LR}}^{(n+1)}$ is in the case $l-1$, and,  by inductive hypothesis, we get \eqref{eq:switch0}  and \eqref{eq:x1recur} with $1\le d\le l-1$.
\end{proof}


\begin{ex} $(I)$. $\ell(\mu)=1$.
 $n=4$, $\mu_1=7$,  $v=3$, $d=1$,
$$\YT{0.15in}{}{
 {\color{red}1,\color{red}1,{\color{red}1}|,\color{red}1,\color{red}1,\color{red}1,\color{red}1,1,1,1,1,1,1},
 {1,1,1|,1,1,2,2,2,2},
 {2,2,2|,2,3,3,3},
 {3,3,3|,4,4,4},
 {4,4,4|,5,5},
 }
 \underset {s}\rightarrow \YT{0.15in}{}{
 {\color{red}1,\color{red}1,{\color{red}1}|,\color{red}1,\color{red}1,1,1,1,1,1,1,\color{red}1,\color{red}1},
 {1,1,1|,1,1,2,2,2,2},
 {2,2,2|,2,3,3,3},
 {3,3,3|,4,4,4},
 {4,4,4|,5,5},
 }$$
$$
 \underset{s}\rightarrow
  \YT{0.15in}{}{
 {\color{red}1,\color{red}1,{\color{red}1}|,1,1,1,1,1,1,1,1,\color{red}1,\color{red}1},
 {1,1,1|,\color{red}1,2,2,2,2,\color{red}1},
 {2,2,2|,2,3,3,3},
 {3,3,3|,4,4,4},
 {4,4,4|,5,5},
 }
 \underset{s}\rightarrow
 \YT{0.15in}{}{
 {\color{red}1,\color{red}1,{\color{red}1}|,1,1,1,1,1,1,1,1,\color{red}1,\color{red}1},
 {1,1,1|,2,2,2,2,2,\color{red}1},
 {2,2,2|,\color{red}1,3,3,3},
 {3,3,3|,4,4,4},
 {4,4,4|,5,5},
 }
 $$
 $$
 \underset{s}\rightarrow
 \YT{0.15in}{}{
 {\color{red}1,\color{red}1,{\color{red}1}|,1,1,1,1,1,1,1,1,\color{red}1,\color{red}1},
 {1,1,1|,2,2,2,2,2,\color{red}1},
 {2,2,2|,3,3,3,\color{red}1},
 {3,3,3|,4,4,4},
 {4,4,4|,5,5},
 }
  \underset{s}\rightarrow\YT{0.15in}{}{
 {1,1,{1}|,1,1,1,1,1,1,1,1,\color{red}1,\color{red}1},
 {2,2,2|,2,2,2,2,2,\color{red}1},
 {3,3,3|,3,3,3,\color{red}1},
 {4,4,4|,4,4,4},
 {\color{red}1,\color{red}1,\color{red}1|,5,5},
 }$$
 $$
 \underset{s}\rightarrow\YT{0.15in}{}{
 {1,1,1|,1,1,1,1,1,1,1,1,\color{red}1,\color{red}1},
 {2,2,2|,2,2,2,2,2,\color{red}1},
 {3,3,3|,3,3,3,\color{red}1},
 {4,4,4|,4,4,4},
 {5,5,\color{red}1,\color{red}1,\color{red}1},
 }
 ,$$
$F=4^3$, $D=1^3$ and $Q=1^3\ast\emptyset\ast 1\ast 1\ast
 1^2$ of skew shape $\lambda/(11,8,6,6,2)$.

$(II)$. $l(\mu)=2$,
 $n=3$, $v= 4$, $d=2$
$$\YT{0.15in}{}{
 {\color{red}1,\color{red}1,\color{red}1,\color{red}1,\color{red}1,\color{red}1,1,1,1,1,1},
 {\color{red}2,\color{red}2,\color{red}2,\color{red}2,1,1,2,2,2},
 {1,2,2,2,2,3},
 {3,3,3,3},
 }
 \underset {s}\rightarrow
 \YT{0.15in}{}{
 {\color{red}1,\color{red}1,\color{red}1,\color{red}1,\color{red}1,\color{red}1,1,1,1,1,1},
 {\color{red}2,1,1,\color{red}2,\color{red}2,2,2,2,\color{red}2},
 {1,2,2,2,2,3},
 {3,3,3,3},
 }
 \underset {s}\rightarrow
 \YT{0.15in}{}{
 {\color{red}1,\color{red}1,\color{red}1,\color{red}1,\color{red}1,\color{red}1,1,1,1,1,1},
 {1,1,1,2,2,2,2,2,\color{red}2},
 {\color{red}2,2,2,\color{red}2,\color{red}2,3},
 {3,3,3,3},
 }
 $$
$$ \underset {s}\rightarrow
 \YT{0.15in}{}{
 {\color{red}1,\color{red}1,\color{red}1,\color{red}1,\color{red}1,\color{red}1,1,1,1,1,1},
 {1,1,1,2,2,2,2,2,\color{red}2},
 {2,2,\color{red}2,\color{red}2,\color{red}2,3},
 {3,3,3,3},
 }
\underset {s}\rightarrow
 \YT{0.15in}{}{
 {\color{red}1,\color{red}1,\color{red}1,\color{red}1,\color{red}1,\color{red}1,1,1,1,1,1},
 {1,1,1,2,2,2,2,2,\color{red}2},
 {2,2,3,3,3,\color{red}2},
 {3,3,\color{red}2,\color{red}2},
 },
 $$
 $F=3^4$, $\nu_4=0$, $\widehat F=3^2$, $D=2^2$ and $Q=2^2\ast 1\ast 1\ast \emptyset$ of skew shape $\lambda/(11,8,5,2)$.
\end{ex}

We first recall the following property.
{ The major result of the next statement is that the reading word $\widehat F(n+1)^{\nu_{n+1}} GCn^{\hat\nu_n}$ \eqref{cor1} of the two last rows of  $Y_{(\mu_1,\dots,\mu_{d-1})}\cup S$, and the concatenation of the $(n+1)$th row $F(n+1)^{\nu_{n+1}}$ \eqref{detachch} of $Y_\mu\cup T$ with the $n$th row  $\widehat GCn^{\hat\nu_n}$ \eqref{detachch} of     $Y_{(\mu_1,\dots,\mu_{d-1})}\cup R$,  that is, $F(n+1)^{\nu_{n+1}}GCn^{\hat\nu_n}$,  are related through  Knuth transformations. Since  $GC$ is  a row word  on the alphabet $[n-1]$,  $\widehat G$ a subword of $G$ and $F$ is a word in the alphabet $[n]$ and $\widehat F$ a subword of $F$,     from Lemma \ref{lem:congr}, $(a)$, \eqref{ls}, it means that $\widehat F GCn^{\hat\nu_n}\equiv F\widehat GCn^{\hat\nu_n}$ and from Lemma \ref{lem:congr}, $(c)$, 
 $ \widehat F GC\equiv F\widehat GC$ or $\widehat F G\equiv F\widehat G$.}

\begin{cor} \label{corol}Let $n\ge 1$, $Y_\mu\cup T\in {\cal{LR}}^{(n+1)}$ and  assume the assumptions of previous theorem on $Y_\mu\cup T$. Consider the equality \eqref{eq:switch0}, for some $1\le d\le n$ and $ \mu_d>0$, where we detach the two last rows of $S$ and $Q$,
\begin{eqnarray}Y_\mu\cup T&=F(n+1)^{\nu_{n+1}}\ast (Y\cup T)^{[n]}\underset{s}\rightarrow [Y_{(\mu_1,\dots,\mu_{d-1})}\cup S]\cup Q \label{detach-}\qquad\qquad\qquad\quad\\
&=[\widehat F(n+1)^{\nu_{n+1}}\ast GCn^{\hat\nu_n}\ast (Y_{(\mu_1,\dots,\mu_{d-1})}\cup S^{[n-1]})]\cup [D\ast X \ast Q^{[n-1]}],\label{cor1}
 \end{eqnarray}
with $Y_{(\mu_1,\dots,\mu_{d-1})}\cup S^{[n]}=GC n^{\hat\nu_n}\ast (Y_{(\mu_1,\dots,\mu_{d-1})}\cup S^{[n-1]})$, $Q^{[n]}=X \ast Q^{[n-1]}$
such that $X$ is a row word on the alphabet $[d,n]$, and $GC$, a row word  on the alphabet $[n-1]$,
with the factor $G$  satisfying $|G|=|F|=|\widehat F|+|D|$.
 Then \begin{eqnarray}(Y_\mu\cup T)^{[n]}&\underset{s}\rightarrow &[Y_{(\mu_1,\dots,\mu_{d-1})}\cup R]\cup P,\label{detachch}
\end{eqnarray}
where $P=  DX\ast  Q^{[n-1]}\equiv Q$, and $R\equiv T^{[n]}$ is such that  $Y_{(\mu_1,\dots,\mu_{d-1})}\cup R=\widehat GCn^{\hat\nu_n}\ast
(Y_{(\mu_1,\dots,\mu_{d-1})}\cup S^{[n-1]}),$
with $\widehat G$ a row subword of $G$   so that $|\widehat G|=|\widehat F|$, ($|\widehat G|+|D|=|F|=|G|$) and
$\widehat FG\equiv F\widehat G$ are Knuth equivalent. Also
\begin{align}\rho_1^{(n)}[(Y_\mu\cup T)^{[n]}]&=\rho_1^{(n)}(Y_{(\mu_1,\dots,\mu_{d-1})}\cup R)
\cup (DX\ast Q^{[n-1]})\nonumber\\
&=\rho_1^{(n)}[\widehat GCn^{\hat\nu_n}\ast
(Y_{(\mu_1,\dots,\mu_{d-1})}\cup S^{[n-1]})]\cup (DX\ast Q^{[n-1]}).\label{eq:cor}\end{align}
\end{cor}
\begin{proof}
Switching back in the two last rows of $Y_(\mu_1,\dots,\mu_{d-1})\cup S\cup Q$, \eqref{cor1}, and factoring $G=AB$ into two row words $A$ and $B$,  with $|A|=|\widehat F|$ and $|B|=|D|$,  one has
\begin{eqnarray}\label{switch1-}
[\widehat F(n+1)^{\nu_{n+1}}\cup D]\ast [GCn^{\hat\nu_n}\cup X]=\YT{0.40in}{}{
{A,B,C,\scriptstyle{n^e},\scriptstyle\cdots,\scriptstyle{n^h},\color{red}X},
{\widehat F,\scriptstyle{n+1^{|D|}},\scriptstyle\cdots,\scriptstyle{n+1^e},\color{red}{D}},
}
\overset{s}\longleftrightarrow\\
\nonumber\\ \label{switch2-}
\YT{0.45in}{}{
{G_1,F_2,G_3,F_4,\scriptstyle\cdots,G_{\scriptstyle k-1},F_k,C,\scriptstyle {n^e},\scriptstyle\cdots,\scriptstyle{n^h},\color{red}X},
{ F_1,\color{red}{D_2},F_3,\color{red}{D_4},\scriptstyle\dots,F_{\scriptstyle k-1},\color{red}{D_k},{\scriptstyle n+1^{|C|}},\scriptstyle\cdots,\scriptstyle{n+1^{|D|}}},
}
\overset{s}\longleftrightarrow
\end{eqnarray}
\begin{eqnarray}\label{switch3-}
\YT{0.45in}{}{
{G_1,\color{red}{D_2},G_3,\color{red}{D_4},\scriptstyle\cdots,G_{\scriptstyle k-1},\color{red}{D_k},C,\scriptstyle {n^e},\scriptstyle\cdots,\scriptstyle{n^h},\color{red}X},
{ F_1,F_2,F_3,F_4,\scriptstyle\dots,F_{\scriptstyle k-1},F_k,{\scriptstyle n+1^{|C|}},\scriptstyle\cdots,\scriptstyle{n+1^{|D|}}},
}
\end{eqnarray}
where ${\color{red} D}=\color{red}D_2D_4\cdots D_k$, $F=F_1F_2\cdots F_k$, $\widehat F=F_1F_3\cdots F_{k-1}$ and $\widehat G:=G_1G_3\cdots G_{k-1}$ a row subword of $G=AB=G_1F_2G_3\cdots G_{k-1}F_k$, with $|\widehat G|+|D|=|G|=|F|=|\widehat F|+|D|$. The subword of $F$ with  $\widehat F$ suppressed, $F\setminus \widehat F$, and $\widehat G$ are complementary row subwords of $G$.

  In every step of the {\em jeu de taquin} sliding  or {\em reverse jeu de taquin} sliding process the reading word is transformed into a Knuth equivalent one \cite{fulton, stanley}. Looking at the red letters as holes, the row reading words $\widehat FABC=\widehat FGC$ of \eqref{switch1-} (or \eqref{switch2-}), and $F \widehat G C$ of \eqref{switch3-}, restricted to the alphabet $[n]$ (in black), are Knuth equivalent.
 $$\widehat FGC\equiv F \widehat GC\Leftrightarrow \widehat FG\equiv F \widehat G.$$
The last row of \eqref{switch3-} is $F(n+1)^{\nu_{n+1}}$. Let $Y':=Y_{(\mu_1,\dots,\mu_{d-1})}$. It  then follows from \eqref{detach-},
\begin{eqnarray}
(Y_\mu\cup T)^{[n]}&\underset{s}\rightarrow&
\YT{0.21in}{}{
{G_1,\color{red}{D_2},G_3,\color{red}{D_4},\scriptstyle\cdots,\color{red}{D_k},C,\scriptstyle {n},\scriptstyle\cdots,\scriptstyle{n},\color{red}X}}
\ast [Y'
\cup S^{[n-1]}\cup Q^{[n-1])}]\nonumber\\
&\underset{s}\rightarrow&
\YT{0.21in}{}{
{\widehat G,C,\scriptstyle {n},\scriptstyle\cdots,\scriptstyle{n},
\color{red}D,\color{red}X},
}
\ast [Y'
\cup S^{[n-1]}\cup Q^{[n-1]}]\nonumber\\
&=&Y'\cup R\cup P=
[\widehat GCn^{\hat\nu_n}\cup {D}X]
\ast [Y'\cup S^{[n-1]}\cup Q^{[n-1]}],\label{eq:rows3}\end{eqnarray}
where $Y'\cup R=\widehat GC n^{\hat\nu_n}\ast(Y'\cup S^{[n-1]})$, and  $P$ is $Q^{[n]}$ with $D$ attached to the left of its  $n$th row $X$.
\end{proof}

\section{Proof of the Main Theorem 
}\label{sec:proof2}
Let $u=u_v\cdots u_1$ be an internal inserting order word of a tableau  $T$.
By $\bar\phi_u$-{\em bumping routes} on $Y\cup T$, we mean the collection of $\bar\phi_{u_i}$-{\em bumping routes} on $\bar\phi_{u_{i-1}\cdots u_1}(Y\cup T)$ for $i=1,\dots,v$.

\begin{lem} Let $u=u_v\cdots u_2u_1$ be a row word to be an internal   insertion order word of $T$. Then

$(a)$ the plactic class of $u$ is reduced to the sole  $u$.

$(b)$ the $\bar\phi_u$-bumping routes on $T$ are pairwise disjoint.

$(c)$   if the  $\bar\phi_{u_i}$-bumping route lands in  row $1\le k\le n+1$, the $\bar\phi_{u_{i+1}}$-bumping route lands strictly to the right in a row $\le k$, for $i=1,\dots,n-1$.
\end{lem}
\begin{proof} $(a)$ The plactic class of a row tableau has a sole element. $(b)$ and $(c)$  follow from Lemma \ref{lem:1} $(a)$.
\end{proof}

If $T\in YT(\lambda/\mu)$ is a ballot tableau and $\ell(\lambda)\le n$, and $1\le k\le \lambda_n-\mu_n$, then $T_{\underset{k}\leftarrow}$ denotes the ballot tableau in $YT((\lambda_1,\dots,\lambda_{n-1}, \lambda_n-k)/\mu)$  obtained from $T$ first by suppressing, in the $n$th row, the first $k$ filled boxes  and then pushing the remaining $\lambda_n-\mu_n-k$ boxes $k$ steps to the left.

$$ n=3, \mu=(6,4,1), \lambda=(11,9,6), \lambda_3-\mu_3=6, k=3<6,$$

\begin{align}Y_\mu\cup T=\YT{0.15in}{}{
 {\color{red}1,\color{red}1,\color{red}1,\color{red}1,\color{red}1,\color{red}1,1,1,1,1,1},
 {\color{red}2,\color{red}2,\color{red}2,\color{red}2,1,1,2,2,2},
 {\color{red}2,1,2,2,2,2,3},
 }& Y_\mu\cup T_{\underset 3\leftarrow}=\YT{0.15in}{}{
 {\color{red}1,\color{red}1,\color{red}1,\color{red}1,\color{red}1,\color{red}1,1,1,1,1,1},
 {\color{red}2,\color{red}2,\color{red}2,\color{red}2,1,1,2,2,2},
 {\color{red}2,2,2,3},
 }
 \end{align}
We   are now ready to prove the main result.

\bigskip

\noindent{\em Proof of Theorem \ref{2}.
}

 We may reduce the statement to the case $\mu=(\mu_1,\dots,$ $\mu_{n-1},$ $\mu_n=0)$. Let $\hat\mu:=\mu-(\mu_n^n)$. Performing horizontal switches in $Y\cup T=Y_{(\mu_n^n)}\cup Y_{\hat\mu}\cup T\underset{s}\rightarrow Y_{\hat\mu}\cup T\cup Z$, where $Z$ is the unique ballot tableau of shape $\lambda/(\lambda-(\mu_n^n))$ and content $(\mu_n^n)$. Thus
\begin{equation}\rho_1^{(n)}(Y\cup T)=\rho_1^{(n)}(Y_{\hat\mu}\cup T)\cup Z=\bar\chi_n^{\mu_n}[\rho_1^{(n)}(Y_{\hat\mu}\cup T)\cup (\emptyset\ast Z^{[n-1]})].\label{eqn:a}\end{equation}
 Similarly, $(Y\cup T)^{[n-1]}\underset{s}\rightarrow(Y_{\hat\mu}\cup  T)^{[n-1]}\cup Z^{[n-1]}$. Hence,
$\rho_1^{(n-1)}[(Y\cup T)^{[n-1]}]=\rho_1^{(n-1)}[(Y_{\hat\mu}\cup T)^{[n-1]}]\cup Z^{[n-1]}. $
Assuming that \eqref{introd:mainrecursionx} has been proved in the case of $\mu_n=0$, and using \eqref{eqn:a},  we  then may write
\begin{eqnarray}
\rho_1^{(n)}(Y\cup T)&=&\rho_1^{(n)}(Y_{\hat\mu}\cup T)\cup Z=\bar\chi_n^{\mu_n}[\rho_1^{(n)}(Y_{\hat\mu}\cup T)\cup (\emptyset\ast Z^{[n-1]})]\nonumber\\
&=&\bar\chi_n^{\mu_n}\{\bar\phi_{V_n}\bar\omega_n^{\nu_n}\rho_1^{(n-1)}[(Y_{\hat\mu}\cup T)^{[n-1]}]\cup Z^{[n-1]}\}\nonumber\\
&=&\bar\chi_n^{\mu_n}\bar\phi_{V_n} \bar\omega_n^{\nu_n}\left(\rho_1^{(n-1)}[(Y_{\hat\mu}\cup T)^{[n-1]}]\cup Z^{[n-1]}\right)\label{eqn:1}\\
&=&\bar\chi_n^{\mu_n}\bar\phi_{V_n}\bar\omega_n^{\nu_n}\rho_1^{(n-1)}[(Y\cup T)^{[n-1]}].\label{eqn:2}\nonumber
\end{eqnarray}
 Observe   \eqref{eqn:1} just says that the bumping routes of $\bar\phi_{V_n}$, all of them landing in the $n$th row, will not change when $\rho_1^{(n-1)}[(Y_{\hat\mu}\cup T)^{[n-1]}]\in{\cal{LR}}^{(n-1)}$ is extended  with the tableau $ Z^{[n-1]}$. Each route of $\phi_V$ follows the available path in $\rho_1^{(n)}[(\bar Y\cup T)^{[n-1]} ]$ and remains there  until landing in the $n$th row.

Let $\mu_n=0$. We now show, by induction on $n\ge 1$ and   $|V_n|\ge 0$,  that
\begin{equation}\label{eqn:recursion}
\rho_1^{(n)}(Y\cup T)=\bar\phi_{V_n}\bar\omega^{\widehat\nu_n}\rho_1^{(n-1)}[(Y\cup T)^{[n-1]}]=\bar\omega^{\widehat\nu_n}\bar\phi_{V_n}\rho_1^{(n-1)}[(Y\cup T)^{[n-1]}],
\end{equation}
where $Y=Y(\mu_1,\dots,\mu_{n-1},0)$ and $V_n$ is the word of the $n$th row of $T$ restricted to the alphabet $[n-1]$. In addition, all   bumping routes of $\bar\phi_{V_n}$  terminate in the $n$th row.

If $n=1$ then $|V_1|=0$, $\mu_1=0$, $T=Y_{({\nu_1})}$,  $\widehat\nu_1=\nu_1$, and  $\rho_1^{(1)}(\emptyset\cup Y_{({\nu_1}))}=Y_{({\nu_1})}\cup \emptyset$.
Then $\rho_1^{(1)}(\emptyset\cup Y_{(\nu_1)})=Y_{({\nu_1})}\cup \emptyset$ is obtained from $ \emptyset$ by adding  the row $1^{\nu_1}$.  Thus
 $\rho_1^{(1)}(\emptyset\cup Y_{({\nu_1}))}$ $=$ $\bar\omega_1^{\nu_1}\emptyset=\bar\omega_1^{\nu_1}\rho_1^{(0)}[(\emptyset\cup Y_{(\nu_1)})^{[0]}],$ with
$\bar\phi_{V_1}=id$.

Suppose that \eqref{eqn:recursion} holds for $n\ge 1$, and let us prove for $n+1$. Assume $Y\cup T$ with $n+1$ rows where $Y=Y(\mu_1,\dots,\mu_{n},0)$.

Let $v:=|V_{n+1}|\ge 0$ and  let $F:=V_{n+1}$ be the $(n+1)$th row word of $T$ restricted to the alphabet $[n]$.  Since $\mu_{n+1}=0$, then $F (n+1)^{\nu_{n+1}}$ is the  $(n+1)$th row of  $Y\cup T$ and detaching the $(n+1)$th row,
$Y\cup T=
F(n+1)^{\nu_{n+1}}
\ast (Y\cup T)^{[n]}.
$
We want to show that
\begin{equation}\label{goal}\rho_ 1^{(n+1)}(Y\cup T)=\bar\omega_{n+1}^{\nu_{n+1}}\bar\phi_F\rho_1^{(n)}[(Y\cup T)^{[n])}].\end{equation}
If $v=0$,   $Y\cup T=
(n+1)^{\nu_{n+1}}
\ast (Y\cup T)^{[n]}$.
 Therefore, by Lemma \ref{lem:2-1} 
 $$\rho_1^{(n+1)}(Y\cup T)=
(n+1)^{\nu_{n+1}}\ast \rho_1^{(n)}[(Y\cup T)^{[n]}]=\bar\omega_{n+1}^{\nu_{n+1}}\rho_1^{(n)}[(Y\cup T)^{[n]}],$$
  with $\bar\phi_F$ the identity.

If $v\ge 1$,  the $(n+1)$th  row of $Y\cup T$  is the word $F(n+1)^{\nu_{n+1}}$ with $F$ a no empty word on the alphabet $[n]$. We shall now use induction on $v$.

{\bf Step 1.}\textit{We pass from the ballot pair $Y\cup T$ to a ballot pair $Y'\cup S$ with $(n+1)$th row word  $\hat F(n+1)^{\nu_{n+1}}$ so that $\hat F$ is a strict subword of $F$.}

From Theorem \ref{recursion}, there exists $Y':=Y(\mu_1,\dots,\mu_{d-1})$, for some $1\le d\le n$, with $\mu_d> 0$, such that
\begin{align}\label{eq:switch}Y\cup T&=
F(n+1)^{\nu_{n+1}}
\ast (Y\cup T)^{[n]}\underset{s}\rightarrow Y'\cup S\cup Q\qquad\qquad\qquad\qquad\qquad\qquad\nonumber\\
&=
[\widehat F(n+1)^{\nu_{n+1}}\cup D]
\ast [(Y'\cup S)^{[n]}\cup Q^{[n]}],
\end{align}
 where $S=\widehat F(n+1)^{\nu_{n+1}}
\ast  S^{(n)}\equiv T$, $\widehat F$ is a strict subword of $F$, and  $Q={ D}\ast Q^{[n]}\equiv Y(\mu_d,\dots,\mu_n)$ is  over the alphabet $\{{d},\dots,{n}\}$ and  has  the $(n+1)$th row $D=d^{|D|}$ such that  $|D|=|F|-|\widehat F|>0$.
Therefore
\begin{equation}\label{eq:x1}
\rho_ 1^{(n+1)}(Y\cup T)=\rho_1^{(n+1)}(Y'\cup S)\cup Q.
\end{equation}
Since $\widehat F$ is the $(n+1)$th row of $Y'\cup S$, restricted to the entries $\le n$, and $0\le|\widehat F|<v$,
by induction on $v$, we may write
\begin{equation}\label{eq:x2}
\rho_1^{(n+1)}(Y'\cup S)=\bar\omega_{n+1}^{\nu_{n+1}}\bar\phi_{\widehat F}\rho_1^{(n)}[(Y'\cup S)^{[n]}],
\end{equation}
where all $\bar\phi_{\widehat F}$-{\em bumping routes}  terminate in the $(n+1)$th row.

To reach \eqref{goal}, one has, so far, from \eqref{eq:x1} and \eqref{eq:x2},
\begin{equation}\label{eq:x2*}
\rho_1^{(n+1)}(Y\cup T)=\bar\omega_{n+1}^{\nu_{n+1}}\bar\phi_{\widehat F}\rho_1^{(n)}[(Y'\cup S)^{[n]}]\cup (D\ast Q^{[n]}),
\end{equation}
with $\hat F$ a strict subword of $F$ such that $|\widehat F|+|D|=|F|$.

 {\bf Step 2.}  {\em One has to relate the $n$th row of $(Y'\cup S)^{[n]}$ with the $n$th row of  $(Y'\cup T)^{[n]}$.}

This requires Corollary \ref{corol} and, in particular, the analysis of the reading words in \eqref{switch1-}, \eqref{switch2-} and \eqref{switch3-}.

{\bf Step 2.1.}  \textit{We  first analyse $\rho_1^{(n)}[(Y'\cup S)^{[n]}]$ of \eqref{eq:x2*}.}

From Corollary \ref{corol} one has $Y'\cup S^{[n]}= GCn^{\hat\nu_{n}}\ast(Y'\cup S^{[n-1]})$ with  $G$ and  $C$  row words in the alphabet $[n-1]$, such that $|G|=|\widehat F|+|D|=v=|F|$.
 In addition $Y'\cup (S^{[n]})_{{\underset{v}\leftarrow}}=Cn^{\hat\nu_{n}}\ast(Y'\cup S^{[n-1]})\in {\cal{LR}}^{(n)}$ with $C$  a row word in the alphabet $[n-1]$.
By induction on $n$,   one has  \begin{eqnarray*}\rho_1^{(n)}(Y'\cup (S^{[n]})_{{\underset{v}\leftarrow}})=\bar\phi_{C}\omega_{n}^{\hat\nu_n}\rho_1^{(n-1)}(Y'\cup S^{[n-1]}),
\end{eqnarray*} and
\begin{align}
\rho_1^{(n)}(Y'\cup S^{[n]})&=\rho_1^{(n)}[GCn^{\hat\nu_{n}}\ast(Y'\cup S^{[n-1]})]\nonumber\\
&=\bar\phi_{G}\bar\phi_{C}\bar\omega_{n}^{\hat\nu_n}\rho_1^{(n-1)}(Y'\cup S^{[n-1]})\nonumber \\
&=\bar\phi_{G}\rho_1^{(n)}(Y'\cup (S^{[n]})_{{\underset{v}\leftarrow}}), \label{eq:x2v}
\end{align}
\noindent where  all $\bar\phi_G$-{\em bumping routes} (also $\bar\phi_C$-{\em bumping routes}) will end up in the $n$th row.

{\bf Step 2.2.} \textit{We now analyse $\rho_1^{(n)}(Y\cup T^{[n]})$.}

From Corollary \ref{corol}, one has
\begin{align}
Y\cup T^{[n]}\underset{s}\rightarrow
Y'\cup R\cup P=
[\widehat G C{n}^{\hat\nu_n}\cup {\color{red}{D}\color{red}X}]\ast [Y'\cup S^{[n-1]}\cup Q^{[n-1]}],
\nonumber\end{align}

where  $Y'\cup R=\widehat GC n^{\hat\nu_n}\ast(Y'\cup S^{[n-1]})$ with $\widehat G$ a subword of $G$, and  $P=DX\ast Q^{[n-1]}$  such that $|\widehat G|=|F|-|D|=v-|D|$ and
$\widehat FG\equiv F\widehat G$ Knuth equivalent.


 Observe that \begin{equation}\label{eq}Y'\cup (S^{[n]})_{\underset{v}\leftarrow}=Y'\cup R_{\underset{\scriptstyle v-|D|}\longleftarrow}=C n^{\hat\nu_n}\ast(Y'\cup S^{[n-1]}).\end{equation}
Again by induction on $n$,
  and using the identity \eqref{eq}, one has
\begin{align}
\rho_1^{(n)}(Y'\cup R_{\underset{\scriptstyle v-|D|}\longleftarrow})&=\bar\phi_{C}\bar\omega_{n}^{\tilde\nu_n}\rho_1^{(n-1)}(Y'\cup S^{[n-1]})\nonumber\\
\rho_1^{(n)}(Y'\cup R)&= \rho_1^{(n-1)}[\widehat GC n^{\nu_n}\ast(Y'\cup S^{[n-1]})]\nonumber\\
&=\bar\phi_{\widehat G}\bar\phi_{C}\bar\omega_{n}^{\hat\nu_n}\rho_1^{(n-1)}(Y'\cup S^{[n-1]})\nonumber\\
&=\bar\phi_{\widehat G}\rho_1^{(n)}(Y'\cup R_{\underset{\scriptstyle v-|D|}\longleftarrow}), \label{eq:x3v}
\end{align}
where all $\bar\phi_{\widehat G}$-{\em bumping routes } will end up in the $n$th row. Therefore from Corollary \ref{corol}, \eqref{eq:cor}, and \eqref{eq:x3v},
\begin{align}
 \rho_1^{(n)}(Y\cup T^{[n]})&=\rho_1^{(n)}(Y'\cup R)\cup P=\bar\phi_{\widehat G}\rho_1^{(n)}(Y'\cup R_{\underset{\scriptstyle v-|D|}\longleftarrow})\cup P\nonumber\\
&=\bar\phi_{\widehat G}\rho_1^{(n)}(Y'\cup R_{\underset{\scriptstyle v-|D|}\longleftarrow})\cup (DX\ast Q^{[n-1]}).\label{eq:final0}
\end{align}

 We are now in conditions to go back to  \eqref{eq:x2*}.

 {\bf Step 4.} \textit{Going back to \eqref{eq:x2*}}.
\begin{eqnarray}\label{eq:final}
\rho_1^{(n+1)}(Y\cup T)&=&\bar\omega_{n+1}^{\nu_{n+1}}\bar\phi_{\widehat F}\rho_1^{(n)}(Y'\cup S^{[n]})\cup ({\color{red} D\ast Q^{[n]}}),\,\text{using \eqref{eq:x2*}},\nonumber\\
&=&\bar\omega_{n+1}^{\nu_{n+1}}\bar\phi_{\widehat F}[\bar\phi_{G}\rho_1^{(n)}(Y'\cup (S^{[n]})_{{\underset{v}\leftarrow}})]\cup
 ({\color{red} D\ast Q^{[n]}}), \,\text{using \eqref{eq:x2v}},\nonumber\\
&=&\Big(\bar\omega_{n+1}^{\nu_{n+1}}\bar\phi_{\widehat F}\bar\phi_{G}\rho_1^{(n)}(Y'\cup R_{\underset{\scriptstyle v-|D|}\longleftarrow})\Big)\cup ({\color{red} D\ast Q^{[n]}}),\;\text{using \eqref{eq}},\quad\qquad\label{eq:final1}
\end{eqnarray}
with   all $\bar\phi_G$-{\em bumping routes}
ending up in the $n$th row, and the $\bar\phi_{\hat F}$-{\em bumping routes}  in the $(n+1)$th row. The cardinality of the $\bar\phi_G$-{\em bumping routes} is  $v=v-|D|+ |D|=$ $|\widehat G|+|D|=$ $|\widehat F|+|D|$.

{\bf Key Step 5.} \textit{Knuth relations of internal insertion operators}.

From Proposition \ref{propp:knuth}, since $\widehat F G\equiv F\widehat G$, one has

\begin{equation}\label{key}\bar\phi_{\widehat F G}\rho_1^{(n)}(Y'\cup R_{\underset{\scriptstyle v-|D|}\longleftarrow})=\bar\phi_{F\widehat  G}\rho_1^{(n)}(Y'\cup R_{\underset{\scriptstyle v-|D|}\longleftarrow}).\end{equation}

Another key fact and a simple observation is that  the number $v$ of bumping routes landing in the $n$th row and the number $v-|D|$ of bumping routes landing in the $(n+1)$th row is the same for $\bar\phi_{F}\bar\phi_{\widehat G}$ and $\bar\phi_{\widehat F}\bar\phi_{ G}$ when acting on $\rho_1^{(n)}(Y'\cup R_{\underset{\scriptstyle v-|D|}\longleftarrow})$. Thereby,
\begin{align}
\eqref{eq:final}
&=\bar\omega_{n+1}^{\nu_{n+1}}\bar\phi_{F}\bar\phi_{\widehat G}\rho_1^{(n)}(Y'\cup R_{\underset{\scriptstyle v-\ell(D)}\longleftarrow})\cup ({\color{red} D\ast  Q^{[n]}}),\;\text{using \eqref{key}},\label{order1}\\
&=\bar\omega_{n+1}^{\nu_{n+1}}\bar\phi_{F}\rho_1^{(n)}(Y'\cup R)\cup ({\color{red} D\ast X\ast  Q^{[n-1]}}),\;\text{using \eqref{eq:x3v}}.\label{order2}
\end{align}
From \eqref{order1} to \eqref{order2}, $|\widehat G|=v-|D|$ bumping routes are  executed by $\bar\phi_{\widehat G}$ and land in row $n$. This implies that in the action of $\bar\phi_{F}$ over $\rho_1^{(n)}(Y'\cup R)$, {\em $|D|$ of the $v=|F|$ pairwise disjoint bumping routes will  still land and settle in the $n$th the row, which   means when settling  to adding $|D|$ new entries at the end of the $n$th row of $\rho_1^{(n)}(Y'\cup R)$}, while $v-|D|$ of them will land in the $(n+1)$th row.
  Recall that since $F$ is a row word, from Lemma \ref{lem:1}, $(a)$, the $v$ bumping routes are pairwise disjoint and, more importantly, the $|D|$ bumping routes settling in the $n$th row are necessarily the last to be executed.
This means  that if we attach $D$ at the end of the $n$th row $\rho_1^{(n)}(Y'\cup R)$, the rightmost $|D|$ bumping routes of $\bar\phi_F$ when landing  to the $n$th row will meet the entire row $\color{red} D$ and bumps it out
 to the $(n+1)$th row. Thus  recalling that $P=DX\ast Q^{[n-1]}$ and identity \eqref{eq:final0}, 
\begin{align}
\eqref{order2}&=\bar\omega_{n+1}^{\nu_{n+1}}\bar\phi_{F}[\rho_1^{(n)}(Y'\cup R)\cup (DX\ast Q^{[n-1]})]\\
&=\bar\omega_{n+1}^{\nu_{n+1}}\bar\phi_{F}[\rho_1^{(n)}(Y'\cup R)\cup P], 
\\
&=\bar\phi_{F}\bar\omega_{n+1}^{\nu_{n+1}}\rho_1^{(n)}(Y\cup T^{[n]}),~~\text{using  \eqref{eq:final0}.}\label{end1}
\end{align}

We have shown  $\rho_ 1^{(n+1)}(Y\cup T)=\eqref{end1}$, that is,  identity \eqref{goal},
$$\qquad\qquad\qquad\rho_ 1^{(n+1)}(Y\cup T)=\bar\omega_{n+1}^{\nu_{n+1}}\bar\phi_F\rho_1^{(n)}[(Y\cup T)^{[n]}].
\qquad\qquad\qquad\qquad\qquad\qquad\square$$


\bibliography{sample17}
\bibliographystyle{alpha}

\end{document}